# INVARIANCE PRINCIPLES FOR RANDOM BIPARTITE PLANAR MAPS

By Jean-François Marckert and Grégory Miermont

*CNRS, Université de Bordeaux 1 and Université de Paris-Sud*

Random planar maps are considered in the physics literature as the discrete counterpart of random surfaces. It is conjectured that properly rescaled random planar maps, when conditioned to have a large number of faces, should converge to a limiting surface whose law does not depend, up to scaling factors, on details of the class of maps that are sampled. Previous works on the topic, starting with Chassaing and Schaeffer, have shown that the radius of a random quadrangulation with $n$ faces, that is, the maximal graph distance on such a quadrangulation to a fixed reference point, converges in distribution once rescaled by $n^{1/4}$ to the diameter of the Brownian snake, up to a scaling constant.

Using a bijection due to Bouttier, Di Francesco and Guitter between bipartite planar maps and a family of labeled trees, we show the corresponding invariance principle for a class of random maps that follow a Boltzmann distribution putting weight $q_k$ on faces of degree $2k$: the radius of such maps, conditioned to have $n$ faces (or $n$ vertices) and under a criticality assumption, converges in distribution once rescaled by $n^{1/4}$ to a scaled version of the diameter of the Brownian snake. Convergence results for the so-called profile of maps are also provided. The convergence of rescaled bipartite maps to the Brownian map, in the sense introduced by Marckert and Mokkadem, is also shown. The proofs of these results rely on a new invariance principle for two-type spatial Galton–Watson trees.

## 1. Introduction, motivations and main results.

1.1. *Motivation.* An *embedded graph* $\mathcal{G}$ is an embedding of a connected graph in the two-dimensional sphere $\mathbb{S}^2$, in which edges do not intersect except possibly at their endpoints (the vertices). A *face* of $\mathcal{G}$ is a connected component of $\mathbb{S}^2 \setminus \mathcal{G}$. Faces are homeomorphic to open disks, and the *degree* of a given face is the number of edges that are included in the closure of this face, with the convention that cut edges are counted twice, where cut edges are those edges whose removal disconnects the graph. If the graph is the vertex graph with only one vertex and no edges, we adopt the convention that it bounds one face with degree 0. The degree of a vertex is the number of edges adjacent to that vertex, where self-loops are counted twice, according to the usual graph-theoretic definition. Unlike faces, it depends only on the underlying graph rather than its embedding in $\mathbb{S}^2$.

We say that two embedded graphs are equivalent if there exists an orientation-preserving homeomorphism of $\mathbb{S}^2$ that maps the first embedding to the second. Equivalence classes of embedded graphs are called *planar maps*, and their set is denoted by $\mathcal{M}_0$. When considering a planar map $\mathbf{m} \in \mathcal{M}_0$, we will slightly improperly speak of its vertices, edges, faces and their respective degrees (we should first take an element of the class $\mathbf{m}$ to be completely accurate). We let $S(\mathbf{m}), A(\mathbf{m}), F(\mathbf{m})$ be the sets of vertices, edges and faces of $\mathbf{m}$. The degree of an element $u \in S(\mathbf{m})$ or $f \in F(\mathbf{m})$ will be denoted by $\deg(u)$, respectively $\deg(f)$. We denote the class of the vertex graph by $\dagger$.

If $u, v$ are vertices in a planar map $\mathbf{m} \in \mathcal{M}_0$, and $e_1, \ldots, e_n$ are oriented edges, we say that $e_1, \ldots, e_n$ is a *path* from $u$ to $v$ of length $n$ if the source of $e_1$ is $u$, the target of $e_n$ is $v$, and the target of $e_i$ is the source of $e_{i+1}$ for all $1 \le i \le n-1$. The *graph distance* associated with a planar map $\mathbf{m} \in \mathcal{M}_0$ is the function $d_{\mathbf{m}} : S(\mathbf{m}) \times S(\mathbf{m}) \to \mathbb{Z}_+$ defined by letting $d_{\mathbf{m}}(u, v)$ be the least $n$ such that there exists a path of length $n$ leading from $u$ to $v$. This can be interpreted by saying that we turn $\mathbf{m}$ into a metric space, by endowing edges with lengths all equal to 1.

Planar maps have been of particular interest to physicists in the last decade as they can be considered as discretized versions of surfaces. In order to give a mathematical ground to the stochastic quantization of two-dimensional gravity, in which an integral with respect to an ill-defined uniform measure on Riemannian surfaces is involved, a possible attempt is to replace the integral by a finite sum over distinct discrete geometries, whose role is performed by planar maps [3]. Informally, it is believed that:

- A random map chosen in some class of planar maps with size $n$ (e.g., a quadrangulation with $n$ faces, i.e., a map whose $n$ faces are all of degree 4), whose edge lengths are properly rescaled, should converge in distribution as $n \to \infty$ to a limiting random surface.
- The limiting random surface should not depend, up to scale factors, on details of the class of maps which is randomly sampled.

The second property is called *universality*. A similar situation is well known to probabilists: the role of a Lebesgue measure on paths is performed by Brownian motion, which is the scaling limit of discretized random paths (random walks) whose step distributions have a finite variance.

In a pioneering work, Chassaing and Schaeffer [8] made a very substantial progress in answering the first question, by establishing that the largest distance to the root in a uniform rooted quadrangulation with $n$ faces (see definition below) divided by $n^{1/4}$ converges in distribution to some random variable (which is, up to a multiplicative constant, the diameter of the range of the so-called Brownian snake with lifetime process the normalized Brownian excursion). By using an invariance principle for discrete labeled trees satisfying a positivity constraint, Le Gall [16] has given an alternative proof of the results of [8]. This involves a new random object, called the Brownian snake conditioned to be positive, that was introduced in Le Gall and Weill [17]. Marckert and Mokkadem [20] gave a description of quadrangulations by gluing two trees, and showed that these trees converge when suitably normalized as $n$ goes to $\infty$. They introduced the notion of Brownian map, and showed that under a certain topology, rescaled quadrangulations converge in distribution to the Brownian map. All these results have been obtained by using bijective methods which take their source in the work of Schaeffer [23], and which allow to study random quadrangulations in terms of certain labeled trees. The nice feature of this method is that the labels allow to keep track of geodesic distances to a reference vertex in the map, so that some geometric information on the maps is present in the associated labeled trees.

On the other hand, the second question has not been addressed up to now in a purely probabilistic form, and in the context of scaling limits of planar maps. Angel [4] and Angel and Schramm [5] give evidence that the large-scale properties of large planar maps should not depend on the local details of the map (like the degree of faces), but these remarks hold in the context of *local limits* of random maps, where all edges have a length fixed to 1 as the number of faces of the map goes to infinity (this is an infinite volume limit), rather than in the context of scaling limits, where edge lengths tend to 0 as the number of faces goes to infinity (so that the total volume is kept finite). In a recent article, Bouttier, Di Francesco and Guitter [6] have given a generalization of Schaeffer's bijection to general planar maps. They obtain identities for the generating series of the most general family of (weighted) planar maps, and infer a number of clues for the universality of the pure 2D gravity model, for example, by computing certain scaling exponents with a combinatorial approach.

Their bijection suggests a path to prove invariance principles (the probabilistic word for universality) for random maps. The present work explores this path in the case of bipartite maps, by first giving a probabilistic interpretation of the identities of [6].

1.2. *Boltzmann laws on planar maps.* A planar map is said to be *bipartite* if all its faces have even degree. In this paper, we will only be concerned with bipartite maps, notice † is bipartite with our convention.

Every edge of a map can be given two orientations. A bipartite *rooted* planar map is a pair $(\mathbf{m}, e)$ where $\mathbf{m}$ is a bipartite map and $e$ is a distinguished oriented edge of $\mathbf{m}$. The basic objects that are considered in this article are bipartite planar maps which are rooted and pointed, that is, triples $(\mathbf{m}, e, \mathfrak{r})$ where $(\mathbf{m}, e)$ is a bipartite rooted planar map and $\mathfrak{r}$ is a vertex of $\mathbf{m}$. We let $\mathcal{M}$ be the set of rooted, pointed, bipartite planar maps. The map † cannot be rooted and can be pointed only at its unique vertex, but is still considered as an element of $\mathcal{M}$. By abuse of notation, we will often denote a generic element of $\mathcal{M}$ by $\mathbf{m}$ without referring to $(e, \mathfrak{r})$ when it is free of ambiguity.

By the bipartite nature of elements of $\mathcal{M}$, we have $|d_\mathbf{m}(\mathfrak{r}, u) - d_\mathbf{m}(\mathfrak{r}, v)| = 1$ whenever $u, v \in S(\mathbf{m})$ are neighbors. Therefore, if $(\mathbf{m}, e, \mathfrak{r}) \in \mathcal{M} \setminus \{\dagger\}$, we have either $d_\mathbf{m}(\mathfrak{r}, e_+) > d_\mathbf{m}(\mathfrak{r}, e_-)$ or $d_\mathbf{m}(\mathfrak{r}, e_+) < d_\mathbf{m}(\mathfrak{r}, e_-)$, where $e_-$ and $e_+$ are the source and the target of the oriented edge $e$. We let

$$\mathcal{M}_+ = \{(\mathbf{m}, e, \mathfrak{r}) \in \mathcal{M} : d_\mathbf{m}(\mathfrak{r}, e_+) > d_\mathbf{m}(\mathfrak{r}, e_-)\} \cup \{\dagger\}.$$

All probability distributions on maps in this paper are going to be defined on the set $\mathcal{M}_+$. Notice that an alternative definition for this set is to consider it as the set of pointed maps where a nonoriented edge has been distinguished.

Let $\mathbf{q} = (q_i, i \geq 1)$ be a sequence of nonnegative weights such that $q_i > 0$ for at least one $i > 1$. By convention, let $q_0 = 1$. Consider the $\sigma$-finite measure $W_\mathbf{q}$ on $\mathcal{M}_+$ that assigns to each map $\mathbf{m} \in \mathcal{M}_+$ a weight $q_i$ per face of degree $2i$:

$$\tag{1} W_\mathbf{q}(\mathbf{m}) = \prod_{f \in F(\mathbf{m})} q_{\deg(f)/2},$$

with the convention $W_\mathbf{q}(\dagger) = q_0 = 1$. This multiplicative form is reminiscent of the measures associated with the so-called simply generated trees, which are of the form $w(\mathbf{t}) = \prod_{u \in \mathbf{t}} q_{c_\mathbf{t}(u)}$ for any tree $\mathbf{t}$, where $c_\mathbf{t}(u)$ is the number of children of a vertex $u$ in $\mathbf{t}$, and where $(q_i, i \geq 0)$ is a sequence of nonnegative numbers (see [1], pages 27–28).

Let $Z_\mathbf{q} = W_\mathbf{q}(\mathcal{M}_+)$ be the "partition function of $\mathbf{q}$." Notice that $Z_\mathbf{q} \in (1, \infty]$ since $W_\mathbf{q}(\dagger) = 1$. If $Z_\mathbf{q} < \infty$, we say that $\mathbf{q}$ is *admissible*, and introduce the *Boltzmann* distribution on $\mathcal{M}_+$ with *susceptibility* $\mathbf{q}$ by letting

$$P_\mathbf{q} = \frac{W_\mathbf{q}}{Z_\mathbf{q}}.$$

For $k \geq 1$, let $N(k) = \binom{2k-1}{k-1}$. For any weight sequence $\mathbf{q}$ (not necessarily admissible) define

$$f_\mathbf{q}(x) = \sum_{k \geq 0} x^k N(k+1) q_{k+1} \in [0, \infty], \qquad x \geq 0.$$

The function $f_{\mathbf{q}}:[0,\infty) \to [0,\infty]$ is a completely positive power series, that is, its derivatives of every order are nonnegative, and since $(q_i, i > 1)$ is not identically zero, $f_{\mathbf{q}}$ is strictly positive on $(0,\infty)$, and strictly increasing on the interval $[0, R_{\mathbf{q}}]$, where $R_{\mathbf{q}} \in [0,\infty]$ is the radius of convergence of $f_{\mathbf{q}}$. Moreover, $f_{\mathbf{q}}(x)$ converges to $\infty$ as $x \to \infty$, and the monotone convergence theorem entails that the function $f_{\mathbf{q}}$ is continuous from $[0, R_{\mathbf{q}}]$ to $[0,\infty]$, that is, $f_{\mathbf{q}}(R_{\mathbf{q}}) = f_{\mathbf{q}}(R_{\mathbf{q}}-)$ [which can be finite or infinite, while by definition $f_{\mathbf{q}}(R_{\mathbf{q}}+) = \infty$]. In the sequel, we understand that $f'_{\mathbf{q}}(R_{\mathbf{q}}) \in (0,\infty]$ stands for the left derivative of $f_{\mathbf{q}}$ at $R_{\mathbf{q}}$ (when $R_{\mathbf{q}} > 0$).

Consider the equation

$$(2) \qquad f_{\mathbf{q}}(x) = 1 - 1/x, \qquad x > 0.$$

Since $x \mapsto 1 - x^{-1}$ is nonpositive on $(0,1]$ and $f_{\mathbf{q}}$ is infinite on $(R_{\mathbf{q}}, \infty]$, a solution of (2) always belongs to $(1, R_{\mathbf{q}}]$. Since $x \mapsto 1 - x^{-1}$ is strictly concave on $(0, +\infty)$, with derivative $x \mapsto x^{-2}$, and $f_{\mathbf{q}}$ is convex, strictly increasing and continuous on $[0, R_{\mathbf{q}}]$, we can classify the configurations of solutions for (2) by the following four exclusive cases:

1. There are no solutions.
2. There are exactly two solutions $z_1 < z_2$ in $(1, R_{\mathbf{q}}]$, in which case $f'_{\mathbf{q}}(z_1) < z_1^{-2}$ and $f'_{\mathbf{q}}(z_2) > z_2^{-2}$.
3. There is exactly one solution $z_1$ in $(1, R_{\mathbf{q}}]$ with $f'_{\mathbf{q}}(z_1) < z_1^{-2}$.
4. There is exactly one solution $z$ in $(1, R_{\mathbf{q}}]$ with $f'_{\mathbf{q}}(z) = z^{-2}$.

As will be shown in Section 2.3, the admissibility of $\mathbf{q}$ can be formulated in terms of $f_{\mathbf{q}}$ as follows.

PROPOSITION 1. *The weight sequence $\mathbf{q}$ is admissible if and only if equation (2) has at least one solution. In this case, $Z_{\mathbf{q}}$ is the solution of (2) that satisfies $Z_{\mathbf{q}}^2 f'_{\mathbf{q}}(Z_{\mathbf{q}}) \leq 1$.*

In this paper, we consider case 3 above, and cases when one of the solutions of (2) is equal to $R_{\mathbf{q}}$, as nonregular cases. Also, note that the case 4 in the above classification plays a singular role compared to the others. These remarks motivate the following.

DEFINITION 1. An admissible weight sequence $\mathbf{q}$ is said to be *critical* if case 4 of the above classification is satisfied, that is,

$$(3) \qquad Z_{\mathbf{q}}^2 f'_{\mathbf{q}}(Z_{\mathbf{q}}) = 1.$$

Equivalently, $\mathbf{q}$ is critical if and only if the graphs of $x \mapsto f_{\mathbf{q}}(x)$ and $x \mapsto 1 - 1/x$ are tangent to the left of $x = Z_{\mathbf{q}}$.

We say that $\mathbf{q}$ is *regular critical* if it is critical and $Z_{\mathbf{q}} < R_{\mathbf{q}}$, that is, the graphs are tangent at $Z_{\mathbf{q}}$ both to the left and to the right.

Notice that a critical weight sequence **q** is automatically regular in the case where $f_\mathbf{q}(R_\mathbf{q}) = \infty$; in this case, **q** is regular critical if and only if equation (2) admits a unique solution (because case 3 in the above classification cannot happen).

1.3. *Snakes.* In order to state our main theorem, we first briefly describe the limiting random objects that are involved. Let $B^{\mathrm{exc}}$ be a standard Brownian excursion. Then, given $B^{\mathrm{exc}}$, we let $S^{\mathrm{exc}}$ be a centered Gaussian process whose covariance function is given by

$$\mathrm{cov}(S_s^{\mathrm{exc}}, S_t^{\mathrm{exc}}) = \inf_{s \wedge t \leq u \leq s \vee t} B_u^{\mathrm{exc}}, \qquad 0 \leq s, t \leq 1. \tag{4}$$

It is known (see, e.g., [15], Section IV.6) that $(B^{\mathrm{exc}}, S^{\mathrm{exc}})$ has a continuous version, which is the one we choose to work with. We let $\mathbb{N}^{(1)}$ be the law of the pair $(B^{\mathrm{exc}}, S^{\mathrm{exc}})$. In the sequel, we will let $((\mathrm{e}_s)_{0 \leq s \leq 1}, (\mathrm{r}_s)_{0 \leq s \leq 1})$ be the canonical process for the space $\mathcal{C}(\mathbb{R}_+, \mathbb{R})^2$ of continuous functions with values in $\mathbb{R}^2$. The process $(\mathrm{e}, \mathrm{r})$ under $\mathbb{N}^{(1)}$ is called the "head of the Brownian snake" driven by a Brownian excursion in the literature [12, 13, 19]. We let

$$\Delta_+(\mathrm{r}) = \sup_{t \geq 0} \mathrm{r}_t, \qquad \Delta_-(\mathrm{r}) = \inf_{t \geq 0} \mathrm{r}_t \quad \text{and} \quad \Delta(\mathrm{r}) = \Delta_+(\mathrm{r}) - \Delta_-(\mathrm{r}),$$

the positive and negative range of r, and the diameter of the range of r.

1.4. *Main results.* For $(\mathbf{m}, e, \mathfrak{r}) \in \mathcal{M}_+$, let

$$\mathscr{R}(\mathbf{m}, e, \mathfrak{r}) = \max_{u \in S(\mathbf{m})} d_\mathbf{m}(\mathfrak{r}, u)$$

be the *radius* of $(\mathbf{m}, e, \mathfrak{r}) \in \mathcal{M}_+$. Also, let $\mathscr{I}^{(\mathbf{m}, e, \mathfrak{r})}$ be the *normalized profile* of the map **m**, which is the probability measure on $\mathbf{Z}_+$ such that

$$\mathscr{I}^{(\mathbf{m}, e, \mathfrak{r})}(k) = \frac{\#\{u \in S(\mathbf{m}) : d_\mathbf{m}(\mathfrak{r}, u) = k\}}{\#S(\mathbf{m})}, \qquad k \geq 0.$$

For simplicity we will usually denote these quantities by $\mathscr{R}(\mathbf{m}), \mathscr{I}^\mathbf{m}$. If $n \geq 1$, we also let $\mathscr{I}_n^\mathbf{m}$ be the rescaled measure on $\mathbb{R}_+$ which is defined by $\mathscr{I}_n^\mathbf{m}(A) = \mathscr{I}^\mathbf{m}(n^{1/4} A)$, for $A$ a Borel subset of $\mathbb{R}_+$.

Last, if **q** is a regular critical weight sequence, we let

$$\rho_\mathbf{q} = 2 + Z_\mathbf{q}^3 f_\mathbf{q}''(Z_\mathbf{q}). \tag{5}$$

Letting $M : \mathcal{M}_+ \to \mathcal{M}_+$ be the identity mapping, our main result states as follows.

THEOREM 2. *Let **q** be a regular critical weight sequence. Then:*

(i) *The distribution of the random variable $n^{-1/4}\mathscr{R}(M)$, under $P_\mathbf{q}(\cdot|\#F(M) = n)$, converges weakly as $n \to \infty$ to the law under $\mathbb{N}^{(1)}$ of*

$$\left(\frac{4\rho_\mathbf{q}}{9(Z_\mathbf{q} - 1)}\right)^{1/4} \Delta(\mathrm{r}).$$

(ii) *The distribution of $n^{-1/4} d_M(\mathfrak{r}, \mathfrak{r}')$ under $P_\mathbf{q}(\cdot|\#F(M) = n)$, where $\mathfrak{r}' \in S(M) \setminus \{\mathfrak{r}\}$ is picked uniformly at random conditionally on $M$, converges weakly as $n \to \infty$ to the law under $\mathbb{N}^{(1)}$ of*

$$\left(\frac{4\rho_\mathbf{q}}{9(Z_\mathbf{q} - 1)}\right)^{1/4} \Delta_+(\mathrm{r}).$$

(iii) *The distribution of the random measure $\mathscr{I}_n^M$ under $P_\mathbf{q}(\cdot|\#F(M) = n)$ converges weakly to the law under $\mathbb{N}^{(1)}$ of the random probability measure $\mathscr{I}_\mathbf{q}^\mathrm{r}$ on $\mathbb{R}_+$, defined by*

$$\langle \mathscr{I}_\mathbf{q}^\mathrm{r}, g \rangle = \int_0^1 ds\, g\left(\left(\frac{4\rho_\mathbf{q}}{9(Z_\mathbf{q} - 1)}\right)^{1/4}\left(\mathrm{r}_s - \inf_{0 \le u \le 1} \mathrm{r}_u\right)\right),$$

*for any continuous and bounded function $g : \mathbb{R} \to \mathbb{R}$.*

Notice that Boltzmann distributions always put a positive mass on the set of maps with exactly $n$ faces for all $n$, so that the conditional distributions $P(\cdot|\#F(M) = n)$ are well defined. There exists also a counterpart of this result in which we condition on the number of vertices rather than the number of faces, which states as

PROPOSITION 3. *Let $\mathbf{q}$ be a regular critical weight sequence. Then the previous theorem remains true when considering the laws $P_\mathbf{q}(\cdot|\#S(M) = n)$ instead of $P_\mathbf{q}(\cdot|\#F(M) = n)$, where it must be understood that $n \to \infty$ along values for which $P(\#S(M) = n) > 0$, and the rescaling constant $(4\rho_\mathbf{q}/(9(Z_\mathbf{q} - 1)))^{1/4}$ appearing in* (i), (ii), (iii) *must be replaced by $(4\rho_\mathbf{q}/9)^{1/4}$.*

The last two results are stated under the assumption of admissibility and regular criticality for the weight sequence. However, since the probability laws that appear in the statements are conditioned measures, and thus make a slightly indirect use of the probability $P_\mathbf{q}$, this assumption can be loosened a bit. If

$$\mathcal{M}_+^{F=n} = \{\mathbf{m} \in \mathcal{M}_+ : \#F(\mathbf{m}) = n\},$$

and $\mathbf{q}$ is any weight sequence, then the hypothesis $Z_\mathbf{q}^{F=n} = W_\mathbf{q}(\mathcal{M}_+^{F=n}) < \infty$ allows to define a probability measure on $\mathcal{M}_+^{F=n}$ by

$$P_\mathbf{q}^{F=n}(\cdot) = \frac{W_\mathbf{q}(\cdot \cap \mathcal{M}_+^{F=n})}{Z_\mathbf{q}^{F=n}}.$$

If **q** is admissible, we are clearly in this case, and $P_{\mathbf{q}}(\cdot|\#F(M)=n) = P_{\mathbf{q}}^{F=n}$, but the converse is not true: there can be (and there are in many interesting cases) weight sequences that are not admissible, but for which $P_{\mathbf{q}}^{F=n}$ makes sense. Now, notice that if $\alpha > 0$ and $\alpha \mathbf{q} = (\alpha q_i, i \geq 1)$, and $\mathbf{m} \in \mathcal{M}_+^{F=n}$,

$$W_{\alpha \mathbf{q}}(\{\mathbf{m}\}) = \alpha^n W_{\mathbf{q}}(\{\mathbf{m}\}), \qquad Z_{\alpha \mathbf{q}}^{F=n} = \alpha^n Z_{\mathbf{q}}^{F=n}.$$

Therefore, if **q** is such that $Z_{\mathbf{q}}^{F=n} < \infty$, then $\alpha \mathbf{q}$ is also such a weight sequence, and $P_{\alpha \mathbf{q}}^{F=n} = P_{\mathbf{q}}^{F=n}$ is independent of $\alpha > 0$.

In a similar way, we let $\mathcal{M}_+^{S=n}$ to be the set of maps with $n$ vertices, and define as above $Z_{\mathbf{q}}^{S=n}$, and $P_{\mathbf{q}}^{S=n}$ if the latter is finite and $> 0$. For $\beta > 0$ let $\beta \bullet \mathbf{q} = (\beta^{i-1} q_i, i \geq 1)$. Then for $\mathbf{m} \in \mathcal{M}_+^{S=n}$,

$$W_{\beta \bullet \mathbf{q}}(\{\mathbf{m}\}) = \prod_{f \in F(\mathbf{m})} \beta^{\deg(f)/2 - 1} q_{\deg(f)/2} = \beta^{n-2} W_{\mathbf{q}}(\{\mathbf{m}\}),$$

and $Z_{\beta \bullet \mathbf{q}}^{S=n} = \beta^{n-2} Z_{\mathbf{q}}^{S=n}$, where we used $\sum_{f \in F(\mathbf{m})} \deg(f) = 2\#A(\mathbf{m})$, and Euler's formula $\#F(\mathbf{m}) - \#A(\mathbf{m}) + \#S(\mathbf{m}) = 2$. Thus, Theorem 2 and Proposition 3 can be restated as follows.

COROLLARY 4. (i) *Let* **q** *be a weight sequence such that* $Z_{\mathbf{q}}^{F=n} < \infty$ *for every* $n \geq 1$ *and such that there exists some* $\alpha_c > 0$ *such that* $\alpha_c \mathbf{q}$ *is regular critical. Then all conclusions of Theorem 2 remain true, when replacing the probabilities* $P_{\mathbf{q}}(\cdot|\#F(M)=n)$ *in the statement by* $P_{\mathbf{q}}^{F=n}$, *and where the normalizing constant of* (i), (ii), (iii) *is computed for the weight sequence* $\alpha_c \mathbf{q}$.

(ii) *Let* **q** *be such that* $Z_{\mathbf{q}}^{S=n} < \infty$ *and there exists some* $\beta_c > 0$ *with* $\beta_c \bullet \mathbf{q}$ *regular critical, then the conclusion of Proposition 3 remain true, when considering* $P_{\mathbf{q}}^{S=n}$ *instead of* $P_{\mathbf{q}}(\cdot|\#S(M)=n)$, *and computing the scaling constants for the weight sequence* $\beta_c \bullet \mathbf{q}$.

It is also true that conditioning both on the number of faces and vertices is insensitive to termwise multiplication of **q** by $(\alpha \beta^{i-1}, i \geq 1)$, so this would lead to finding a critical curve of $(\alpha_c, \beta_c)$'s such that $(\alpha_c \beta_c^{i-1} q_i, i \geq 1)$ is critical. We do not concentrate on this last point, as our methods are inefficient in conditioning on both these data.

1.5. *Two illustrating examples.* We illustrate Theorem 2 by explicitly computing the various constants involved there in two natural particular cases.

1.5.1. *$2\kappa$-angulations.* Consider the case when $\mathbf{q} = \alpha \delta_\kappa$, for some integer $\kappa \geq 2$, and some constant $\alpha > 0$. The resulting distributions are the Boltzmann distributions on the set of maps with faces of fixed degree $2\kappa$. These distributions appear in [7] in the case $\kappa = 2$ of quadrangulations (they also appear in [4], but for triangulations).

In that case, $f_\mathbf{q}$ takes the simple form of a monomial $f_\mathbf{q}(x) = \alpha N(\kappa) x^{\kappa-1}$, which satisfies $R_\mathbf{q} = \infty$. According to Proposition 1, Definition 1 and the fact that $f_\mathbf{q}(R_\mathbf{q}) = \infty$, the weight sequence $\mathbf{q}$ is critical, and thus regular critical, if and only if and the system of equations

$$f_\mathbf{q}(z) = 1 - 1/z, \qquad z^2 f'_\mathbf{q}(z) = 1$$

has a real solution. This system, considered in the variables $\alpha, z$ admits the unique solution

$$\alpha_\kappa = \frac{(\kappa-1)^{\kappa-1}}{\kappa^\kappa N(\kappa)}, \qquad z_\kappa = \frac{\kappa}{\kappa-1}.$$

Since $f_\mathbf{q}(x)$ increases strictly with $\alpha$ for every $x > 0$, it is straightforward to see that $\mathbf{q}$ is admissible if and only if $\alpha \leq \alpha_\kappa$, and is (regular) critical if and only if $\alpha = \alpha_\kappa$. In the critical case $\alpha = \alpha_\kappa$, the partition function $Z_\mathbf{q}$ is given by $z_\kappa$ as defined above, and $\rho_\mathbf{q} = \kappa$ (see Figure 1).

Notice that when $\alpha \leq \alpha_\kappa$, the conditional law $P_\mathbf{q}(\cdot|\#F(M) = n)$, as considered in Theorem 2, coincides with the *uniform distribution* on the set

$$\{\mathbf{m} \in \mathcal{M}_+ : \deg(f) = 2\kappa \text{ for all } f \in F(\mathbf{m}), \#F(\mathbf{m}) = n\}$$

of $2\kappa$-angulations with $n$ faces, since $W_\mathbf{q}$ puts the same weight $\alpha^n$ on all the elements of this set. At the light of the discussion leading to Corollary 4, a more natural way to define this uniform distribution would have been to take the nonadmissible weight sequence $\mathbf{q} = \delta_\kappa$ in the first place, so $W_\mathbf{q}$ puts mass 1 on every $2\kappa$-angulation, and $P_\mathbf{q}^{F=n}$ is indeed uniform.

By further specialization of these results to the $\kappa = 2$ case of quadrangulations, we check that $\alpha_2 = 1/12, Z_{\alpha_2 \delta_2} = 2$, which is consistent with the results of [7]. Furthermore, the constant $(4\rho_\mathbf{q}/(9(Z_\mathbf{q}-1)))^{1/4}$ appearing in Theorem 2 is $(8/9)^{1/4}$ for $\kappa = 2$, which is consistent with the results of [8, 16].

1.5.2. *$q_i = \beta^i$.* Let $\beta > 0$, and let $q_i = \beta^i$ for $i \geq 1$, so that the weight of a map $\mathbf{m}$ is

$$W_\mathbf{q}(\{\mathbf{m}\}) = \prod_{f \in F(\mathbf{m})} q_{\deg(f)/2} = \beta^{1/2 \sum_{f \in F(\mathbf{m})} \deg(f)} = \beta^{\#A(\mathbf{m})},$$

(when summing the degrees of faces, each edge is counted twice). In this case,

$$f_\mathbf{q}(x) = \sum_{i \geq 0} x^i \beta^{i+1} N(i+1) = \beta \sum_{i \geq 0} (\beta x)^i \frac{(2i+1)!}{(i+1)! i!},$$

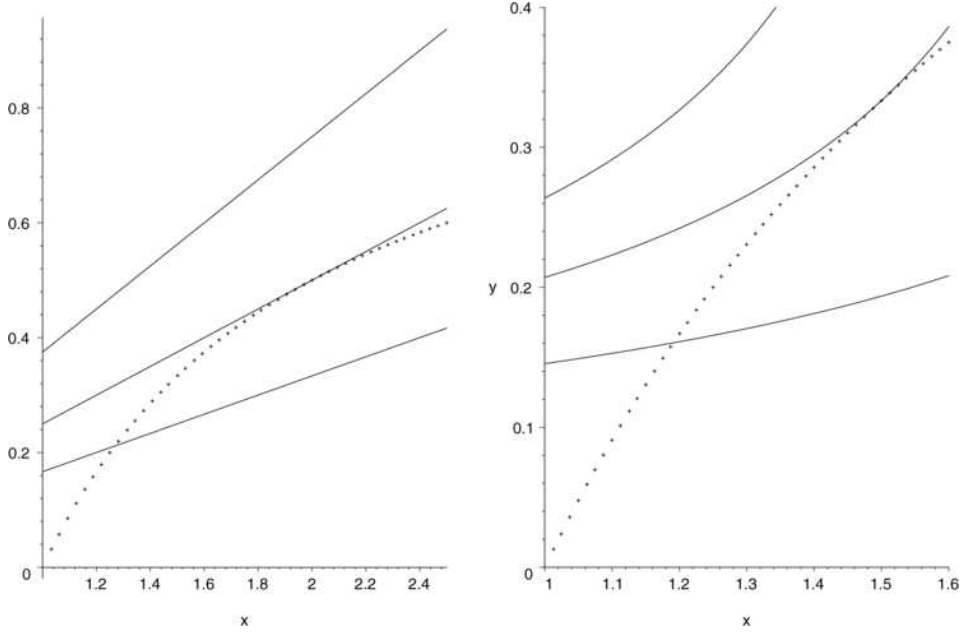

FIG. 1. *Example* 1.5.1: *drawing $f_{\mathbf{q}}$ for $\alpha = 1/18, 1/12, 1/8$ and $x \mapsto 1 - 1/x$ (dashed) in the case $\kappa = 2$ of quadrangulations. Example* 1.5.2: *drawing $f_{\mathbf{q}}$ for $\beta = 1/7, 1/8, 1/10$ and $x \mapsto 1 - 1/x$ (dashed).*

which is equal to

$$f_{\mathbf{q}}(x) = \frac{1}{2x} \sum_{j \geq 1} (\beta x)^j \frac{(2j)!}{j!^2} = \frac{1}{2x}((1 - 4\beta x)^{-1/2} - 1).$$

We see that $R_{\mathbf{q}} = (4\beta)^{-1}$ and that $f_{\mathbf{q}}(R_{\mathbf{q}}) = \infty$. Since we are looking for solutions of equation (2), which must be $> 1$, we see that the only interesting cases are when $\beta < 1/4$. More precisely, one can check that the equation $f_{\mathbf{q}}(z) = 1 - 1/z$ has real solutions if and only if $\beta \leq 1/8$, and these are given by

$$\frac{1 + 4\beta - \sqrt{1 - 8\beta}}{8\beta} \quad \text{and} \quad \frac{1 + 4\beta + \sqrt{1 - 8\beta}}{8\beta}.$$

These two solutions merge into a unique solution $3/2$ at $\beta = 1/8$, which is the value making $\mathbf{q}$ (regular) critical. This can be double checked by solving $z^2 f'_{\mathbf{q}}(z) = 1$, whose solution is $3/(16\beta)$. This gives $Z_{\mathbf{q}} = 3/2$ in the critical case $\beta = 1/8$, while $\rho_{\mathbf{q}} = 27/4$, and the value $(4\rho_{\mathbf{q}}/9)^{1/4}$ of Proposition 3 is $3^{1/4}$ (see Figure 1). Conditioning with respect to the number of vertices is indeed a bit more natural here: notice that we can rewrite $\mathbf{q} = \beta \bullet (\beta, i \geq 1)$. We thus obtain that $P^{S=n}_{\mathbf{q}}$ is equal to $P^{S=n}_{\beta}$, where $\beta$ stands (a bit improperly) for the constant sequence $q_i = \beta, i \geq 1$.

1.6. *Comments and organization of the paper.* As discussed in Section 1.5.1, the asymptotic behavior of the radius and profile of quadrangulations that are uniformly chosen in the set

$$\mathcal{Q}^n_+ = \{\mathbf{m} \in \mathcal{M}_+ : \#F(\mathbf{m}) = n, \deg(f) = 4 \text{ for all } f \in F(\mathbf{m})\}$$

is obtained as a particular case of Theorem 2 for $\mathbf{q} = 12^{-1}\delta_2$. Therefore, our results encompass in principle the results of Chassaing and Schaeffer [8] and Le Gall [16]. The reason why "in principle" is that these two papers deal with slightly different objects, namely rooted maps which are not pointed, and use the base point of the root edge as the reference point with respect to which geodesic distances are measured. Considering these objects would lead us to extra nontrivial complications. Roughly speaking, both pointing and rooting will allow us to study maps thanks to freely labeled trees, while simple rooting leads to considerations on labeled trees with a positivity constraint on labels. It is fortunate, however, that the scaling limits are the same for our model as in [8, 16]. On a very informal level, this indicates that the base vertex of the root edge in a uniform rooted quadrangulation with $n$ faces plays asymptotically the same role as a randomly picked base vertex. This is natural, since if we believe that a scaling limit for random maps exist, then a desirable feature of the limit would be that it (statistically) looks the same everywhere, and the singular role of the root in the discrete setting should vanish as the size of the map goes to infinity. On the other hand, we stress that rooting maps is not just a technical annoyance, but is really a crucial requirement in the methods used in most articles on the topic.

A natural question would be to ask whether similar techniques as ours could be used to prove similar invariance principles in nonbipartite cases (e.g., triangulations), using the more elaborate version of Bouttier, Di Francesco and Guitter's bijection for Eulerian planar maps. Although this makes the study slightly more intricate, this is indeed possible and will be addressed elsewhere.

The rest of the article is organized as follows. Section 2 introduces basic definitions for deterministic and random spatial trees, and shows how the bijection of Bouttier, Di Francesco and Guitter allows to interpret features of Boltzmann random bipartite maps in terms of functionals of certain two-type Galton–Watson (GW) trees coupled with a spatial motion. Section 3 provides the proof of Theorem 2 and Proposition 3, by introducing a new invariance principle for such spatial two-type GW trees (Theorem 8), in which the increments of the spatial motion can depend both on the type of the current vertex and on the local structure of the tree around the current vertex. This result is interesting in its own right. The proof of this invariance principle occupies the remaining Sections 4 and 5.

Finally, in Section 6, we show that under the hypothesis of Theorem 2 and Proposition 3, scaled bipartite maps converge to the Brownian map, introduced in [20]. This generalization is more or less straightforward, and then we just outline the procedure leading to this result.

## 2. Pushing Boltzmann planar maps to two-type spatial GW trees.

2.1. *Planar spatial trees.* Let $\mathbb{N}$ be the set of positive integers, and by convention let $\mathbb{N}^0 = \{\varnothing\}$. We define

$$\mathcal{U} = \bigsqcup_{n \geq 0} \mathbb{N}^n$$

(here and in the sequel, the symbol $\sqcup$ stands for the disjoint union) the set of all finite words with alphabet $\mathbb{N}$, using the notation $u = u_1 \cdots u_k \in \mathcal{U}$ where $u_1, \ldots, u_k \in \mathbb{N}$. If $u = u_1 \cdots u_k \in \mathcal{U}$ is such a word, we let $k = |u|$ be its length, with $|\varnothing| = 0$. If $u = u_1 \cdots u_k, v = v_1 \cdots v_{k'}$ are words, we let $uv$ be the concatenated word $u_1 \cdots u_k v_1 \cdots v_{k'}$, with the convention $\varnothing u = u\varnothing = u$. If $u = vw$ is a decomposition of a word $u$ as a concatenation, we say that $v$ is a prefix of $u$, and write $v \vdash u$. If $A$ is a subset of $\mathcal{U}$ and $u \in \mathcal{U}$, we let $uA = \{uv : v \in A\}$. The set $\mathcal{U}$ comes with the natural total lexicographical order $\preceq$, such that $u \preceq v$ if and only if either $u \vdash v$, or $u = wu', v = wv'$ with nonempty words $u', v'$ such that $u'_1 < v'_1$.

DEFINITION 2. A (rooted, planar) tree is a finite subset $\mathbf{t}$ of $\mathcal{U}$ that contains $\varnothing$, and such that $ui \in \mathbf{t}$ (with $u \in \mathcal{U}$ and $i \in \mathbb{N}$) implies that $u \in \mathbf{t}$ and $uj \in \mathbf{t}$ for all $1 \leq j \leq i$. We let $\mathcal{T}$ be the set of trees.

It is well known that this definition of rooted planar trees is equivalent to the graph-theoretic definition (a rooted planar map with no cycle), by associating every element $u \in \mathbf{t}$ with a vertex of a graph, and drawing edges from the vertex associated to $u$ to the ones associated to $u1, \ldots, uk \in \mathbf{t}$ "from left to right." If $\mathbf{t} \neq \{\varnothing\}$, the embedded graph thus obtained is rooted at the oriented edge from $\varnothing$ to 1. We call (with a slight abuse of notation) $\varnothing$ the *root* of $\mathbf{t}$. We call *vertices* the elements of a tree $\mathbf{t} \in \mathcal{T}$, the number $|u|$ is called the height of $u$, and the order $\preceq$ will be called the *depth-first order* on $\mathbf{t}$.

The set $\{ui : ui \in \mathbf{t}\}$ is interpreted as the set of *children* of $u \in \mathbf{t}$, and its cardinality is denoted by $c_{\mathbf{t}}(u)$. If $v = ui$ with $v, u \in \mathcal{U}$ and $i \in \mathbb{N}$, we say that $u$ is the father of $v$ and note $u = \neg v$. If $v \vdash u$ for $u, v \in \mathbf{t}$, we say that $v$ is an *ancestor* of $u$. If $\mathbf{t} \in \mathcal{T}$ and $u \in \mathbf{t}$ is a vertex, we let $\mathbf{t}_u = \{v \in \mathcal{U} : uv \in \mathbf{t}\}$ be the *fringe subtree* of $\mathbf{t}$ rooted at $u$. It is easily seen to be an element of $\mathcal{T}$. We also let $[\mathbf{t}]_u = \{u\} \cup (\mathbf{t} \setminus u\mathbf{t}_u)$ be the subtree of $\mathbf{t}$ which is *pruned* at $u$.

Next, let $\mathcal{T}_0, \mathcal{T}_1$ be two copies of $\mathcal{T}$. The picture that we have in mind is that if $\mathbf{t} \in \mathcal{T}_i$ for $i \in \{0, 1\}$, the mark $i$ is interpreted as a color (white 0 or black 1) that we assign to the root. All vertices at even height $|u|$ then earn the same color, while those at odd height earn the color $i+1 \bmod 2$. Although we should differentiate elements of $\mathcal{T}_0, \mathcal{T}_1, \mathcal{T}$ to be completely accurate, we keep the same notation $\mathbf{t}$ for elements of either of these sets. For $\mathbf{t} \in \mathcal{T}_i$, we let $\mathbf{t}^{(j)} = \{u \in \mathbf{t} : |u| = i + j \bmod 2\}$ to be the set of vertices of $\mathbf{t}$ with color $j$ [e.g., $\mathbf{t}^{(0)}$ is the set of vertices with even height if $\mathbf{t} \in \mathcal{T}_0$, and with odd height if $\mathbf{t} \in \mathcal{T}_1$]. This notation is the only one that actually distinguishes $\mathcal{T}_0$ from $\mathcal{T}_1$. In the sequel, we will often omit the mention of $\bmod 2$ when dealing with marks. For example, it is understood that $\mu_k, m_k$ stand for $\mu_{k \bmod 2}$ and its mean.

The definitions of children of a vertex, fringe subtrees and pruned subtrees extend naturally to $\mathcal{T}_0, \mathcal{T}_1$. The minor change is that if $\mathbf{t} \in \mathcal{T}_i$ for $i \in \{0, 1\}$, we take the convention that $\mathbf{t}_u \in \mathcal{T}_{|u|+i}$ (this should be clear from the intuitive picture that $i$ is the color of vertices at even heights in $\mathbf{t}$), and if $\mathbf{t} \in \mathcal{T}_i$, we still let $[\mathbf{t}]_u \in \mathcal{T}_i$ (the color of the root does not change).

DEFINITION 3. A spatial tree is a pair $(\mathbf{t}, \ell)$ where $\mathbf{t} \in \mathcal{T}$ and $\ell : \mathbf{t} \to \mathbb{R}$ is a labeling function that attributes a spatial position to every vertex. We let $\mathbb{T}$ be the set of spatial trees. Notice that for a fixed $\mathbf{t} \in \mathcal{T}$, taking a labeling $\ell$ is equivalent to attributing a label $\ell(\varnothing)$ to the root and determining the increments $\ell(u) - \ell(\neg u), u \in \mathbf{t} \setminus \{\varnothing\}$.

Again, we consider two copies $\mathbb{T}_0, \mathbb{T}_1$ of $\mathbb{T}$, that assign white or black color to the root, and alternate color between generations.

2.2. *Two-type spatial GW trees.* We now want to consider a particular family of multitype GW trees, in which vertices of type 0 only give birth to vertices of type 1 and vice versa. The following construction and discussion on a.s. finiteness of the tree is not the most economic one, but allows us to introduce some of the tools that will be needed later.

Let $\mu = (\mu_0, \mu_1)$ be a pair of probability distributions on $\mathbb{Z}_+$ with means $m_0$ and $m_1$, respectively. We make the basic assumption that $\mu$ is nondegenerate, that is, $\mu_0(1) + \mu_1(1) < 2$, and we exclude the trivial case $m_0 m_1 = 0$. We say that $\mu$ is subcritical if $m_0 m_1 < 1$, critical if $m_0 m_1 = 1$ and supercritical if $m_0 m_1 > 1$.

Consider a family of independent random variables $(X_u, u \in \mathcal{U})$ on some probability space $(\Omega, \mathcal{A}, P)$, such that $X_u$ with $|u|$ even all have law $\mu_0$ and $X_u$ with $|u|$ odd all have law $\mu_1$, and define

$$\xi = \{u = u_1 \cdots u_k \in \mathcal{U} : u_i \leq X_{u_1 \cdots u_{i-1}}, 1 \leq i \leq k\} \cup \{\varnothing\}.$$

It is elementary to prove that $\xi$ is a random subset of $\mathcal{U}$ that satisfies the properties of a tree. The only difference is that it might be infinite, though every vertex still has a finite number of children (local finiteness). We let $\widehat{\mathcal{T}}$ be the set of such possibly infinite trees which are locally finite, and keep the notation $c_{\hat{\mathbf{t}}}(u)$ for the number of children of $u \in \hat{\mathbf{t}} \in \widehat{\mathcal{T}}$. We also let $\widehat{\mathcal{T}}_0, \widehat{\mathcal{T}}_1$ be two copies of $\widehat{\mathcal{T}}$, and consider $\xi$ as a random element of $\widehat{\mathcal{T}}_0$. As before, if $\hat{\mathbf{t}} \in \widehat{\mathcal{T}}_i$ we let $\hat{\mathbf{t}}^{(j)}$ be the set $\{u \in \hat{\mathbf{t}} : |u| + i = j \bmod 2\}$. Notice that if $\mathbf{t} \in \mathcal{T}_0$, then we have

$$(6) \qquad P(\xi = \mathbf{t}) = \prod_{u \in \mathbf{t}^{(0)}} \mu_0(c_{\mathbf{t}}(u)) \prod_{u \in \mathbf{t}^{(1)}} \mu_1(c_{\mathbf{t}}(u)) = \prod_{u \in \mathbf{t}} \mu_{|u|}(c_{\mathbf{t}}(u)),$$

and these probabilities sum to 1 if and only if $\xi$ is a.s. an element of $\mathcal{T}_0$.

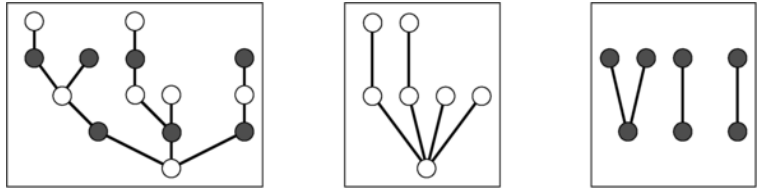

FIG. 2.   *The first frame depicts a tree $\mathbf{t} \in \mathcal{T}_0$: the white (resp. black) vertices stand for vertices of $\mathbf{t}^{(0)}$ (resp. $\mathbf{t}^{(1)}$). The second frame represents $\Gamma(\mathbf{t})$, and the third frame, $\Gamma'(1\mathbf{t})$ to be introduced later in Lemma* 11.

Now, for $\hat{\mathbf{t}} \in \widehat{\mathcal{T}}_0$, we introduce the mapping $\Gamma_{\hat{\mathbf{t}}} : \hat{\mathbf{t}}^{(0)} \to \mathcal{U}$ that associates with $\hat{\mathbf{t}}$ the tree having as number of vertices $\#\hat{\mathbf{t}}^{(0)}$, and which skips the odd generations of $\hat{\mathbf{t}}$, going straight from a vertex of $\hat{\mathbf{t}}^{(0)}$ to its grandsons. Formally, it is defined recursively by $\Gamma_{\hat{\mathbf{t}}}(\varnothing) = \varnothing$, and if $v \in \hat{\mathbf{t}}^{(0)}$ has grandchildren $vw_1, \ldots, vw_k \in \hat{\mathbf{t}}^{(0)}$, where $w_1, \ldots, w_k$ are words of two letters such that $w_1 \prec \cdots \prec w_k$, and $k = \sum_{1 \leq i \leq c_{\hat{\mathbf{t}}}(v)} c_{\hat{\mathbf{t}}}(vi)$ is the number of grandchildren of $v$, then $\Gamma_{\hat{\mathbf{t}}}(vw_l) = \Gamma_{\hat{\mathbf{t}}}(v)l$ for $1 \leq l \leq k$. We extend this to a mapping $\Gamma_{\hat{\mathbf{t}}} : \hat{\mathbf{t}} \to \mathcal{U}$ by letting $\Gamma_{\hat{\mathbf{t}}}(u) = \Gamma_{\hat{\mathbf{t}}}(\neg u)$ whenever $u \in \hat{\mathbf{t}}^{(1)}$.

We simply denote the tree $\Gamma_{\hat{\mathbf{t}}}(\hat{\mathbf{t}})$ by $\Gamma(\hat{\mathbf{t}})$. In particular, it is indeed an element of $\widehat{\mathcal{T}}$, and the root has $c_{\Gamma(\hat{\mathbf{t}})}(\varnothing) = \sum_{1 \leq k \leq c_{\hat{\mathbf{t}}}(\varnothing)} c_{\hat{\mathbf{t}}}(k)$ children. The tree is unmarked, because what we have done is to get rid of the vertices with color 1. Moreover, an easy recursion shows that $2|\Gamma_{\hat{\mathbf{t}}}(u)| = |u|$ for $u \in \hat{\mathbf{t}}$ with even height.

Now, it is elementary that $\Gamma(\xi)$ has the same law as the random element $\overline{\xi}$ of $\widehat{\mathcal{T}}$ that is defined as follows. Let $(\overline{X}_u, u \in \mathcal{U})$ be an i.i.d. sequence of random variables that have same distribution as

$$(7) \qquad \sum_{1 \leq k \leq X_\varnothing} X_k,$$

where $(X_u, u \in \mathcal{U})$ are the variables used to construct $\xi$. Then let

$$\overline{\xi} = \{u = u_1 \cdots u_k \in \mathcal{U} : u_i \leq \overline{X}_{u_1 \cdots u_{i-1}}, 1 \leq i \leq k\} \cup \{\varnothing\}.$$

By construction, $\overline{\xi}$ is a random variable in $\widehat{\mathcal{T}}$, and the process $(\#\{u \in \overline{\xi} : |u| = n\}, n \geq 0)$ is a GW process, whose offspring distribution is the law $\overline{\mu}$ of $\overline{X}_\varnothing$.

Moreover, the process is nondegenerate, that is, $\overline{\mu}(1) < 1$, as is easily deduced from the nondegeneracy condition on $\mu_0, \mu_1$. In particular, the process becomes extinct (i.e., $\overline{\xi}$ is finite) a.s. if and only if the mean $\overline{m}$ of $\overline{\mu}$ satisfies $\overline{m} \leq 1$. For any distribution $\nu$ on $\mathbb{Z}^+$, denote by $\mathbf{G}_\nu$ the generating function of $\nu$. We see from (7) that $\mathbf{G}_{\overline{\mu}} = \mathbf{G}_{\mu_0} \circ \mathbf{G}_{\mu_1}$. Differentiating this shows that the mean of the new offspring distribution is $\overline{m} = m_0 m_1$. Therefore, $\Gamma(\xi)$ is a.s. finite if and only if $m_0 m_1 \leq 1$, and the fact that $\xi$ is locally finite implies that the finiteness of $\Gamma(\xi)$ is equivalent to that of $\xi$.

By recalling formula (6), and considering as well the case where the roles of $\mu_0$ and $\mu_1$ are interchanged, we have proved:

PROPOSITION 5. *The formulas*

$$P_\mu^{(0)}(T = \mathbf{t}) = \prod_{u \in \mathbf{t}^{(0)}} \mu_0(c_\mathbf{t}(u)) \prod_{u \in \mathbf{t}^{(1)}} \mu_1(c_\mathbf{t}(u)), \qquad \mathbf{t} \in \mathcal{T}_0,$$

$$P_\mu^{(1)}(T = \mathbf{t}) = \prod_{u \in \mathbf{t}^{(0)}} \mu_0(c_\mathbf{t}(u)) \prod_{u \in \mathbf{t}^{(1)}} \mu_1(c_\mathbf{t}(u)), \qquad \mathbf{t} \in \mathcal{T}_1$$

*both sum to* 1 *when adding over* $\mathbf{t} \in \mathcal{T}_0$, *respectively* $\mathcal{T}_1$, *if and only if* $(\mu_0, \mu_1)$ *is* (*sub*)-*critical. In this case,* $P_\mu^{(0)}, P_\mu^{(1)}$ *are probability distributions, called the law of a* (*n alternating*) *two-type* (*sub*)-*critical GW tree, with root of type* 0, *respectively* 1.

Notice that the case $\mu_0 = \mu_1$ is that of a single type GW tree, and we do reobtain the usual a.s. extinction criterion $m \leq 1$ where $m$ is the expectation of $\mu_1$. In the sequel, by a two-type GW tree, we will always mean a random variable with a law of the form $P_\mu^{(0)}$ or $P_\mu^{(1)}$, as we will not be interested in the more general nonalternating cases.

We now couple the trees with a spatial displacement, in order to turn them into random elements of $\mathbb{T}$. For our purposes, we need to consider the case when the increments of the spatial motion depend both on the type and the degree of the neighboring vertex. To this end, let $(\nu_0^k, \nu_1^k, k \geq 1)$ be a family of probability distributions such that $\nu_0^k, \nu_1^k$ are defined on $\mathbb{R}^k$. Given $\mathbf{t} \in \mathcal{T}_i, i \in \{0, 1\}$, we let $(\mathbf{Y}_u, u \in \mathbf{t})$ be a family of independent random variables, such that for $u \in \mathbf{t}$ with $c_\mathbf{t}(u) = k$, $\mathbf{Y}_u = (Y_{u1}, \ldots, Y_{uk})$ has law $\nu_0^k$ for ever $u \in \mathbf{t}^{(0)}$ and $\nu_1^k$ for every $u \in \mathbf{t}^{(1)}$. This yields a family of random variables $(Y_u^\mathbf{t}, u \in \mathbf{t} \setminus \{\varnothing\})$, which we use as increments of a random labeling function on $\mathbf{t}$, that is, we set $\ell_\varnothing^\mathbf{t} = 0$ and

$$\ell_u^\mathbf{t} = \sum_{v \vdash u, v \neq \varnothing} Y_v^\mathbf{t}, \qquad u \in \mathbf{t}.$$

We denote by $\Lambda_\nu^\mathbf{t}$ the law of $(\ell_u^\mathbf{t}, u \in \mathbf{t})$. We let $\mathbb{P}_{\mu,\nu}^{(i)}$ be the law on $\mathbb{T}_i$ such that

$$\mathbb{P}_{\mu,\nu}^{(i)}(d\mathbf{t}\, d\ell) = P_\mu^{(i)}(d\mathbf{t}) \Lambda_\nu^\mathbf{t}(d\ell).$$

We usually let $(T,L)\colon \mathbb{T}_i \to \mathbb{T}_i$ be the identity mapping, so that under $\mathbb{P}^{(i)}_{\mu,\nu}$, $T$ has distribution $P^{(i)}_\mu$, and given $T = \mathbf{t}$, the labeling $L$ is $\Lambda^{\mathbf{t}}_\nu$-distributed. To avoid trivial degenerate cases, we will always implicitly suppose that there exists $i \in \{0,1\}, k \geq 1$ with $\mu_i(k) > 0$ and $\nu_i^k$ is not the Dirac mass at 0. We then say that the displacements laws $\nu_i^k$ are nondegenerate. We now have all the necessary background to describe the push-forward of the Boltzmann measures $P_\mathbf{q}$ under the bijection of Bouttier–Di Francesco–Guitter.

2.3. *The Bouttier–Di Francesco–Guitter bijection and its consequences.* The basic bijection presented in [6], Sections 2.1 and 2.2, is a bijection between the set of pointed *unrooted* bipartite planar maps (i.e., planar maps with a distinguished vertex), and the set of so-called *well-labeled mobiles*. These objects are *unrooted* planar trees together with a bipartite coloration of vertices (black or white say), such that white vertices carry positive integer labels, which satisfy a set of constraints. The nice feature of this bijection, aside from providing enumerative formulas, is that the faces of the initial map with degree $k$ are in one-to-one correspondence with the black vertices of the mobile with degree $k/2$, while the vertices of the map that are at distance $d > 0$ from the distinguished vertex are in one-to-one correspondence with white vertices of the mobile with label $d$.

It is explained in [6], Section 2.4, how a further rooting of the pointed map (giving a map of $\mathcal{M}_+$) allows to root the associated mobile at a white vertex, and lift the constraint that the labels are positive by subtracting the label of the root vertex to all other labels (recovering the initial labels amounts to subtracting the minimal label to every label and adding 1). We may reformulate their result as follows.

Let $\overline{\mathbb{T}} \subset \mathbb{T}_0$ be the set of pairs $(\mathbf{t},\ell)$, where the mark of the root of $\mathbf{t}$ is 0 and where the labeling function satisfies the following constraints:

- $\ell(\varnothing) = 0$,
- $\ell$ takes its values in $\mathbf{Z}$,
- $\ell(u) = \ell(\neg u)$ if $u \in \mathbf{t}^{(1)}$ (i.e., $|u|$ is odd),
- if $u \in \mathbf{t}^{(1)}$ has children $u1, \ldots, uk$, with $k = c_\mathbf{t}(u)$, then with the conventions $\ell(u0) = \ell(u) = \ell(u(k+1))$,

$$(8) \qquad \ell(uj) - \ell(u(j-1)) \in \{-1, 0, 1, 2, 3, \ldots\}, \qquad 1 \leq j \leq k+1.$$

For a given $\mathbf{t} \in \mathcal{T}_0$, a labeling $\ell$ satisfying these constraints [i.e., such that $(\mathbf{t},\ell) \in \overline{\mathbb{T}}$] is called compatible with $\mathbf{t}$, and the set of compatible labelings with $\mathbf{t}$ is denoted by $\mathcal{L}_\mathbf{t}$. Note that our conventions are slightly different from those of [6], where vertices of type 1 would be unlabeled. The difference is minor, since we consider that these vertices earn the label of their father.

PROPOSITION 6 (BDFG (Bouttier, Di Francesco and Guitter) bijection [6]). *There exists a bijection between the sets $\mathcal{M}_+$ and $\overline{\mathbb{T}}$, which we denote by $\Psi : \mathcal{M}_+ \to \overline{\mathbb{T}}$, that sends $\dagger$ on $\{\varnothing\}$ and satisfies the following extra properties. If $\mathbf{m} \in \mathcal{M}_+ \setminus \{\dagger\}$ and $(\mathbf{t}, \ell) = \Psi(\mathbf{m})$:*

- *Faces $f$ of $\mathbf{m}$ with degree $2k$ are in one-to-one correspondence with vertices $u \in \mathbf{t}^{(1)}$ (i.e., with $|u|$ odd) that have $k-1$ children. In particular, $\#F(\mathbf{m}) = \#\mathbf{t}^{(1)}$.*
- *Vertices $v$ of $\mathbf{m}$ such that $d_{\mathbf{m}}(v, \mathfrak{r}) = d > 0$ are in one-to-one correspondence with vertices $u \in \mathbf{t}^{(0)}$ (i.e., with $|u|$ even) with $\ell(u) - \min_{u' \in \mathbf{t}} \ell(u') + 1 = d$. In particular, $\#S(\mathbf{m}) = \#\mathbf{t}^{(0)} + 1$,*

$$(9) \qquad \mathscr{R}(\mathbf{m}) = \max_{u \in \mathbf{t}} \ell(u) - \min_{u \in \mathbf{t}} \ell(u) + 1,$$

*and for any $k \geq 0$*

$$(10) \quad \mathscr{I}^{\mathbf{m}}(k) = \frac{1}{\#\mathbf{t}^{(0)} + 1} \left( \# \left\{ u \in \mathbf{t}^{(0)} : \ell(u) - \min_{u' \in \mathbf{t}} \ell(u') + 1 = k \right\} + \mathbb{1}_{\{k=0\}} \right).$$

A short description of $\Psi^{-1}$ can be found in Section 6.

Except from the trivial difference explained before the statement of this proposition, the only difference with [6], Section 2.4, is that the case of the vertex-map is not considered there, and the mobiles always have at least one white vertex and one black vertex. This distinction is important in our study, as we will see after the next proposition. The key observation of this paper is given by the following statement, which gives the image measure of the Boltzmann distributions $P_{\mathbf{q}}$ on $\mathcal{M}_+$ by $\Psi$. We let $\Pi(\mathbf{m}) = \mathbf{t}, \Pi'(\mathbf{m}) = \ell$ whenever $\Psi(\mathbf{m}) = (\mathbf{t}, \ell)$.

PROPOSITION 7. *Let $\mathbf{q}$ be an admissible weight sequence, and define two probability distributions $(\mu_0, \mu_1)$ by*

$$\mu_0(k) = Z_{\mathbf{q}}^{-1} f_{\mathbf{q}}(Z_{\mathbf{q}})^k, \qquad k \geq 0,$$

*the geometric law with parameter $f_{\mathbf{q}}(Z_{\mathbf{q}})$ (as defined in the Introduction), and*

$$\mu_1(k) = \frac{Z_{\mathbf{q}}^k N(k+1) q_{k+1}}{f_{\mathbf{q}}(Z_{\mathbf{q}})}, \qquad k \geq 0.$$

*Also, for every $k \geq 1$, let $\nu_0^k$ be the Dirac mass at $0 \in \mathbb{R}^k$, and $\nu_1^k$ be the law on $\mathbf{Z}^k$ of $(X_1, X_1 + X_2, \ldots, X_1 + X_2 + \cdots + X_k)$ where $(X_1, \ldots, X_{k+1})$ is uniform in the set*

$$\{(x_1, \ldots, x_{k+1}) \in (\mathbf{Z}_+ \cup \{-1\})^{k+1} : x_1 + \cdots + x_{k+1} = 0\}.$$

*Then the two-type GW tree associated with $\mu_0, \mu_1$ is (sub)-critical, and $\Psi(M)$ under $P_{\mathbf{q}}$ has law $\mathbb{P}_{\mu,\nu}^{(0)}$.*

*Moreover, the weight sequence* **q** *is critical in the sense of Definition* 1 *if and only if* $\Pi(M)$ *under* $P_{\mathbf{q}}$ *is a critical two-type GW tree, and* **q** *is regular critical if and only if it is critical and* $\mu_1$ *admits small exponential moments, namely* $\langle \mu_1, \exp(a\cdot)\rangle < \infty$ *for some* $a > 0$.

This explains why considering † as a map is important for our concern: otherwise, the previous statement would not be true, since the root vertex of $\Pi(M)$ under $P_{\mathbf{q}}$ would be constrained to have at least one child, and the tree would not enjoy the GW property (in fact it would enjoy it everywhere but at the root, which would make the forthcoming discussion tedious). It is also this convention that allows neat statements in Definition 1 and Proposition 1.

The first step in proving Proposition 7 is to compute the cardinality of the set $\mathcal{L}_{\mathbf{t}}$ of labelings that are compatible with some $\mathbf{t} \in \mathcal{T}_0$. For such $\mathbf{t}$, the constraint (8) says that for every $u \in \mathbf{t}^{(1)}$, the label differences $\ell(uj) - \ell(u(j-1)), 1 \leq j \leq k+1$ must be in $\mathbf{Z}_+ \cup \{-1\}$ and sum to 0 because of the convention $\ell(u0) = \ell(u) = \ell(u(k+1))$. This is the same as the number of $k+1$-tuples $(\ell(uj) - \ell(u(j-1)) + 2, 1 \leq j \leq k+1)$ forming a composition of the integer $2k+2$ with $k+1$ positive parts. The number of such compositions is equal to $\binom{2k+1}{k+1} = N(k+1)$ with the conventions of Section 1. Since the label of the root of $\mathbf{t}$ is fixed to 0, the number of admissible labelings $\ell$ of $\mathbf{t}$ is therefore equal to

$$\#\mathcal{L}_{\mathbf{t}} = \prod_{u \in \mathbf{t}^{(1)}} N(c_{\mathbf{t}}(u) + 1). \tag{11}$$

Next, let $\mathbf{q} = (q_i, i \geq 1)$ be any nonnegative weight sequence, not necessarily admissible. Let $\mathbf{m} \in \mathcal{M}_+$. Then, by letting $(\mathbf{t}, \ell) = \Psi(\mathbf{m})$ and using Proposition 6 we get that

$$W_{\mathbf{q}}(M = \mathbf{m}) = \prod_{f \in F(\mathbf{m})} q_{\deg(f)/2}$$
$$= \prod_{u \in \mathbf{t}^{(1)}} q_{c_{\mathbf{t}}(u)+1} = W_{\mathbf{q}}(\Psi(M) = (\mathbf{t}, \ell)). \tag{12}$$

This quantity is independent of the values taken by $\ell$, for any $\ell$ compatible with $\mathbf{t}$, so by (11),

$$W_{\mathbf{q}}(\Psi(M) \in \{(\mathbf{t}, \ell) : \ell \in \mathcal{L}_{\mathbf{t}}\}) = \sum_{\ell \in \mathcal{L}_{\mathbf{t}}} \prod_{u \in \mathbf{t}^{(1)}} q_{c_{\mathbf{t}}(u)+1}$$
$$= \prod_{u \in \mathbf{t}^{(1)}} N(c_{\mathbf{t}}(u) + 1) q_{c_{\mathbf{t}}(u)+1}.$$

Otherwise said,

$$W_{\mathbf{q}}(\Pi(M) = \mathbf{t}) = \prod_{u \in \mathbf{t}^{(1)}} N(c_{\mathbf{t}}(u) + 1) q_{c_{\mathbf{t}}(u)+1}. \tag{13}$$

We are now ready to prove Proposition 7. Some of the computations appear implicitly in [6].

PROOF OF PROPOSITION 7. Suppose $\mathbf{q}$ is admissible, so $Z_\mathbf{q} = W_\mathbf{q}(\mathcal{M}_+) < \infty$. Notice that for any tree $\mathbf{t} \in \mathcal{T}_0$, one has

$$\sum_{u \in \mathbf{t}^{(0)}} c_\mathbf{t}(u) = \#\mathbf{t}^{(1)}, \qquad \sum_{u \in \mathbf{t}^{(1)}} c_\mathbf{t}(u) = \#\mathbf{t}^{(0)} - 1,$$

since $\varnothing$ is the only vertex which has no ancestor. Therefore, we may redisplay (13) as

$$W_\mathbf{q}(\Pi(M) = \mathbf{t}) = \left(\frac{1}{Z_\mathbf{q}}\right)^{\#\mathbf{t}^{(0)}-1} \prod_{u \in \mathbf{t}^{(1)}} (Z_\mathbf{q})^{c_\mathbf{t}(u)} N(c_\mathbf{t}(u) + 1) q_{c_\mathbf{t}(u)+1},$$

so finally

$$(14) \quad P_\mathbf{q}(\Pi(M) = \mathbf{t}) = \left(\frac{1}{Z_\mathbf{q}}\right)^{\#\mathbf{t}^{(0)}} \prod_{u \in \mathbf{t}^{(1)}} (Z_\mathbf{q})^{c_\mathbf{t}(u)} N(c_\mathbf{t}(u) + 1) q_{c_\mathbf{t}(u)+1}.$$

We know that summing this formula over $\mathbf{t} \in \mathcal{T}_0$ gives 1. But notice that any tree $\mathbf{t} \in \mathcal{T}_0$ can be written as $\mathbf{t} = \{\varnothing\} \cup 1\mathbf{t}_1 \cup \cdots \cup k\mathbf{t}_k$, if $k = c_\mathbf{t}(\varnothing)$, and with $\mathbf{t}_i \in \mathcal{T}_1$, $1 \le i \le k$. So summing the last formula over $\mathbf{t} \in \mathcal{T}_0$ amounts to sum over $k \geq 0$ and $\mathbf{t}_{(1)}, \ldots, \mathbf{t}_{(k)} \in \mathcal{T}_1$, and factorize the term $1/Z_\mathbf{q}$ that involves the root of $\mathbf{t}$, so

$$1 = \sum_{k \geq 0} \frac{1}{Z_\mathbf{q}} \sum_{\mathbf{t}_{(1)}, \ldots, \mathbf{t}_{(k)} \in \mathcal{T}_1} \prod_{i=1}^{k} \left( \left(\frac{1}{Z_\mathbf{q}}\right)^{\#\mathbf{t}_{(i)}^{(1)}} \prod_{u \in \mathbf{t}_{(i)}^{(0)}} (Z_\mathbf{q})^{c_\mathbf{t}(u)} N(c_\mathbf{t}(u) + 1) q_{c_\mathbf{t}(u)+1} \right)$$

(15)
$$= \frac{1}{Z_\mathbf{q}} \sum_{k \geq 0} \left( \sum_{\mathbf{t} \in \mathcal{T}_1} \left(\frac{1}{Z_\mathbf{q}}\right)^{\#\mathbf{t}^{(1)}} \prod_{u \in \mathbf{t}^{(0)}} (Z_\mathbf{q})^{c_\mathbf{t}(u)} N(c_\mathbf{t}(u) + 1) q_{c_\mathbf{t}(u)+1} \right)^k.$$

But the quantity which is raised to the successive integer powers can be decomposed by a similar method, and is equal to

$$\sum_{k \geq 0} (Z_\mathbf{q})^k N(k+1) q_{k+1} \left( \sum_{\mathbf{t} \in \mathcal{T}_0} \left(\frac{1}{Z_\mathbf{q}}\right)^{\#\mathbf{t}^{(0)}} \prod_{u \in \mathbf{t}^{(1)}} (Z_\mathbf{q})^{c_\mathbf{t}(u)} N(c_\mathbf{t}(u) + 1) q_{c_\mathbf{t}(u)+1} \right)^k.$$

This time, the right-most quantity which is raised to the power $k$ is nothing but the sum of (14) over $\mathbf{t} \in \mathcal{T}_0$, that is $P_\mathbf{q}(\mathcal{M}_+) = 1$, with which we started. Thus, the last expression is nothing but $f_\mathbf{q}(Z_\mathbf{q})$. Plugging this in (15), this leads to

$$1 = \frac{1}{Z_\mathbf{q}} \sum_{k \geq 0} f_\mathbf{q}(Z_\mathbf{q})^k.$$

This yields both that $f(Z_\mathbf{q}) < 1$ and that $Z_q$ is solution of equation (2). Therefore, the definition of $\mu_1$ in the statement of the theorem makes sense and defines a probability distribution.

With this in hand, we can rewrite (14) and easily get

$$P_\mathbf{q}(\Pi(M) = \mathbf{t}) = \prod_{u \in \mathbf{t}^{(0)}} \mu_0(c_\mathbf{t}(u)) \prod_{u \in \mathbf{t}^{(1)}} \mu_1(c_\mathbf{t}(u)).$$

Since these probabilities sum to 1 when summing over $\mathbf{t}$, we get that $\Pi(M)$ under $P_\mathbf{q}$ is indeed a (sub) critical two-type GW tree by Proposition 5, and with the claimed offspring distributions. Obtaining the law of the labeling $\Pi'(M)$ given $\Pi(M)$ is then easy as (12) may be rewritten

(16) $\qquad P_\mathbf{q}(\Psi(M) = (\mathbf{t}, \ell)) = P_\mathbf{q}(\Pi(M) = \mathbf{t}) \prod_{u \in \mathbf{t}^{(1)}} \frac{1}{N(c_\mathbf{t}(u) + 1)}.$

Therefore, given $\Pi(M) = \mathbf{t}$, the labeling is uniform among all compatible labelings $\mathcal{L}_\mathbf{t}$. We can reexpress this by saying that [still under $P_\mathbf{q}(\cdot|\Pi(M) = \mathbf{t})$] the increments $(\Pi'(M)(uj) - \Pi'(M)(u(j-1)), 1 \leq j \leq c_\mathbf{t}(u) + 1), u \in \mathbf{t}^{(1)}$, with the cyclic convention of (8), are independent as $u$ varies, and uniform among all the $N(c_\mathbf{t}(u) + 1)$ increments sequences that are respectively allowed. Equivalently, under $P_\mathbf{q}(\cdot|\Pi(M) = \mathbf{t})$, the increments $(\Pi'(M)(uj) - \Pi'(M)(u), 1 \leq j \leq c_\mathbf{t}(u))$ are independent as $u$ varies in $\mathbf{t}^{(1)}$ and have the law $\nu_1^{c_\mathbf{t}(u)}$ of the statement. Increments $(\Pi'(M)(uj) - \Pi'(M)(u), 1 \leq j \leq c_\mathbf{t}(u))$ for vertices $u \in \mathbf{t}^{(0)}$ are a.s. equal to 0, and contribute to an invisible factor of 1 to (16), which explains the definition of $\nu_0^k$.

To prove the criticality statement, it suffices to compute the expectations $m_0$ and $m_1$ of $\mu_0$ and $\mu_1$. The expectation of the geometric law $\mu_0$ is equal to $Z_\mathbf{q} - 1 = Z_\mathbf{q} f_\mathbf{q}(Z_\mathbf{q})$, while

$$m_1 = \frac{1}{f_\mathbf{q}(Z_\mathbf{q})} \sum_{k \geq 0} k Z_\mathbf{q}^k N(k+1) q_{k+1} = \frac{Z_\mathbf{q} f'_\mathbf{q}(Z_\mathbf{q})}{f_\mathbf{q}(Z_\mathbf{q})}.$$

The product $m_0 m_1$ is thus $Z_\mathbf{q}^2 f'_\mathbf{q}(Z_\mathbf{q})$, which must be $\leq 1$ (the tree is subcritical), which shows that $Z_\mathbf{q}$ must be the smallest solution of (2), by the classification of solutions of (2) given in the introduction. The weight sequence $\mathbf{q}$ is then critical in the sense of Definition 1 if and only if $m_0 m_1 = 1 = Z_\mathbf{q}^2 f'_\mathbf{q}(Z_\mathbf{q})$, that is, $\Pi(M)$ under $P_\mathbf{q}$ is critical. If $\mathbf{q}$ is critical, it is regular critical if and only if $Z_\mathbf{q} < R_\mathbf{q}$ where $R_\mathbf{q}$ is the radius of convergence of $f_\mathbf{q}$, and it is easy to see that this is equivalent to $\langle \mu_1, \exp(a\cdot) \rangle < \infty$ for some $a > 0$. $\square$

PROOF OF PROPOSITION 1. We have already noticed that if $\mathbf{q}$ is admissible, then $Z_\mathbf{q}$ satisfies (2) and is the smallest solution, that is, the one satisfying $Z_\mathbf{q}^2 f_\mathbf{q}(Z_\mathbf{q}) \leq 1$.

Conversely, suppose that (2) admits a solution. Then thanks to the classification of solutions of the Introduction, we know that one of the solutions, say $z$, satisfies $z^2 f'_\mathbf{q}(z) \leq 1$. In a similar way as in the proof of Proposition 7, we can write (13) as

$$\frac{W_\mathbf{q}(\Pi(M) = \mathbf{t})}{z} = \left(\frac{1}{z}\right)^{\#\mathbf{t}^{(0)}} f_\mathbf{q}(z)^{\#\mathbf{t}^{(1)}} \prod_{u \in \mathbf{t}^{(1)}} \frac{z^{c_\mathbf{t}(u)} N(c_\mathbf{t}(u) + 1) q_{c_\mathbf{t}(u)+1}}{f_\mathbf{q}(z)}$$
$$= \prod_{u \in \mathbf{t}^{(0)}} \mu'_0(c_\mathbf{t}(u)) \prod_{u \in \mathbf{t}^{(1)}} \mu'_1(c_\mathbf{t}(u)),$$

where $\mu'_0(k) = z^{-1} f_\mathbf{q}(z)^k$ and $\mu'_1(k) = z^k N(k+1) q_{k+1}/f_\mathbf{q}(z)$, for $k \geq 0$, are two probability distributions [for $\mu'_0$, use the fact that $f_\mathbf{q}(z) = 1 - z^{-1}$]. Since moreover these distributions have means $m'_0 = z - 1 = z f_\mathbf{q}(z)$ and $m'_1 = z f'_\mathbf{q}(z)/f_\mathbf{q}(z)$, whose product is $z^2 f'_\mathbf{q}(z) \leq 1$, we finally recognize that the image of $W_\mathbf{q}/z$ under $\Pi$ is the (probability) law of a (sub) critical two-type GW tree. This shows $z = Z_\mathbf{q} < \infty$, hence the result. $\square$

**3. An invariance principle for spatial GW trees.** In view of Proposition 6, and in particular the formulas (9) and (10), and Proposition 7, the asymptotic behavior of the radius and profile of random maps under $P_\mathbf{q}(\cdot|\#F(M) = n)$ [resp. $P_\mathbf{q}(\cdot|\#S(M) = n)$], with $\mathbf{q}$ critical boils down to that of the labels distribution in a critical spatial GW tree with law $\mathbb{P}^{(0)}_{\mu,\nu}(\cdot|\#T^{(1)} = n)$ [resp. $\mathbb{P}^{(0)}_{\mu,\nu}(\cdot|\#T^{(0)} = n - 1)$]. We now state an invariance principle for such trees.

3.1. *The invariance principle.* For $\mathbf{t} \in \mathcal{T}$, let $\varnothing = u(0) \prec u(1) \prec \cdots \prec u(\#\mathbf{t} - 1)$ be the list of vertices of $\mathbf{t}$ in depth-first order. We let $H^\mathbf{t}_k = |u(k)|, 0 \leq k \leq \#\mathbf{t} - 1$, and we construct a continuous piecewise linear process $(H^\mathbf{t}_t, 0 \leq t \leq \#\mathbf{t} - 1)$ by linear interpolation between integer points. The process $H^\mathbf{t}$ is called the *height process* of $\mathbf{t}$.

Next, for a labeled tree $(\mathbf{t}, \ell) \in \mathbb{T}$, we let $(S^{\mathbf{t},\ell}_k = \ell(u(k)), 0 \leq k \leq \#\mathbf{t} - 1)$ be the head of the discrete snake associated with $(\mathbf{t}, \ell)$. We extend this process into a piecewise linear continuous process $(S^{\mathbf{t},\ell}_t, 0 \leq t \leq \#\mathbf{t} - 1)$ by interpolating between integer values.

Let $(\mu_0, \mu_1)$ be a nondegenerate critical two-type offspring distribution, and $(\nu^k_i, i \in \{0, 1\}, k \geq 1)$ be a centered spatial displacement law ($\langle \nu^k_i, x \rangle = 0$). We let $\mathbb{P}^{(i)} = \mathbb{P}^{(i)}_{\mu,\nu}$ for simplicity. Let $m_0, m_1, \sigma_0^2, \sigma_1^2$ be the means and variances of $\mu_0, \mu_1$, and define

$$(17) \qquad \sigma = \frac{1}{2}\sqrt{\sigma_0^2 \frac{1 + m_1}{m_0} + \sigma_1^2 \frac{1 + m_0}{m_1}} \in (0, \infty].$$

Also, for $i \in \{0,1\}$, $k \geq 1$ and $1 \leq l \leq k$, let $\Sigma_i^{k,l} = \sqrt{\langle \nu_i^k, x_l^2 \rangle}$ be the square root of the variance of the $l$th component of a random vector with law $\nu_i^k$, and

$$\Sigma_i^k = \sqrt{\langle \nu_i^k, |x|^2 \rangle} = \left( \sum_{l=1}^k (\Sigma_i^{k,l})^2 \right)^{1/2}, \tag{18}$$

where $|x|$ is the Euclidean norm of $x \in \mathbb{R}^k$. We define

$$\Sigma = \sqrt{\frac{1}{2} \sum_{k \geq 1} \left[ \frac{\mu_0(k)}{m_0} (\Sigma_0^k)^2 + \frac{\mu_1(k)}{m_1} (\Sigma_1^k)^2 \right]}. \tag{19}$$

Recall the definition of $\mathbb{N}^{(1)}$, Section 1.3. We endow $\mathcal{C}(\mathbb{R}_+, \mathbb{R})$ with the uniform topology, and $\mathcal{C}(\mathbb{R}_+, \mathbb{R})^2$ with the product topology. The invariance principle states as:

THEOREM 8.  *Let $(\mu_0, \mu_1)$ be a critical nondegenerate offspring distribution, and suppose it admits some exponential moments. Let $(\nu_0^k, \nu_1^k, k \geq 1)$ be nondegenerate spatial displacement laws which are centered, and such that there exists some $\eta > 0$ such that for $i \in \{0,1\}$ and $k \geq 1$,*

$$M_i^k := \int_{\mathbb{R}^k} |x|^{4+\eta} \nu_i^k(dx) < \infty.$$

*Last, assume that for some $D > 0$, as $k \to \infty$,*

$$M_0^k \vee M_1^k = O(k^D). \tag{20}$$

*Then, the constants $\sigma, \Sigma > 0$ are finite, and the following convergence in distribution holds on $\mathcal{C}(\mathbb{R}_+, \mathbb{R})^2$, for $i, j \in \{0, 1\}$:*

$$\left( \left( \frac{H^T_{(\#T-1)t}}{n^{1/2}} \right)_{0 \leq t \leq 1}, \left( \frac{S^{T,L}_{(\#T-1)t}}{n^{1/4}} \right)_{0 \leq t \leq 1} \right) \quad \text{under } \mathbb{P}^{(i)}(\cdot | \#T^{(j)} = n)$$

$$\xrightarrow[n \to \infty]{(d)} \left( \left( \frac{2\sqrt{1+m_j}}{\sigma} \mathrm{e}_t \right)_{0 \leq t \leq 1}, \left( \frac{\sqrt{2}\Sigma(1+m_j)^{1/4}}{\sigma^{1/2}} \mathrm{r}_t \right)_{0 \leq t \leq 1} \right) \quad \text{under } \mathbb{N}^{(1)},$$

*where by convention, $n$ goes to $+\infty$ along the values for which the conditioning event has positive probability.*

One of the key ingredients in the proof of this result is the forthcoming Lemma 15, which deals with the repartition between vertices of either type in a conditioned two-type GW tree. In order to be able to prove Theorem 2 right away, we give a simpler statement for now. Let $\mathbf{t} \in \mathcal{T}_i$ for some $i \in \{0,1\}$. For $0 \leq k \leq \#\mathbf{t} - 1$ and $j \in \{0,1\}$, we let

$$J^{(j)}_{\mathbf{t}}(k) = \mathrm{Card}(\mathbf{t}^{(j)} \cap \{u(0), \ldots, u(k)\})$$

be the counting process for the ranks of the vertices of $\mathbf{t}^{(j)}$, when $\mathbf{t}$ is visited in depth-first order. We extend it into a right-continuous nondecreasing function on $[0, \#\mathbf{t} - 1]$ by letting $J_\mathbf{t}^{(j)}(t) = J_\mathbf{t}^{(j)}([t])$. The renormalized function $\overline{J}_\mathbf{t}^{(j)} = (J_\mathbf{t}^{(j)}((\#\mathbf{t} - 1)t)/\#\mathbf{t}^{(j)}, 0 \leq t \leq 1)$ is the distribution function for the probability measure putting equal mass on each number $k/(\#\mathbf{t} - 1)$ with $k \in \mathbf{Z}_+$ such that $u(k) \in \mathbf{t}^{(j)}$. The following result says that vertices of either type are homogeneously displayed in a GW tree conditioned to be large.

LEMMA 9. *Let $\mu_0, \mu_1$ be nondegenerate critical, and admitting small exponential moments. Then for $i, j \in \{0, 1\}$, under $P^{(i)}(\cdot | \#T^{(j)} = n)$, the processes $(\overline{J}_T^{(c)}(t), 0 \leq t \leq 1)$, $c \in \{0, 1\}$ converge in probability to the identity $(t, 0 \leq t \leq 1)$, for the uniform norm.*

We end the present section by showing how Theorem 8 and Lemma 9 allow to prove Theorem 2 and Proposition 3.

3.2. *Computation of the scaling constants associated with random maps.* Let $\mathbf{q}$ be a regular critical admissible weight sequence. Then we know that $\Psi(M)$ under $P_\mathbf{q}$ has law $P_{\mu,\nu}^{(0)}$, where $\mu, \nu$ are defined as in Proposition 7. We also know from this proposition that $\mu_0, \mu_1$ admit some exponential moments ($\mu_0$ because the law is geometric, and $\mu_1$ because $\mathbf{q}$ is regular critical). Also, it is plainly nondegenerate.

On the other hand, we have to check that the $\nu_i^k$ are centered and satisfy the moments conditions of Theorem 8. For $\nu_0^k$ it is trivial (these are Dirac masses at 0). Since $\nu_1^k$ is carried by the set $[-k, k]^k$, it is straightforward that its marginals have moments of order 5 which grow at most like $k^5$.

Next, we compute the constants $\sigma, \Sigma$ associated with $\mu, \nu$. On the one hand, $\mu_0$ has mean $m_0 = Z_\mathbf{q} - 1 = m_1^{-1}$, and variance $\sigma_0^2 = Z_\mathbf{q}(Z_\mathbf{q} - 1)$. Also, $\mathbf{G}_{\mu_1}(x) = f_\mathbf{q}(xZ_\mathbf{q})/f_\mathbf{q}(Z_\mathbf{q})$ by definition, and by differentiating, $\mu_1$ has variance

$$\sigma_1^2 = \frac{Z_\mathbf{q}^2 f_\mathbf{q}''(Z_\mathbf{q})}{f_\mathbf{q}(Z_\mathbf{q})} + \frac{Z_\mathbf{q} - 2}{(Z_\mathbf{q} - 1)^2} = \frac{Z_\mathbf{q}^3 f_\mathbf{q}''(Z_\mathbf{q})}{Z_\mathbf{q} - 1} + \frac{Z_\mathbf{q} - 2}{(Z_\mathbf{q} - 1)^2}.$$

This gives, after some simplifications,

$$\sigma = \frac{\sqrt{Z_\mathbf{q} \rho_\mathbf{q}}}{2},$$

where $\rho_\mathbf{q}$ is defined at (5).

On the other hand, we have to compute $\sum_{1 \leq l \leq k} (\Sigma_1^{k,l})^2$ to give the value of $\Sigma$ (notice that $\Sigma_0^{k,l} = 0$ for every $k, l$). Recall that $\nu_1^k$ is the law of $(X_1, X_1 + X_2, \ldots, X_1 + \cdots + X_k)$, where $(X_1, \ldots, X_{k+1})$ has a uniform law in $\{(x_1, \ldots, x_{k+1}) \in (\mathbf{Z}_+ \cup \{-1\})^{k+1} : x_1 + \cdots + x_{k+1} = 0\}$. But then, $(X_1, \ldots, X_{k+1})$ is exchangeable, and $E[X_1] = (k+1)^{-1}E[X_1 + \cdots + X_{k+1}] = 0$, so the variables $X_l, 1 \leq l \leq k$, are centered, as well as the marginals of $\nu_1^k$. Moreover, it holds by exchangeability (this argument was suggested by a referee) that

$$(\Sigma_1^{k,l})^2 = \text{Var}(X_1 + \cdots + X_l) = l\,\text{Var}(X_1) + l(l-1)\,\text{Cov}(X_1, X_2).$$

Since $\text{Var}(X_1 + \cdots, X_{k+1}) = 0$, we obtain that $\text{Cov}(X_1, X_2) = -\text{Var}(X_1)/k$. It remains to compute the variance of $X_1$. Using the interpretation in terms of compositions, one finds easily that

$$(21) \qquad \mathbb{P}(X_1 = l) = \binom{2k-l-1}{k-1} \Big/ \binom{2k+1}{k+1} \qquad \text{for } -1 \leq l \leq k$$

and then, since $\frac{a+1}{b+1}\binom{a}{b} = \binom{a+1}{b+1}$,

$$\frac{E((2k-X_1)(2k+1-X_1))}{k(k+1)} = \sum_{l=-1}^{k} \binom{2k+1-l}{k+1} \Big/ \binom{2k+1}{k+1}$$

$$= \sum_{i=0}^{k+1} \binom{k+1+i}{k+1} \Big/ \binom{2k+1}{k+1}$$

$$= \binom{2k+3}{k+2} \Big/ \binom{2k+1}{k+1}$$

from which we get $\text{Var}(X_1) = 2k/(k+2)$. Finally, this gives $(\Sigma_1^{k,l})^2 = 2l(k-l+1)/(k+2)$, and by summing this for $1 \leq l \leq k$,

$$(\Sigma_1^k)^2 = \frac{k(k+1)}{3}.$$

We obtain

$$\Sigma = \sqrt{\frac{(Z_\mathbf{q}-1)}{2} \sum_{k \geq 1} \mu_1(k)\frac{k(k+1)}{3}} = \sqrt{\frac{Z_\mathbf{q}^3 f_\mathbf{q}''(Z_\mathbf{q}) + 2}{6}} = \sqrt{\frac{\rho_\mathbf{q}}{6}}.$$

Finally, we obtain that the scaling constants $D_\mathbf{q}$ and $C_\mathbf{q}$ appearing respectively in front of e and r in Theorem 8, for $j = 1$, are

$$D_\mathbf{q} = 4(\rho_\mathbf{q}(Z_\mathbf{q}-1))^{-1/2}, \qquad C_\mathbf{q} = \frac{\sqrt{2}\Sigma(1+m_1)^{1/4}}{\sigma^{1/2}} = \left(\frac{4\rho_\mathbf{q}}{9(Z_\mathbf{q}-1)}\right)^{1/4}.$$

3.3. *Proof of Theorem* 2. (i) From (9), we know that $\mathscr{R}(M)$ under $P_\mathbf{q}(\cdot|\#F(M) = n)$ has the same law as $\max_{u \in T} L(u) - \min_{u \in T} L(u) + 1$, under $\mathbb{P}^{(0)}(\cdot|\#T^{(1)} = n)$. In turn, this is equal to $1 + \max_{0 \leq t \leq \#T-1} S_t^{T,L} -$ $\min_{0 \leq t \leq \#T-1} S_t^{T,L}$. The convergence of the second component in Theorem 8 entails that $n^{-1/4}\mathscr{R}(M)$ converges in distribution to $C_{\mathbf{q}}\Delta(\mathrm{r})$, as claimed. Indeed, the convergence holds for the uniform topology, under which the mapping $f \mapsto \sup_{0 \leq t \leq 1} f(t)$ [resp. $f \mapsto \inf_{0 \leq t \leq 1} f(t)$] is continuous.

(ii) Let $\mathfrak{r}'$ be picked at random in $S(M) \setminus \{\mathfrak{r}\}$ conditionally on $M$ under $P_{\mathbf{q}}(\cdot | \#F(M) = n)$. Then by Proposition 6, the law of $d_M(\mathfrak{r}, \mathfrak{r}')$ is the same as that of $1 + L(V) - \min_{v' \in T} L(v')$ under $\mathbb{P}^{(0)}(\cdot | \#T^{(1)} = n)$, where $V$ is uniformly picked among the vertices of $T^{(0)}$ conditionally on $T, L$. To be completely rigorous, this involves an enlarging of the probability space $\mathbb{T}$, and we do it in the following convenient way. We endow the space $\mathbb{T}^{(0)} \times [0,1]$ with the law $\widetilde{\mathbb{P}} = \mathbb{P}^{(0)}(\cdot | \#T^{(1)} = n) \times dx$, where $dx$ is Lebesgue measure on $[0, 1]$. If $((T, L), U)$ is the identity map on this space, then under $\widetilde{\mathbb{P}}$, $U$ is a uniform random variable in $[0, 1]$, independent of $H^T, S^{T,L}, \overline{J}_T^{(0)}$. Then, let $U_T^{(0)} = (\overline{J}_T^{(0)})^{-1}(U)$, where $(\overline{J}_T^{(0)})^{-1}$ is the right-continuous inverse of $\overline{J}_T^{(0)}$. By definition of $\overline{J}_T^{(0)}$, it holds that $(\#T-1)U_T^{(0)}$ is the rank in depth-first order of a uniform random vertex of $T^{(0)}$. Otherwise said, $u((\#T-1)U_T^{(0)})$ is uniform in $T^{(0)}$ given $(T, L)$. Hence, the law of the distance in a $P_{\mathbf{q}}(\cdot | \#F(M) = n)$-chosen random map from the root to a uniformly chosen nonroot vertex is the same as that of $S_{(\#T-1)U_T^{(0)}}^{T,L} - \min S^{T,L} + 1$ under $\widetilde{\mathbb{P}}$.

On the other hand, the convergence of $(\overline{J}_T^{(0)}(t), 0 \leq t \leq 1)$ to the deterministic identity function, which is described in Lemma 9, must hold jointly with that of the height and snake processes under $\mathbb{P}^{(0)}(\cdot | \#T^{(1)} = n)$. So, by Skorokhod's representation theorem, we may find a probability space on which the convergence holds almost surely, that is, we can find processes $(H^n, S^n, \overline{J}^n)$ with the same law as

(22) $(n^{-1/2} H_{(\#T-1)\cdot}^T, n^{-1/4} S_{(\#T-1)\cdot}^{T,L}, \overline{J}_T^{(0)})$ under $\mathbb{P}^{(0)}(\cdot | \#T^1 = n),$

and which converge uniformly a.s. to a triple $(B, S, \mathrm{id}_{[0,1]})$ where $(B, S)$ is distributed as $(D_{\mathbf{q}}\mathrm{e}, C_{\mathbf{q}}\mathrm{r})$ under $\mathbb{N}^{(1)}$. We take a uniform random variable $\widetilde{U}$ on $[0, 1]$, independent of all these processes, and let $U^n = (\overline{J}^n)^{-1}(\widetilde{U})$, which has the same law as $U_T^{(0)}$ with the above notation. Since $\overline{J}^n$ converges uniformly to the identity, $U_n$ converges to $\widetilde{U}$ a.s., and therefore, $S_{U^n}^n$ converges a.s. to $S_{\widetilde{U}}$ as $n \to \infty$. Moreover, $\inf S^n$ converges to $\inf S$ as $n \to \infty$, so finally, we obtain that $n^{-1/4}(S_{(\#T-1)U_T^{(0)}}^{T,L} - \min S^{T,L} + 1)$ under $\widetilde{\mathbb{P}}$ converges in distribution to $S_{\widetilde{U}} - \inf S$, which has the law of $C_{\mathbf{q}}(\mathrm{r}_U - \inf \mathrm{r})$ under $\mathbb{N}^{(1)} \times dx$. By the rerooting properties of the Brownian snake of [17, 20], $(\mathrm{r}_{s+t \bmod 1} - \mathrm{r}_t, 0 \leq s \leq 1)$ has same law as r under $\mathbb{N}^{(1)}$, for every $t$. So under

$\mathbb{N}^{(1)} \times dx$,

$$r_U - \inf_{0 \leq s \leq 1} r_s = r_U - \inf_{0 \leq s \leq 1}(r_{s+U \bmod 1} - r_U + r_U)$$

$$\stackrel{(d)}{=} - \inf_{0 \leq s \leq 1} r_s = -\Delta_-(r),$$

which by symmetry has the same law as $\Delta_+(r)$, as claimed.

(iii) For every $k > 0$, we have, using (10),

$$\mathscr{I}^{\mathbf{m}}(k) = (1 + \#\mathbf{t}^{(0)})^{-1}(\#\{u \in \mathbf{t}^{(0)} : \ell(u) - \min \ell + 1 = k\})$$

whenever $(\mathbf{t}, \ell) = \Psi(\mathbf{m})$. We can rewrite this as

$$\frac{\#\mathbf{t}^{(0)}}{1 + \#\mathbf{t}^{(0)}} \int_0^1 \mathbb{1}_{\{S^{\mathbf{t},\ell}_{(\#T-1)t} - \inf S^{\mathbf{t},\ell} + 1 = k\}} \, d\overline{J}^{(0)}_{\mathbf{t}}(t),$$

Thus, for every $g$ which is Lipschitz and bounded,

$$\langle \mathscr{I}^{\mathbf{m}}_n, g \rangle = \frac{\#\mathbf{t}^{(0)}}{1 + \#\mathbf{t}^{(0)}} \int_0^1 g\left(\frac{S^{\mathbf{t},\ell}_{(\#T-1)t} - \inf S^{\mathbf{t},\ell} + 1}{n^{1/4}}\right) d\overline{J}^{(0)}_{\mathbf{t}}(t) + \frac{g(0)}{1 + \#\mathbf{t}^{(0)}}.$$

Note that because of the convergence of $\overline{J}^{(0)}_T$ in Lemma 9, the quantity $\#T^{(0)}$ under $P^{(0)}(\cdot | \#T^{(1)} = n)$ converges to infinity in probability. We then use again the Skorokhod representation theorem, and suppose given processes $H^n, S^n, \overline{J}^n$ with respective laws that of (22), which converge almost surely for the uniform norm to $(B, S, \mathrm{id}_{[0,1]})$, where $(B, S)$ is distributed as $(D_{\mathbf{q}}e, C_{\mathbf{q}}r)$ under $\mathbb{N}^{(1)}$. Then, the measures $d\overline{J}^n$ converge weakly to the uniform law on $[0, 1]$. We have

$$\int_0^1 g(S^n_t - \inf S^n + n^{-1/4}) \, d\overline{J}^n(t) - \int_0^1 g(S_t - \inf S) \, dt$$

$$= \int_0^1 (g(S^n_t - \inf S^n + n^{-1/4}) - g(S_t - \inf S)) \, d\overline{J}^n(t)$$

$$+ \int_0^1 g(S_t - \inf S)(d\overline{J}^n(t) - dt).$$

The first term on the right-hand side converges to 0, because $S^n - \inf S^n$ converges uniformly to $S - \inf S$, and $g$ is Lipschitz. The second term converges to 0 because $g(S_t - \inf S), 0 \leq t \leq 1$, is continuous and bounded, and $d\overline{J}^n$ converges weakly to $dt$. Since $\int_0^1 g(S_t - \inf S) \, dt$ has the law of $\langle \mathscr{I}^{\mathrm{r}}_{\mathbf{q}}, g \rangle$ under $\mathbb{N}^{(1)}$, this ends the proof. $\square$

REMARK. A somewhat simpler proof for (ii), using (iii), could be obtained following the same lines as Le Gall [16]. We thought however that the present approach, which for example can be easily extended to handle the case of several sampled points, was worth mentioning.

The proof of Proposition 3 is entirely similar to the previous proof, the only significant difference being that one should use the probability distributions $P^{(0)}(\cdot|\#T^{(0)}=n)$ rather than $P^{(0)}(\cdot|\#T^{(1)}=n)$. This tacitly implies that we must take $n$ along values for which this conditioning is well defined. Except from that, there is a minor change due to the fact that the scaling constants in the limit are different, namely

$$\frac{\sqrt{2}\Sigma(1+m_0)^{1/4}}{\sqrt{\sigma}} = \left(\frac{4\rho_{\mathbf{q}}}{9}\right)^{1/4}.$$

Details are left to the reader.

**4. Convergence of the height process.** The goal of this section is to prove the convergence of the first component in Theorem 8. This involves a couple of lemmas, which we now describe.

4.1. *GW forests.* A forest $\mathbf{f}$ is a subset of $\mathcal{U}$ that is of the form

$$\mathbf{f} = \bigcup_k k\mathbf{t}_{(k)},$$

where $(\mathbf{t}_{(k)})$ is a finite of infinite sequence of trees, called the tree components of $\mathbf{f}$. We let $\mathcal{F}$ be the set of forests. If $\mathbf{f} \in \mathcal{F}$ and $u \in \mathbf{f}$, we define the fringe subtree $\mathbf{f}_u \in \mathcal{T}$ by $\{v \in \mathcal{U} : uv \in \mathbf{f}\}$ as above, and $[\mathbf{f}]_u = \{u\} \cup (\mathbf{f} \setminus u\mathbf{f}_u) \in \mathcal{F}$ the pruned forest. With this notation, observe that the tree components of $\mathbf{f}$ are $\mathbf{f}_1, \mathbf{f}_2, \ldots$. For $\mathbf{f} \in \mathcal{F}$ and $u \in \mathbf{f}$, we let $c_{\mathbf{f}}(u) = c_{\mathbf{f}_u}(\varnothing)$ be the number of children of $u$ in $\mathbf{f}$. If it is understood that $u$ is an element of $\mathbf{f}$, for $\mathbf{f} \in \mathcal{F}$, we call $|u|-1$ the height of $u$. It differs from the convention on trees because we want the roots of the forest components to be at height 0. For $\mathbf{f} \in \mathcal{F}$ and $u \in \mathbf{f}$, $\Upsilon_{\mathbf{f}}(u)$ be the first letter of $u$, that is, the rank of the tree component of $\mathbf{f}$ containing $u$.

We also want to consider forests of marked trees, that is, sets of the form $\bigcup_k k\mathbf{t}_{(k)}$ with $\mathbf{t}_{(k)} \in \mathcal{T}_i, i \in \{0,1\}$. We then define, for $i \in \{0,1\}$,

$$\mathbf{f}^{(i)} = \bigcup_k k\mathbf{f}^{(i)}_k.$$

We let $\mathcal{F}_0$ be the set of forests constituted only of trees marked 0 at the root, and $\mathcal{F}_1$ the set of forests constituted only of trees marked 1 at the root.

If $i \in \{0,1\}$, and $(\mu_0, \mu_1)$ is a (sub) critical pair of offspring distributions as in Section 2.2, and for $r \in \mathbb{N} \sqcup \{\infty\}$, we let $P^{(i)}_r$ be the image law on $\mathcal{F}_i$ of $(P^{(i)}_\mu)^{\otimes r}$ under the map

$$(\mathbf{t}_{(1)}, \mathbf{t}_{(2)}, \ldots) \mapsto \bigcup_k k\mathbf{t}_{(k)},$$

going from the set of sequences of $r$ trees in $\mathcal{T}_i$ to $\mathcal{F}_i$. We do not refer to $\mu$ in the definition of $P_r^{(i)}$, but the value of $\mu_0, \mu_1$ should be clear according to the context. We let $F:\mathcal{F} \to \mathcal{F}$ be the identity mapping.

In the sequel, if $\mathbf{t} \in \mathcal{T}$ or $\mathbf{f} \in \mathcal{F}$, we let $u(0) \prec u(1) \prec \cdots$ be the list of vertices of $\mathbf{t}$ or $\mathbf{f}$ in depth-first order. Similarly, for $i, j \in \{0, 1\}$ and $\mathbf{t} \in \mathcal{T}_i$ or $\mathbf{f} \in \mathcal{F}_i$, we let $u^{(j)}(0) \prec u^{(j)}(1) \prec \cdots$ be the list of vertices of $\mathbf{t}^{(j)}$ or $\mathbf{f}^{(j)}$ listed in depth-first order. Although there is no mention of $\mathbf{t}, \mathbf{f}$ in the notation, it should be unambiguous according to the context.

4.2. *Controlling the height and number of components of forests.* The first technical lemma gives an exponential control on quantities related to the $n$ first vertices in a monotype GW forest.

LEMMA 10. *Let $\mu$ be a critical nondegenerate offspring distribution on $\mathbf{Z}_+$, that is, $\mu(1) < 1$ and $\mu$ has mean 1. Suppose also that $\mu$ has finite variance. Let $P_\infty$ be the law of a GW forest with and infinite number of components and offspring distribution $\mu$ (with the previous notations it is $P_\infty^{(i)}$ whenever $\mu_0 = \mu_1 = \mu$, in that case, the role of $i \in \{0, 1\}$ is irrelevant). Then, there exist constants $0 < C_1, C_2 < \infty$ such that for every $\eta > 0$, for every $n \geq 0$,*

$$(23) \qquad P_\infty\left(\max_{0 \leq k \leq n} |u(k)| \geq n^{1/2+\eta}\right) \leq C_1(n+1)\exp(-C_2 n^\eta)$$

*and*

$$(24) \qquad P_\infty(\Upsilon_F(u(n)) \geq n^{1/2+\eta}) \leq C_1 \exp(-C_2 n^\eta).$$

PROOF. We bound the first probability by $(n+1)\max_{0 \leq k \leq n} P_\infty(|u(k)| \geq n^{1/2+\eta})$. It is known ([10], Section 2.2) that $|u(k)| - 1$ has same distribution as the number of weak records for a random walk with step distribution $\mu(\cdot+1)$ on $\{-1\} \cup \mathbf{Z}_+$, from time 1 up to time $k$. Suppose such a random walk $(W_n, n \geq 0)$ is defined on some probability space $(\widetilde{\Omega}, \widetilde{\mathcal{A}}, \widetilde{P})$. By assumption, the step distribution of this random walk is centered and has finite variance. Therefore, calling $\tau_0 = 0$ and $\tau_i, i \geq 1$, the time of the $i$th weak record of $(W_n, n \geq 0)$, we have from [11] that $(\tau_i - \tau_{i-1}, i \geq 1)$ is i.i.d., and the Laplace exponent of the common distribution satisfies

$$(25) \qquad \widetilde{\phi}(\lambda) = -\log \widetilde{E}[\exp(-\lambda \tau_1)] \underset{\lambda \downarrow 0}{\sim} C'\sqrt{\lambda}$$

for some $C' > 0$. Now, for $k \leq n$, we write $P_\infty(|u(k)| - 1 \geq m)$ as

$$(26) \qquad \widetilde{P}\left(\sum_{i=1}^m (\tau_i - \tau_{i-1}) \leq k\right) \leq e\mathbb{E}\left[\exp\left(-\sum_{i=1}^m \frac{\tau_i - \tau_{i-1}}{k}\right)\right]$$
$$= \exp(1 - m\widetilde{\phi}(1/k)),$$

which by monotonicity of $\widetilde{\phi}$ is less than $\exp(1 - m\widetilde{\phi}(1/n))$, and taking $m = \lceil n^{1/2+\eta} \rceil - 1$ and using (25) gives (23) for large enough $n$, thus for every $n$ up to tuning the constants $C_1, C_2$.

The proof of (24) is very similar. By a well-known application of the Otter–Dwass formula (see, e.g., [22], Chapter 5), the sizes $(\#F_1, \#F_2, \ldots)$ of the components of the forest $F$ under $P_\infty$ are i.i.d. random variables with distribution

$$P_\infty(\#F_1 = n) = \frac{1}{n}\widetilde{P}(W_n = -1).$$

By using once again the fact that the step distribution is centered and has finite variance, the local limit theorem ([11], Theorem XV.5.3) entails that $P_\infty(\#F_1 = n)$ is equivalent to $C'' n^{-3/2}$ (or 0 for lattice type reasons), where $C'' > 0$ is some positive constant. Therefore, an Abelian theorem ([11], Theorem XIII.5.5) entails that the Laplace exponent $\phi$ of the distribution of $\#F_1$ under $P_\infty$ is equivalent to $\lambda^{1/2}$ up to a multiplicative constant as $\lambda \downarrow 0$. Noticing that $\{\Upsilon_F(u(n)) \geq m\} = \{\sum_{i=1}^{m-1} \#F_i \leq n\}$, the result is then obtained by a straightforward analog of (26), replacing $\widetilde{P}$ by $P_\infty$, $(\tau_i - \tau_{i-1})$ by $\#F_i$, and $\widetilde{\phi}$ by $\phi$. We finally adapt the constants $C_1$, $C_2$ so that they match to both cases. $\square$

A consequence of this is the following analogous result for two-type forests.

LEMMA 11. *Let $(\mu_0, \mu_1)$ be a critical nondegenerate offspring distribution, and suppose $\mu_0, \mu_1$ have finite variances. Then, for $i, j \in \{0, 1\}$, and every $\eta > 0$, there exists some $\varepsilon > 0$ such that for every $n$ large enough,*

$$(27) \qquad P_\infty^{(i)}\left(\max_{0 \leq k \leq n} |u^{(j)}(k)| \geq n^{1/2+\eta}\right) \leq \exp(-n^\varepsilon).$$

*Moreover,*

$$(28) \qquad P_\infty^{(i)}(\Upsilon_F(u^{(j)}(n)) \geq n^{1/2+\eta}) \leq \exp(-n^\varepsilon).$$

REMARK. Since $\max_{u \preceq u^{(j)}(n)} |u| \leq \max_{0 \leq k \leq n} |u^{(j)}(k)| + 1$, (27) also yields

$$(29) \qquad P_\infty^{(i)}\left(\max_{u \preceq u^{(j)}(n)} |u| \geq n^{1/2+\eta}\right) \leq \exp(-n^\varepsilon).$$

PROOF OF LEMMA 11. Suppose that $i = j$. For $\mathbf{f} \in \mathcal{F}_i$, define the following analog of the transformation $\Gamma_\mathbf{f} : \mathbf{f} \to \mathcal{U}$ of Section 2.2, by

$$\Gamma_\mathbf{f}(kv) = k\Gamma_{\mathbf{f}_k}(v), \qquad kv \in \mathbf{f}^{(i)},$$

and $\Gamma_{\mathbf{f}}(v) = \Gamma_{\mathbf{f}}(\neg v)$ if $v \in \mathbf{f}^{(i+1)}$, so $\Gamma$ skips odd generations in all of the tree components of the forest $\mathbf{f}$. We let $\Gamma(\mathbf{f})$ be its image, so that

$$\Gamma(\mathbf{f}) = \bigcup_{k \geq 1} k\Gamma(\mathbf{f}_k) \in \mathcal{F}.$$

Notice that $\Gamma_{\mathbf{f}}(u^{(i)}(n))$ is the $(n+1)$st vertex in depth-first order in $\Gamma(\mathbf{f})$. It is then a consequence of the definitions of $\Gamma_{\mathbf{t}}, \Gamma_{\mathbf{f}}$ that

$$2|\Gamma_{\mathbf{f}}(u^{(i)}(n))| - 2 = |u^{(i)}(n)| - 1 \quad \text{and}$$
(30)
$$\Upsilon_{\mathbf{f}}(u^{(i)}(n)) = \Upsilon_{\Gamma(\mathbf{f})}(\Gamma_{\mathbf{f}}(u^{(i)}(n))).$$

As in the discussion leading to Proposition 5, under $P^{(i)}$, the tree $\Gamma(T)$ is a monotype critical GW tree, and therefore, under $P_\infty^{(i)}$, the forest $\Gamma(F)$ is a monotype critical GW forest. Its offspring distribution $\overline{\mu}$ has generating function $\mathbf{G}_{\overline{\mu}} = \mathbf{G}_{\mu_i} \circ \mathbf{G}_{\mu_{i+1}}$, and if $\mu_0, \mu_1$ have finite variances, we can differentiate $\mathbf{G}_{\overline{\mu}}$ twice to obtain that $\overline{\mu}$ itself has finite variance. Hence, Lemma 10 applies and gives that for every $\eta > 0$, there is some $\varepsilon > 0$ such that for $n$ large enough,

(31) $$P_\infty^{(i)}\left( \max_{0 \leq k \leq n} |\Gamma_F(u^{(i)}(k))| \geq n^{1/2 + \eta} \right) \leq \exp(-n^\varepsilon).$$

Therefore, using (30), up to taking a smaller $\varepsilon$, we obtain (27) for some $\varepsilon > 0$ and all $n$ large. Similarly, (28) follows by applying (30) and (24) to $\Gamma(F)$ under $P_\infty^{(i)}$.

It remains to prove the case $i + 1 = j$. To this end, we introduce a transformation on forests that skips the first generation. For $\mathbf{f} \in \mathcal{F}$, we let $\pi(\mathbf{f})$ be the forest with tree components $(\mathbf{f}_{k1}, \ldots, \mathbf{f}_{kc_{\mathbf{f}}(k)}, k \geq 1)$, where these tree components are put in lexicographical order of the index $kl, 1 \leq l \leq c_{\mathbf{f}}(k)$. For every $u \in \mathbf{f}$ with $|u| \geq 2$, there is a unique corresponding $u'$ in $\pi(\mathbf{f})$, that we denote $\pi_{\mathbf{f}}(u) = u'$.

If $\mathbf{f} \in \mathcal{F}_i$, then $\pi(\mathbf{f})$ is considered as an element of $\mathcal{F}_{i+1}$, and we let $\Gamma'_{\mathbf{f}} = \Gamma_{\pi(\mathbf{f})} \circ \pi_{\mathbf{f}}$, and $\Gamma'(\mathbf{f}) = \Gamma'_{\mathbf{f}}(\mathbf{f})$, so $\Gamma'$ first skips the first generation of a forest, and then skips all odd generations of the new forest (see Figure 2 where $\Gamma'$ is applied to a forest with one component). As for $\Gamma$, is is easy to see that if $\mathbf{f} \in \mathcal{F}_i$, then $\Gamma'_{\mathbf{f}}(u^{(i+1)}(k))$ is the $(k+1)$st vertex of $\Gamma'(\mathbf{f})$ in depth-first order, and has height satisfying $2|\Gamma'_{\mathbf{f}}(u^{(i+1)}(k))| = |u^{(i+1)}(k)|$, which mirrors the first half of (30). Moreover, we have

(32) $$\Upsilon_{\Gamma'(\mathbf{f})}(\Gamma'_{\mathbf{f}}(u^{(i+1)}(n)))$$
$$\geq \Upsilon_{\mathbf{f}}(u^{(i+1)}(n)) - \#\{1 \leq k \leq \Upsilon_{\mathbf{f}}(u^{(i+1)}(n)) : c_{\mathbf{f}}(k) = 0\}.$$

Indeed, some of the tree components of $\mathbf{f}$ do not have vertices of type $i+1$, and thus do not count in the construction of $\Gamma'(\mathbf{f})$. This gives the second term in the right-hand side. On the other hand, those tree components which have at least one vertex of type $i+1$ produce at least one tree component in $\Gamma'(\mathbf{f})$, which gives the inequality.

Now, by construction, the law of $\pi(F)$ under $P_\infty^{(i)}$ is $P_\infty^{(i+1)}$, so that the law of $\Gamma'(F)$ under $P_\infty^{(i)}$ is that of a monotype GW forest with offspring distribution $\overline{\mu}'$, and $\mathbf{G}_{\overline{\mu}'} = \mathbf{G}_{\mu_{i+1}} \circ \mathbf{G}_{\mu_i}$. We thus obtain (27) by applying Lemma 10 to $\Gamma'(F)$ under $P_\infty^{(i)}$, just as we did in the case $i=j$.

Obtaining (28) is slightly more delicate, since the second half of (30) is now replaced by (32). Under $P_\infty^{(i)}$, the random variables $c_F(k)$ are identically distributed with law $\mu_i$, so that $B(m) := \#\{1 \leq k \leq m : c_F(k) = 0\}$ is a Binomial random variable with parameters $(m, \mu_i(0))$. By Hoeffding's inequality, if $B(n,p)$ is Binomial with parameters $n \geq 1, 0 < p < 1$, we have

$$(33) \qquad P(|B(n,p) - np| \geq y) \leq 2\exp(-2y^2/n).$$

Therefore, for any $\gamma > 0$

$$(34) \qquad P_\infty^{(i)}\left(\max_{1 \leq k \leq n} |B(k) - k\mu_i(0)| \geq n^{1/2+\gamma}/2\right) \leq 2n\exp(-n^{2\gamma}/2).$$

For simplicity, write $\Upsilon_n = \Upsilon_F(u^{(i+1)}(n))$ and $\Upsilon'_n = \Upsilon_{\Gamma'(F)}(\Gamma'_F(u^{(i+1)}(n)))$. Then, (32) gives $B(\Upsilon_n) \geq \Upsilon_n - \Upsilon'_n$, and since $\Upsilon_n \leq n$ by definition, (34) yields

$$P_\infty^{(i)}(\Upsilon_n - \Upsilon'_n \geq n^{1/2+\eta}/2) \leq 2n\exp(-n^{2\eta}/2).$$

Finally,

$$P_\infty^{(i)}(\Upsilon_n \geq n^{1/2+\eta}) \leq P_\infty^{(i)}(\Upsilon_n - \Upsilon'_n \geq n^{1/2+\eta}/2) + P_\infty^{(i)}(\Upsilon'_n \geq n^{1/2+\eta}/2),$$

and both terms are $\leq \exp(-n^\varepsilon)$ for some $\varepsilon > 0$ and $n$ large enough, the second term because of (24) applied to the monotype GW forest $\Gamma'(F)$ under $P_\infty^{(i)}$. □

4.3. *Ancestral decomposition of a GW forest.* A key result for our study is a multitype version of an ancestral decomposition for GW trees, related to the so-called size-biased GW distribution. It is inspired from [14]. Let $(\mu_0, \mu_1)$ be a nondegenerate critical two-type offspring distribution, and define the associated size-biased distributions

$$\hat{\mu}_0(k) = \frac{k\mu_0(k)}{m_0}, \qquad \hat{\mu}_1(k) = \frac{k\mu_1(k)}{m_1}, \qquad k \geq 0.$$

Notice that these distributions do not charge $\{0\}$. The size-biased GW tree is an infinite tree (an element of $\widehat{\mathcal{T}}$ with the notation of Section 2.2) containing a unique *spine*, that is, an infinite injective path starting from the root. On some probability space, let $(X_u, \widehat{X}_u, j_u, u \in \mathcal{U})$ be a sequence of i.i.d. random variables such that $X_u$ has law $\mu_{|u|}$ ($|u|$ taken modulo 2), $\widehat{X}_u$ has law $\hat{\mu}_{|u|}$, and conditionally on $\widehat{X}_u$, $j_u$ is uniform in $\{1, 2, \ldots, \widehat{X}_u\}$. Then, let $w_0 = \varnothing$, and recursively $w_{j+1} = w_j j_{w_j}$, $j \geq 0$. Let $\widetilde{X}_u = \widehat{X}_u$ if $u \in \{w_0, w_1, \ldots\}$, and $\widetilde{X}_u = X_u$ otherwise. Finally, let $\widehat{\xi}$ be the element of $\widehat{\mathcal{T}}_0$ defined by

$$\widehat{\xi} = \{u = u_1 \cdots u_k \in \mathcal{U} : u_i \leq \widetilde{X}_{u_1 \cdots u_{i-1}}, 1 \leq i \leq k\} \cup \{\varnothing\}.$$

We see that $\widehat{\xi}$ is almost a GW tree, except for one distinguished spine which uses the distributions $\hat{\mu}_0, \hat{\mu}_1$ instead of $\mu_0, \mu_1$. In particular, under the criticality assumption, we see that all fringe subtrees of $\widehat{\xi}$ attached to the spine, that is, of the form $\widehat{\xi}_{w_j k}$ for $k \neq j_{w_j}$, are a.s. finite, so the only infinite simple path starting from the root in $\widehat{\xi}$ is a.s. $(w_0, w_1, \ldots)$. In particular the trees $[\widehat{\xi}]_{w_h}$ are a.s. finite for $h \geq 0$. For every $h \geq 0$, let $\widehat{P}^{(0),h}$ be the law of $([\widehat{\xi}]_{w_h}, w_h)$, where we understand that $[\widehat{\xi}]_{w_h}$ is an element of $\mathcal{T}_0$. It is a law on the set of pointed trees with white root whose distinguished vertex is a leaf (i.e., has no child)

$$\mathcal{T}_0^* = \{(\mathbf{t}, u) : \mathbf{t} \in \mathcal{T}_0, u \in \mathbf{t}, c_\mathbf{t}(u) = 0\}.$$

We let $(T, V)$ be the identity mapping on $\mathcal{T}_0^*$. Similarly we define $\widehat{P}^{(1),h}$ on $\mathcal{T}_1^*$, where $\mathcal{T}_1^*$ is a copy of $\mathcal{T}_0^*$, by switching the roles of $\mu_0$ and $\mu_1$. Finally, for $r \in \mathbb{N}$, $j \in \{1, \ldots, r\}$, $i \in \{0, 1\}$ and $h \geq 0$, we let $\widehat{P}_r^{(i),j,h}$ be the law on

$$\mathcal{F}_i^* = \{(\mathbf{f}, u) : \mathbf{f} \in \mathcal{F}_i, u \in \mathbf{f}, c_\mathbf{f}(u) = 0\}$$

of the random variable

$$\left(\bigcup_{1 \leq k \leq r} k \xi_{(k)}, jw\right),$$

where $(\xi_{(k)}, k \neq j)$ are independent with distribution $P^{(i)}$ and independent of $(\xi_{(j)}, w)$, which has law $\widehat{P}^{(i),h}$. We let $(F, V)$ be the identity mapping on $\mathcal{F}_i^*$.

LEMMA 12 (Ancestral decomposition for GW forests). *Let $(\mu_0, \mu_1)$ be a critical nondegenerate offspring distribution. For every $r \in \mathbb{N}$ and nonnegative measurable functions $G_1, G_2$*

$$(35) \quad E_r^{(i)}\left[\sum_{w \in F} G_1(w, [F]_w) G_2(F_w)\right]$$
$$= \frac{1}{1 + m_{i+1}} \sum_{j=1}^{r} \sum_{h \geq 0} (1 + m_{h+1+i}) \widehat{E}_r^{(i),j,h}[G_1(V, F)] E^{(h+i)}[G_2(T)]$$

*where as usual $i + 1, h + i(+1)$ are taken modulo 2.*

PROOF. We treat the case $i = 0$ only. Let $\mathbf{f} \in \mathcal{F}_0$, let $u$ be a leaf of $\mathbf{f}$ and let $\mathbf{t} \in \mathcal{T}_{|u|-1}$. Then, it is enough to show the result for $G_1 = \mathbb{1}_{\{(u,\mathbf{f})\}}$ and $G_2 = \mathbb{1}_{\mathbf{t}}$, by linearity and monotone convergence. In this case, the left-hand side of (35) is equal to $P_r^{(0)}(F = [\mathbf{f}, u, \mathbf{t}])$, where $[\mathbf{f}, u, \mathbf{t}]$ is the only forest $\mathbf{f}' \in \mathcal{F}_0$ containing $u$ with $[\mathbf{f}']_u = \mathbf{f}$ and $\mathbf{f}'_u = \mathbf{t}$. This probability is

$$\prod_{v \in \mathbf{f}'} \mu_{|v|-1}(c_{\mathbf{f}'}(v)).$$

Let $j = u_1$ be the first letter of $u$. We can redisplay the last expression as

$$\prod_{v \in \mathbf{t}} \mu_{|v|+|u|}(c_{\mathbf{t}}(v)) \prod_{1 \le l \le r, l \ne j} \prod_{v \in \mathbf{f}_l} \mu_{|v|}(c_{\mathbf{f}}(lv))$$
$$\times \prod_{v \in \mathbf{f}_j, jv \nvdash u} \mu_{|v|}(c_{\mathbf{f}}(jv)) \prod_{v \in \mathbf{f}, v \vdash u, v \ne u} \mu_{|v|-1}(c_{\mathbf{f}}(v))$$

(we omit the brackets around the different products for convenience). Let $S = \{v \in \mathbf{f} : \neg v \vdash u, v \nvdash u\}$ be the set of neighbors of the ancestors of $u$, which are not ancestors of $u$. We recognize

$$P_r^{(0)}(F = [\mathbf{f}, u, \mathbf{t}]) = P^{(|u|-1)}(T = \mathbf{t}) \prod_{l \ne j, 1 \le l \le r} P^{(0)}(T = \mathbf{f}_l)$$
$$\times \prod_{v \in S} P^{(|v|-1)}(T = \mathbf{f}_v) \prod_{v \in \mathbf{f}, v \vdash u, v \ne u} \mu_{|v|-1}(c_{\mathbf{f}}(v)).$$

We can also rewrite the last product as

$$\prod_{v \in \mathbf{f}, v \vdash u, v \ne u} \widehat{\mu}_{|v|-1}(c_{\mathbf{f}}(v)) \frac{m_{|v|-1}}{c_{\mathbf{f}}(v)}.$$

After a moment's thought, we see that, letting $u = ju'$,

$$\prod_{v \in S} P^{(|v|-1)}(T = \mathbf{f}_v) \prod_{v \in \mathbf{f}, v \vdash u, v \ne u} \frac{\widehat{\mu}_{|v|-1}(c_{\mathbf{f}}(v))}{c_{\mathbf{f}}(v)}$$
$$= \widehat{P}^{(0),|u|-1}((T, V) = (\mathbf{f}_j, u')).$$

On the other hand, as one can check from the fact that $m_0 m_1 = 1$,

$$\prod_{v \in \mathbf{f}, v \vdash u, v \ne u} m_{|v|-1} = \frac{1 + m_{|u|}}{1 + m_1}$$

and we finally recognize

$$P_r^{(0)}(F = [\mathbf{f}, u, \mathbf{t}]) = \frac{1 + m_{|u|}}{1 + m_1} \widehat{P}^{(0),|u|-1}((T, V) = (\mathbf{f}_j, u'))$$
$$\times \left( \prod_{l \ne j, 1 \le l \le r} P^{(0)}(T = \mathbf{f}_l) \right) P^{(|u|-1)}(T = \mathbf{t})$$

$$= \frac{1+m_{h+1}}{1+m_1}\widehat{E}_r^{(0),j,h}\,[\mathbb{1}_{\{u,\mathbf{f}\}}(V,F)]E^{(h)}[\mathbb{1}_{\{\mathbf{t}\}}(T)],$$

where $h=|u|-1$, which is (35). $\square$

The first corollary we infer from this is a control on the maximum vertex degree in a two-type GW forest.

LEMMA 13. *Assume that the pair $(\mu_0,\mu_1)$ is nondegenerate, critical and has some exponential moments. Then for every $\eta>0$ there exists $\varepsilon>0$ such that for $n$ large enough, and $i,j\in\{0,1\}$,*

$$(36)\qquad P_\infty^{(i)}\left(\max_{u\preceq u^{(j)}(n)} c_F(u)\geq n^\eta\right)\leq \exp(-n^\varepsilon).$$

PROOF. Let $\eta>0$. Then by using Lemma 11, the left-hand side of (36) is equal to

$$(37)\qquad \begin{aligned}&P_\infty^{(i)}\Bigg(\max_{u\preceq u^{(j)}(n)} c_F(u)\geq n^\eta,\\ &\qquad \max_{u\preceq u^{(j)}(n)}|u|\leq n^{1/2+\eta},\Upsilon_F(u^{(j)}(n))\leq n^{1/2+\eta}\Bigg)+R(n)\\ &\leq P_{[n^{1/2+\eta}]}^{(i)}\left(\max_{|u|\leq n^{1/2+\eta}} c_F(u)\geq n^\eta\right)+R(n),\end{aligned}$$

where $R(n)\leq \exp(-n^\varepsilon)$ for some $\varepsilon>0$ and $n$ large enough. But then,

$$P_{[n^{1/2+\eta}]}^{(i)}\left(\max_{|u|\leq n^{1/2+\eta}} c_F(u)\geq n^\eta\right)\leq E_{[n^{1/2+\eta}]}^{(i)}\left(\sum_{u\in F}\mathbb{1}_{\{|u|\leq n^{1/2+\eta}\}}\mathbb{1}_{\{c_F(u)\geq n^\eta\}}\right),$$

and applying Lemma 12 to $G_1(u,\mathbf{f})=\mathbb{1}_{\{|u|\leq n^{1/2+\eta}\}}$ and $G_2(\mathbf{t})=\mathbb{1}_{\{c_\mathbf{t}(\varnothing)\geq n^\eta\}}$, this is equal to

$$(38)\qquad \begin{aligned}&\frac{1}{1+m_{i+1}}\sum_{j=1}^{[n^{1/2+\eta}]}\sum_{h=0}^{[n^{1/2+\eta}]}(1+m_{h+i+1})\mu_{h+i}([n^\eta,\infty))\\ &\leq Cn^{1+2\eta}(\mu_0([n^\eta,\infty))\vee\mu_1([n^\eta,\infty))),\end{aligned}$$

for some $C>0$. But since $\mu_0,\mu_1$ have some exponential moments, it holds that $\mu_i([n^\eta,\infty))\leq \exp(-an^\eta)$, for some $a>0$, and $n$ large enough. Combined with (38) and (37), this yields (36). $\square$

4.4. *An estimate for the size of GW trees.* In order to pass from statements on forests to statements on conditioned trees, we need to estimate the number of vertices of either type in two-type GW trees.

LEMMA 14. *Let $(\mu_0, \mu_1)$ be a critical nondegenerate offspring distribution, and suppose that $\mu_0$ and $\mu_1$ have finite variances. Then for $i, j \in \{0, 1\}$, there exists a finite constant $C_{ij} > 0$ such that*

$$n^{3/2} P^{(i)}(\#T^{(j)} = n) \longrightarrow C_{ij},$$

*where it is understood that $n$ goes to infinity along values for which the quantity on the left-hand side is strictly positive.*

PROOF. Suppose $i = j = 0$. Then $P^{(0)}(\#T^{(0)} = n) = P^{(0)}(\#\Gamma(T) = n)$, where $\Gamma$ denotes the mapping that skips odd generations, as usual, so $\Gamma(T)$ under $P^{(0)}$ is a (monotype) GW tree whose offspring distribution $\overline{\mu}$ has generating function $\mathbf{G}_{\overline{\mu}} = \mathbf{G}_{\mu_0} \circ \mathbf{G}_{\mu_1}$ It results that $\overline{\mu}$ is critical, nondegenerate, and has finite variance (by differentiating twice $\mathbf{G}_{\overline{\mu}}$). The conclusion follows from the Otter–Dwass formula and the local limit theorem that we used in the proof of Lemma 11: we have $P^{(0)}(\#T^{(0)} = n) = n^{-1}\widetilde{P}(W_n = -1) \sim C_{00} n^{-3/2}$ where $W$ under $\widetilde{P}$ is a random walk with step distribution $\overline{\mu}(\cdot + 1)$. The case $i = j = 1$ is similar.

It remains to deal with the case $i = 0, j = 1$. In that case, we have

$$P^{(0)}(\#T^{(1)} = n) = \sum_{r \geq 1} P^{(0)}(c_T(\varnothing) = r) P^{(0)}(\#T^{(1)} = n | c_T(\varnothing) = r)$$

$$= \sum_{r \geq 1} P^{(0)}(c_T(\varnothing) = r) P^{(1)}_r(\#F^{(1)} = n),$$

where $P^{(1)}_r$ is the law of a two-type forest with black roots and $r$ tree components. Using again the map $\Gamma$ on all the tree components, we see that this probability is the same as the probability that a monotype GW forest with $r$ tree components has $n$ vertices, where the offspring distribution $\overline{\mu}'$ has generating function $\mathbf{G}_{\mu_1} \circ \mathbf{G}_{\mu_0}$. The Otter–Dwass formula shows that this is equal to $rn^{-1}\widetilde{P}(W'_n = -r)$, where $W'$ under $\widetilde{P}$ is a random walk with step distribution $\overline{\mu}'$ (which has finite variance). Hence,

$$n^{3/2} P^{(0)}(\#T^{(1)} = n) = \sum_{r \geq 1} r\mu_0(r) \, n^{1/2} \widetilde{P}(W'_n = -r).$$

To conclude, notice that $\sum_r r\mu_0(r) = m_0$, and that the local limit theorem of [11], Theorem XV.5.3 shows that $n^{1/2}\widetilde{P}(W'_n = -r)$ converges to a limit $C > 0$ (a multiple of a Gaussian density evaluated at 0), while remaining uniformly bounded as $r$ varies. By the dominated convergence theorem, it results that $n^{3/2} P^{(0)}(\#T^{(1)} = n)$ converges to $Cm_0 = C_{01}$. □

This estimate allows to use a conditioning argument similar to that used in [18], which will be illustrated in the proof of the next lemma. This idea is the following: if $A_n$ is a set of trees such that $P^{(i)}(A_n) \leq \exp(-n^\varepsilon)$ for large $n$, and some $\varepsilon > 0$, then

$$P^{(i)}(A_n | \#T^{(j)} = n) = \frac{P_\infty^{(i)}(F_1 \in A_n, \#F_1^{(j)} = n)}{P_\infty^{(i)}(\#F_1^{(j)} = n)},$$

which by Lemma 14 is less than $\exp(-n^{\varepsilon/2})$ for all $n$ large. Thus, we can obtain a similar exponential control of the event $A_n$ under the law $P^{(i)}(\cdot | \#T^{(j)} = n)$.

4.5. *The convergence of types lemma.* The goal of this section is to give the asymptotic repartition of vertices of either color in large two-type GW forests. This is known as the convergence of types theorem in the literature on multitype GW processes, and we propose a new approach to it.

For $\mathbf{t} \in \mathcal{T}_0 \sqcup \mathcal{T}_1$, let

$$G_\mathbf{t}^{(i)}(n) = \#\{u \in \mathbf{t} : u \prec u^{(i)}(n)\}, \qquad 0 \leq n \leq \#\mathbf{t}^{(i)} - 1.$$

Notice that $u^{(i)}(n)$ is not counted in the set. We also let by convention $G_\mathbf{t}^{(i)}(\#\mathbf{t}^{(i)}) = \#\mathbf{t} - 1$. The definition of $G_\mathbf{f}^{(i)}(k)$ is similar for a forest $\mathbf{f} \in \mathcal{F}_0 \sqcup \mathcal{F}_1$.

LEMMA 15. *Assume that $(\mu_0, \mu_1)$ are nondegenerate critical offspring distribution, that admits some exponential moments. Then, for any $\gamma > 0$ there exists $\varepsilon > 0$ such that for any $i, j \in \{0, 1\}$, for every $n$ large enough,*

$$(39) \qquad P_\infty^{(i)}\left(\sup_{0 \leq k \leq n} |G_F^{(j)}(k) - (1 + m_j)k| > n^{1/2+\gamma}\right) \leq \exp(-n^\varepsilon),$$

*and similarly, if besides $c \in \{0, 1\}$, for $n$ large enough,*

$$(40)\ P^{(i)}\left(\sup_{0 \leq k \leq \#T^{(j)}} |G_T^{(j)}(k) - (1 + m_j)k| > n^{1/2+\gamma} \Big| \#T^{(c)} = n\right) \leq \exp(-n^\varepsilon)$$

*where we take the convention that the conditional probability on the left-hand side is $0$ if $P^{(i)}(\#T^{(c)} = n) = 0$.*

PROOF. Let $\mathbf{f} \in \mathcal{F}_0$. Notice that

$$\begin{aligned}(41)\quad G_\mathbf{f}^{(1)}(m) &= \Upsilon_\mathbf{f}(u^{(1)}(m)) + \sum_{k=0}^{m-1}(1 + c_\mathbf{f}(u^{(1)}(k)))\mathbb{1}_{\{u^{(1)}(k) \not\vdash u^{(1)}(m)\}} \\ &\quad + \sum_{k=0}^{m-1}(1 + c'_\mathbf{f}(u^{(1)}(k)))\mathbb{1}_{\{u^{(1)}(k) \vdash u^{(1)}(m)\}},\end{aligned}$$

where, for $u \vdash u^{(1)}(m)$ in $\mathbf{f}^{(1)}$,

$$c'_{\mathbf{f}}(u) = \#\{v : \neg v = u, v \prec u^{(1)}(m)\}.$$

Indeed, in (41), are counted the number of (type 0) roots of the forest $\mathbf{f}$ before attaining $u^{(1)}(m)$, and the terms $(1 + c_{\mathbf{f}}(u^{(1)}(k)))$ come from counting vertices of $\mathbf{f}$ by groups of parents of type 1, and their children of type 0. One should be careful, however, that if the parent of a group is an ancestor of $u^{(1)}(m)$, then its children that appear after $u^{(1)}(m)$ in depth-first order should not be counted, hence the terms $c'_{\mathbf{f}}$. This shows that

$$\max_{0 \leq k \leq n} \left| G^{(1)}_{\mathbf{f}}(k) - \sum_{l=1}^{k}(1 + c_{\mathbf{f}}(u^{(1)}(l))) \right|$$
$$\leq \Upsilon_{\mathbf{f}}(u^{(1)}(n)) + \max_{0 \leq k \leq n} \sum_{l=0}^{k}(1 + c_{\mathbf{f}}(u^{(1)}(l))) \mathbb{1}_{\{u^{(1)}(l) \vdash u^{(1)}(k)\}}.$$

Moreover, we claim the variables $c_F(u^{(1)}(n)), n \geq 0$ under $P^{(0)}_{\infty}$ are i.i.d. with law $\mu_1$. This is due to the fact that the random variables $X_u, u \in \mathcal{U}$ with $|u|$ odd that are used in the construction of the tree $\xi$ in Section 2.2 are i.i.d. with law $\mu_1$, and that during the exploration of the tree in depth-first order, each vertex is visited only after all his ancestors have been. Thus, on an intuitive level, the depth-first order exploration does not give any information on the number of children of a vertex before it is visited.

Therefore, for every $\eta > 0$,

$$P^{(0)}_{\infty}\left( \max_{0 \leq k \leq n} \left| \sum_{l=1}^{k}(1 + c_F(u^{(1)}(l))) - (1 + m_1)k \right| \geq n^{1/2+\eta} \right) \leq \exp(-n^{\varepsilon}),$$

for some $\varepsilon > 0$ and all $n$ large enough, where we have used a standard moderate deviation inequality for i.i.d. random variables that admit some exponential moments (see [21], Theorem 2.6).

Therefore, by further using Lemmas 11 and 13, if we let

$$A_n = \left\{ \max_{0 \leq k \leq n} |G^{(1)}_F(k) - (1 + m_1)k| \geq n^{1/2+\gamma} \right\}$$

and

$$B_n = \left\{ \max_{0 \leq k \leq n} \left| \sum_{l=0}^{k}(1 + c_F(u^{(1)}(l))) - (1 + m_1)k \right| \leq n^{1/2+\eta} \right\},$$

it holds that

$$P^{(0)}_{\infty}(A_n) = R(n) + P^{(0)}_{\infty}\bigg( A_n, B_n, \max_{u \preceq u^{(1)}(n)} |u| \leq n^{1/2+\eta},$$
$$\Upsilon_F(u^{(1)}(n)) \leq n^{1/2+\eta}, \max_{u \preceq u^{(1)}(n)} c_F(u) < n^{\eta} \bigg),$$

where $R(n) \leq \exp(-n^\varepsilon)$ for some $\varepsilon > 0$ and $n$ large enough. But on the event that $\max_{u \preceq u^{(1)}(n)} |u| \leq n^{1/2+\eta}$ and $\max_{u \preceq u^{(1)}(n)} c_F(u) < n^\eta$ we have for $n$ large

$$\max_{0 \leq k \leq n} \sum_{l=0}^{k} (1 + c_F(u^{(1)}(l))) \mathbb{1}_{\{u^{(1)}(l) \vdash u^{(1)}(k)\}} \leq n^{1/2+\eta}(1 + n^\eta) \leq n^{1/2+3\eta}.$$

If we choose $3\eta < \gamma$, we finally obtain that for $n$ large, $A_n$ is disjoint from the intersection

$$B_n \cap \left\{ \max_{u \preceq u^{(1)}(n)} |u| \leq n^{1/2+\eta} \right\}$$

$$\cap \{\Upsilon_F(u^{(1)}(n)) \leq n^{1/2+\eta}\} \cap \left\{ \max_{u \preceq u^{(1)}(n)} c_F(u) < n^\eta \right\},$$

so that for $n$ large, $P_\infty^{(0)}(A_n) = R(n) \leq \exp(-n^\varepsilon)$.

The case $i = j = 0$ is similar but easier, as the term $\Upsilon_\mathbf{f}(u^{(1)}(n))$ of (41) does not appear anymore. Details are left to the reader.

We now pass to the conditioned statements. We apply the conditioning argument mentioned in Section 4.4. We first treat the case $c = j = 0$. Using Lemma 14, for some constant $C > 0$

$$P_\infty^{(i)}\left( \sup_{0 \leq k \leq \#F^{(0)}} |G_F^{(0)}(k) - (1+m_0)k| > n^{1/2+\gamma} \Big| \#F_1^{(0)} = n \right)$$

$$= \frac{P_\infty^{(i)}(\sup_{0 \leq k \leq n} |G_F^{(0)}(k) - (1+m_0)k| > n^{1/2+\gamma}, \#F_1^{(0)} = n)}{P^{(i)}(\#T^{(0)} = n)}$$

$$\leq Cn^{3/2} P_\infty^{(i)}\left( \sup_{0 \leq k \leq n} |G_F^{(0)}(k) - (1+m_0)k| > n^{1/2+\gamma} \right)$$

$$\leq \exp(-n^\varepsilon),$$

for some $\varepsilon > 0$ and all $n$ large. Notice the little artifact here: rather than considering a single tree, we have considered a forest whose first component is conditioned. This yields the wanted result (40) for $j = c = 0$, by restricting the sup to $0 \leq k \leq \#F^{(0)} - 1$, but it also gives us a little more: namely that $P^{(i)}(|G_F^{(0)}(\#F_1^{(0)}) - (1+m_0)\#F_1^{(0)}| \geq n^{1/2+\gamma} | \#F_1^{(0)} = n) \leq e^{-n^\varepsilon}$ for large $n$. Since $G_F^{(0)}(\#F_1^{(0)} - 1) \leq \#F_1 \leq G_F^{(0)}(\#F_1^{(0)}) + 1$, this shows that $P^{(i)}(|\#T - (1+m_0)n| \geq n^{1/2+\gamma} | \#T^{(0)} = n) \leq e^{-n^\varepsilon}$ for large $n$. Since $\#T^{(0)} + \#T^{(1)} = \#T$,

$$(42) \qquad P^{(i)}(|\#T^{(1)} - m_0 n| \geq n^{1/2+\gamma} | \#T^{(0)} = n) \leq e^{-n^\varepsilon}$$

for large $n$. Thanks to this control on the number of vertices of type 1, we obtain for large $n$,

$$P^{(i)}\left(\max_{0\leq k\leq \#T^{(1)}}|G_T^{(1)}(k)-(1+m_1)k|>n^{1/2+\gamma}\Big|\#T^{(0)}=n\right)$$

$$\leq \frac{P_\infty^{(i)}(\max_{0\leq k\leq (m_0+\varepsilon')n}|G_F^{(1)}(k)-(1+m_1)k|>n^{1/2+\gamma})+\exp(-n^\varepsilon)}{P^{(i)}(\#T^{(0)}=n)},$$

where $0<\varepsilon'$, and the $\exp(-n^\varepsilon)$ term bounds the probability that $\#T^{(1)}$ is larger than $n(m_0+\varepsilon')$. This expression is less than $\exp(-n^{\varepsilon''})$ for some $\varepsilon''$ and large $n$ because of the unconditioned control (39) on $G_F^{(1)}$ (notice that the maximum is taken over $1\leq k\leq Dn$ for some constant $D>0$ rather than 1, but this does not matter up to a change in the constant $\gamma$). The remaining cases for $j,c$ are symmetric. $\square$

PROOF OF LEMMA 9. From Lemma 15, we obtain that $(G_T^{(j)}([\#T^{(j)}t])/\#T^{(j)}, 0<t\leq 1)$ converges in probability to the function $((1+m_j)t, 0\leq t\leq 1)$ under $P^{(i)}(\cdot|\#T^{(c)}=n)$, for the uniform norm. In fact, it holds that $(G_T^{(j)}([\#T^{(j)}t])/\#T, 0<t\leq 1)$ converges in probability to the identity function, because $\#T^{(j)}/\#T$ converges to $(1+m_j)^{-1}$ under $P^{(i)}(\cdot|\#T^{(c)}=n)$ as was shown in the proof of the previous lemma. On the other hand, the process $\overline{J}_T^{(c)}$ of Lemma 9 is the right-continuous inverse function of $(G_T^{(c)}([\#T^{(c)}t])/\#T, 0<t\leq 1)$ [this motivates our convention $G_T^{(c)}(\#T^{(c)})=\#T-1$], so it also converges in probability to the identity function for the uniform norm, as claimed. $\square$

4.6. *Convergence of the height process.* The last ingredient that we need is the fact that the height process of a monotype GW forest with $r$ components, conditioned by the number of its vertices, converges to a scaled Brownian excursion. The case $r=1$ is known (see [2, 9, 18]). The result for $r\geq 1$ is suggested in [22], Chapter 5. Recall that when $\mu_0=\mu_1=\mu$, the index $i$ in the probability $P^{(i)}$ becomes irrelevant, so we let $P_r$ be the law of a (monotype) GW forest with offspring distribution $\mu$ and $r$ tree components.

THEOREM 16. *Let $\mu$ be a critical nondegenerate offspring distribution, admitting small exponential moments, and let $\sigma_\mu$ be its variance. Then for any $r>0$, the process $(n^{-1/2}H_{(n-1)t}^F, 0\leq t\leq 1)$ under $P_r(\cdot|\#F=n)$ converges in distribution for the uniform topology to $2\sigma_\mu^{-1}\mathrm{e}$ under $\mathbb{N}^{(1)}$.*

PROOF. We first note that under $P(\cdot|\#T=n+1, c_T(\varnothing)=r)$, the forest $1T_1\cup\cdots\cup rT_r$ has same distribution as $F$ under $P_r(\cdot|\#F=n)$. Thus, under $P_r(\cdot|\#F = n)$, the height process of $F$ has same law as the process obtained by concatenation of $H^{T_1}, \ldots, H^{T_r}$ under $P(\cdot|\#T = n+1, c_T(\varnothing) = r)$, so the rescaled process $n^{-1/2}H^F_{(n-1)}.$ under $P_r(\cdot|\#F = n)$ has same law as $(n^{-1/2}(H^T_{((n-1)t+1)/n} - 1), 0 \leq t \leq 1)$ under $P(\cdot|\#T = n+1, c_T(\varnothing) = r)$.

Let $n_1, \ldots, n_r$ be positive integers with sum $n$, and let $\mathbf{t}_1, \ldots, \mathbf{t}_r \in \mathcal{T}$ be such that $\#\mathbf{t}_j = n_j$. Using the branching properties of GW trees, we have

$$
(43) \qquad P((T_1, \ldots, T_r) = (\mathbf{t}_1, \ldots, \mathbf{t}_r)|c_T(\varnothing) = r, \#T_j = n_j, 1 \leq j \leq r)
$$
$$
= \prod_{i=1}^r P(T = \mathbf{t}_i|\#T = n_i),
$$

so that given $\#T_i = n_i, 1 \leq i \leq r$, under $P(\cdot|\#T = n+1, c_T(\varnothing) = r)$ the $T_i$'s are independent GW trees, respectively conditioned to have size $n_i$. We next claim that for every $\varepsilon > 0$, and every $r$ such that $\mu(r) > 0$,

$$
(44) \qquad P\left(\max_{i \leq c_T(\varnothing)} \#T_i/n \leq 1 - \varepsilon \Big| \#T = n+1, c_T(\varnothing) = r\right) \underset{n \to \infty}{\to} 0,
$$

which will be proved later on.

Let $T_\star$ be the largest tree among the $T_i$'s under $P(\cdot|\#T = n+1, c_T(\varnothing) = r)$ (or the first largest tree, if several trees have the maximal size). According to (43), (44), and by the known $r = 1$ case of the theorem, $(n^{-1/2}H^{T_\star}_{(\#T_\star - 1)t}, 0 \leq t \leq 1)$ converges in distribution to $2\sigma_\mu^{-1}\mathbf{e}$ under $\mathbb{N}^{(1)}$. Since the number of individuals of the $r - 1$ other subtrees is $o(n)$ in probability, the maximal height of a vertex of these trees is $o(n^{1/2})$ in probability, as can be checked from the case $r = 1$ of the statement, and using (43). Thus, it easily follows that $n^{-1/2}\|H^{T_\star}_{(\#T_\star - 1)\cdot} - H^T_{n\cdot}\|_\infty$ goes to 0 in probability under $P(\cdot|\#T = n+1, c_T(\varnothing) = r)$, which yields the wanted result.

To argue (44), we first observe that the same statement holds without conditioning on $c_T(\varnothing)$. Indeed, the known $r = 1$ case of the theorem shows that under $P(\cdot|\#T = n+1)$, $(n^{-1/2}H^T_{nt}, 0 \leq t \leq 1)$ converges in distribution for the uniform topology to $2\sigma_\mu^{-1}\mathbf{e}$. Recall that the Brownian excursion is a.s. strictly positive on $(0, 1)$, and let $(t_i, 1 \leq i \leq c_T(\varnothing))$ be the ordered list of integers such that $H^T_{t_i} = 1$ [and $t_{c_T(\varnothing)+1} = \#T - 1$]. Then the lengths of the intervals $([t_i, t_{i+1}], 1 \leq i \leq c_T(\varnothing))$, is exactly $(\#T_i, 1 \leq i \leq c_T(\varnothing))$. Now, if $f_n$ is a sequence of continuous functions converging for the uniform topology to $f$ which is positive on $(0, 1)$, then for any $\varepsilon > 0$, for $n$ large enough, $f_n$ is positive on $(\varepsilon, 1 - \varepsilon)$, and it follows that $P(\exists i : t_i/n \in (\varepsilon, 1 - \varepsilon)|\#T = n+1) \to 0$, which is equivalent to the wanted property [observe $t_1/n \to 0$ while $t_{c_T(\varnothing)+1}/n = 1$].

On the other hand, for fixed $r \in \mathbb{N}$, $P(c_T(\varnothing) = r|\#T = n+1) = \mu(r) \times P_r(\#F = n)/P(\#T = n+1) \to r\mu(r)$ as $n \to \infty$, by the Otter–Dwass formula and the local limit theorem, as in Lemma 14. Since $\sum_r r\mu(r) = 1$, the law

of $c_T(\varnothing)$ under $P(\cdot|\#T = n+1)$ converges weakly. Equation (44) is now an elementary consequence of this and the previous paragraph. □

We are now ready to prove the first half of Theorem 8, which we state as:

PROPOSITION 17. *Let $(\mu_0, \mu_1)$ be nondegenerate, critical and admit some exponential moments. Then the process $(n^{-1/2} H^T_{(\#T-1)t}, 0 \leq t \leq 1)$ under $P^{(i)}(\cdot|\#T^{(j)} = n)$ converges in distribution to the process $2\sigma^{-1}\sqrt{1+m_j}\mathrm{e}$ under $\mathbb{N}^{(1)}$.*

PROOF. Suppose first that $i = j = 0$ and recall the definition of the mapping $\Gamma_\mathbf{t}$, for $\mathbf{t} \in \mathcal{T}_0$. Under $P^{(0)}$, we know that $\Gamma(T)$ is a monotype GW tree with offspring distribution $\overline{\mu}$ and $\mathbf{G}_{\overline{\mu}} = \mathbf{G}_{\mu_0} \circ \mathbf{G}_{\mu_1}$. Moreover, under $P^{(0)}(\cdot|\#T^{(0)} = n)$, $\Gamma(T)$ has the law $P_{\overline{\mu}}(\cdot|\#T = n)$ of a conditioned monotype GW tree, because $\Gamma_T$ maps each vertex of $T^{(0)}$ to a vertex of $\Gamma(T)$ in a one-to-one way.

On the other hand, the formula $\mathbf{G}_{\overline{\mu}} = \mathbf{G}_{\mu_0} \circ \mathbf{G}_{\mu_1}$ and the fact that $\mu_0$ and $\mu_1$ admit some exponential moments entail that $\overline{\mu}$ itself admits exponential moments. Therefore Theorem 16 applies and it holds that under $P^{(0)}(\cdot|\#T^{(0)} = n)$, $(n^{-1/2} H^{\Gamma(T)}_{(n-1)t}, 0 \leq t \leq 1)$ converges in distribution to $2\sigma_{\overline{\mu}}^{-1}\mathrm{e}$ under $\mathbb{N}^{(1)}$. We can compute $\sigma_{\overline{\mu}}$ by differentiating twice $\mathbf{G}_{\overline{\mu}}$, and we find

$$\sigma_{\overline{\mu}}^2 = m_0 \sigma_1^2 + m_1^2 \sigma_0^2. \tag{45}$$

Next, for every $\mathbf{t} \in \mathcal{T}_0$, it is an elementary exercise to check that

$$|H^{\mathbf{t}}_k - 2H^{\Gamma(\mathbf{t})}_{J^{(0)}_\mathbf{t}(k)-1}| \leq 2|H^{\Gamma(\mathbf{t})}_{J^{(0)}_\mathbf{t}(k)-1} - H^{\Gamma(\mathbf{t})}_{J^{(0)}_\mathbf{t}(k)}| + 1 \tag{46}$$

for every $0 \leq k \leq \#\mathbf{t} - 1$, with the convention that $H^{\Gamma(\mathbf{t})}_{\#\mathbf{t}^{(0)}} = 0$.

Lemma 9 states that the process $\overline{J}^{(0)}_T = (J^{(0)}_T([(\#T-1)t])/\#T^{(0)}, 0 \leq t \leq 1)$ under $P^{(0)}(\cdot|\#T^{(0)} = n)$ converges in probability to the identity $(t, 0 \leq t \leq 1)$. This convergence holds jointly with that of $(n^{-1/2} H^{\Gamma(T)}_{(n-1)t}, 0 \leq t \leq 1)$ to a rescaled Brownian excursion. Skorokhod's representation theorem ensures that there exists a probability space on which random processes $(H^n, J^n)$ converge a.s. to $(B, (t, 0 \leq t \leq 1))$ for the uniform norm, where $H^n, J^n$ have same law as

$$(n^{-1/2} H^{\Gamma(T)}_{(\#T^{(0)}-1)t}, 0 \leq t \leq 1), \qquad \overline{J}^{(0)}_T$$

under $P^{(0)}(\cdot|\#T^{(0)} = n)$, and $B$ has same law as $2\sigma_{\overline{\mu}}^{-1}\mathrm{e}$ under $\mathbb{N}^{(1)}$. Then, the composed function $(H^n_{(nJ^n(t)-1)/(n-1)}, 0 \leq t \leq 1)$ converges a.s. uniformly to $B$, which says that $n^{-1/2} H^{\Gamma(T)}_{(n\overline{J}^{(0)}_T-1)}$ under $P^{(0)}(\cdot|\#T^{(0)} = n)$ converges to $2\mathrm{e}/\sigma_{\overline{\mu}}$ under $\mathbb{N}^{(1)}$. Also, it holds that $\sup_{0\le t\le 1}|H^n_{(nJ^n(t)-1)/(n-1)} - H^n_{nJ^n(t)/(n-1)}|$ converges to 0 a.s. This, paired with equation (46), implies that for every $\varepsilon > 0$,

$$P^{(0)}\left(\sup_{0\le t\le 1}\frac{1}{\sqrt{n}}|H^T_{(\#T-1)t} - 2H^{\Gamma(T)}_{n\overline{J}^{(0)}_T(t)-1}| > \varepsilon \Big| \#T^{(0)} = n\right) \underset{n\to\infty}{\longrightarrow} 0.$$

Finally, this entails that $(n^{-1/2}H^T_{(\#T-1)t}, 0 \le t \le 1)$ converges to $4\sigma_{\overline{\mu}}^{-1}\mathrm{e}$ under $\mathbb{N}^{(1)}$. The result follows from the fact that $4\sigma_{\overline{\mu}}^{-1} = 2\sigma^{-1}\sqrt{1+m_0}$, as is easily checked from (17) and (45).

We next treat the case where $i=0, j=1$. We apply the transformations $\pi, \Gamma'$ of Lemma 11, that skips the first generation, and then squeezes odd generations. Recall that $\pi, \Gamma'$ were defined on forests, and that they take values in the set of forests (even if the initial forest has only one component). Notice that if $\mathbf{t} \in \mathcal{T}_0$, then $1\mathbf{t} = \{1u : u \in \mathbf{t}\} \in \mathcal{F}_0$. Under $P^{(0)}$, the forest $\pi(1T)$ is a GW forest with a random number of components, which is given by $c_T(\varnothing)$ and is independent of the components of $\pi(1T)$, and it holds that the law of $\pi(1T)$ under $P^{(0)}(\cdot|\#T^{(1)} = n)$ is the probability measure $d\mu_0(r)P_r^{(1)}(d\mathbf{f})$ on $\mathcal{F}_0$, conditioned on the event $\#F^{(1)} = n$. Then, $\Gamma'(1T)$ under this law is a monotype GW forest with law $d\mu_0(r)P_r(d\mathbf{f})$ given $\#F = n$, where $P_r$ is the law of a GW forest with $r$ trees and offspring distribution $\overline{\mu}'$ whose generating function is $\mathbf{G}_{\mu_1} \circ \mathbf{G}_{\mu_0}$.

Notice that under $P^{(0)}(\cdot|\#T^{(1)} = n, c_T(\varnothing) = r)$, the forest $\pi(1T)$ has law $P_r^{(0)}(\cdot|\#F^{(1)} = n)$. Therefore, under $P^{(0)}(\cdot|c_T(\varnothing) = r, \#T^{(1)} = n)$, we obtain that $\Gamma'(1T)$ has law $P_r(\cdot|\#F = n)$. By Theorem 16, under this law, the process $n^{-1/2}H^{\Gamma'(1T)}_{(n-1)\cdot}$ converges to a Brownian excursion scaled by $2/\sigma_{\overline{\mu}'} = \sqrt{1+m_1}/\sigma$. One checks that a companion formula to (46) holds, namely that

$$(47) \qquad |H^{\mathbf{t}}_k - 2H^{\Gamma'(\mathbf{t})}_{J^{(1)}_{\mathbf{t}}(k)-1}| \le 2|H^{\Gamma'(1\mathbf{t})}_{J^{(1)}_{\mathbf{t}}(k)-1} - H^{\Gamma'(1\mathbf{t})}_{J^{(1)}_{\mathbf{t}}(k)}| + 2,$$

for $0 \le k \le \#\mathbf{t} - 1$, with the convention that $H^{\Gamma'(1\mathbf{t})}_{-1} = H^{\Gamma'(1\mathbf{t})}_{\#\mathbf{t}^{(1)}} = 0$. The same arguments as in the case $i=j=0$ entail that $(n^{-1/2}H^F_{(\#F-1)t}, 0 \le t \le 1)$ under $P_r^{(0)}(\cdot|\#F^{(1)} = n)$ converges to a Brownian excursion scaled by $2\sqrt{1+m_1}/\sigma$.

To complete the proof, since we are dealing with a *random* number of components, whose law depends on $n$, that is, a mixture in $r$ of the laws $P^{(0)}(\cdot|c_T(\varnothing) = r, \#T^{(1)} = n)$, it suffices to use that $c_T(\varnothing)$ converges in distribution under $P^{(0)}(\cdot|\#T^{(1)} = n)$, as seen in the proof of Theorem 16. $\square$

A remarkable difference between Proposition 17 and Theorem 16 is that in Proposition 17, we condition only on the total number of a portion of the vertices of $T$, instead of the total size of the tree. We have seen that this is the right conditioning to do when considering random maps conditioned on the number of faces, or vertices, but one may wonder whether the result is still true under the probability laws $P^{(i)}(\cdot|\#T=n)$, with different scaling constants (this would allow to consider random maps conditioned by the number of edges). We expect this to be true, but the methods that are used in the present work are powerless to address this issue.

Let us end this section with a result on the increments of the height processes.

LEMMA 18. *Let $\mu$ be a (monotype) nondegenerate critical offspring distribution. Then for every fixed $r \geq 1$, for every $\gamma > 0$, under $P_r(\cdot|\#F=n)$, the quantity $n^{-\gamma} \sup_{0 \leq k \leq \#F-1} |H_k^F - H_{k+1}^F|$ converges to 0 in probability, with the convention that $H_{\#F}^F = 0$.*

PROOF. Using the same argument as in the proof of Theorem 16, we may prove the statement under $P(\cdot|\#T=n)$, that is for a single GW tree conditioned to have $n$ vertices. The positive jumps of the height process of any tree are $+1$, and then, only the negative jumps have to be controlled. For elementary symmetry reasons, the largest negative jump in $H^T$, plus 1, has the same law as the largest number of consecutive steps $+1$ in $H^T$. But, in a nonconditioned GW tree, a run of steps $+1$ in the height process has a geometric distribution: the probability that a run has length $k$ is $\mu(0)(1-\mu(0))^{k-1}$, and the different runs are independent.

Denote by $G_1, G_2, \ldots, G_K$ the sizes of successive runs of $+1$ in the height process, where $K$ is random and is bounded above by $n$ under $P(\cdot|\#T=n)$. Let $\gamma > 0$ be fixed. Thanks to the conditioning argument (the Otter–Dwass formula and the local limit theorem), since the function max is nondecreasing, we have

$$P\left(\sup_{0 \leq k \leq \#T-1} |H_k^T - H_{k+1}^T| \geq n^{\gamma/2} \Big| \#T=n\right)$$
$$= O\left(n^{3/2} P_\infty\left(\sup_{1 \leq k \leq n} G_k \geq n^{\gamma/2}+1\right)\right)$$
$$= O(n^{5/2} P_\infty(G_1 \geq n^{\gamma/2}+1))$$

and since $\mu(0) \in (0,1)$, this is bounded by $\exp(-n^\varepsilon)$ for some $\varepsilon > 0$, for $n$ large enough. □

**5. Convergence of the label process.** The proof of the convergence of the second components in Theorem 8 will be done by showing that their finite marginal distributions converge, combined by a tightness argument.

5.1. *Controlling the branching in conditioned trees.* The convergence of finite-marginal distributions first needs some improvements and variations around Lemma 13. A conditioned version of Lemma 13 holds:

LEMMA 19. *Assume that the pair $(\mu_0, \mu_1)$ is nondegenerate, critical and has some exponential moments. Then for every $\eta > 0$ there exists $\varepsilon > 0$ such that for $n$ large enough, and $i, j \in \{0, 1\}$,*

$$(48) \qquad P^{(i)}\left(\max_{u \in T} c_T(u) \geq n^\eta \Big| \#T^{(j)} = n\right) \leq \exp(-n^\varepsilon).$$

PROOF. We have

$$P_\infty^{(i)}\left(\max_{u \preceq u^{(j)}(n)} c_F(u) \geq n^\eta \Big| \#F_1^{(j)} = n\right) \leq \frac{P_\infty^{(i)}(\max_{u \preceq u^{(j)}(n)} c_F(u) \geq n^\eta)}{P^{(i)}(\#T^{(j)} = n)}$$

which is smaller than $\exp(-n^\varepsilon)$ for all large $n$, by Lemmas 13 and 14. Since $u \preceq u^{(j)}(n)$ for every $u \in F_1$ given $\{\#F_1^{(j)} = n\}$, this entails (48). $\square$

Let $\mathbf{t} \in \mathcal{T}_i$, $u \in \mathbf{t}$, $0 \leq h \leq |u|$, $k \geq 1$, $1 \leq l \leq k$. We define

$$A_\mathbf{t}^{(j)}(u, k, l, h) = \#\{v \vdash u : c_\mathbf{t}(v) = k, u \in vl\mathbf{t}_{vl}, v \in \mathbf{t}^{(j)}, |v| > |u| - h\}$$

the number of ancestors of $u$ which are at distance at most $h$ from $u$, with type $j$, $k$ children, and such that $u$ is a descendant of the $l$th of these children. We let $A_\mathbf{f}^{(j)}(u, k, l, h)$ be the similar quantity for a forest $\mathbf{f} \in \mathcal{F}_0 \sqcup \mathcal{F}_1$ and $u \in \mathbf{f}$. Note that if $\max_u c_\mathbf{f}(u) \leq K$ then $A_\mathbf{f}(u, k, l, h) = 0$ for any $k > K$, any $l, h$ and any $u \in \mathbf{f}$.

LEMMA 20. *Let $(\mu_0, \mu_1)$ be a nondegenerate critical offspring distribution admitting some exponential moments. For every $\gamma > 0$, $M > 0$, and $i, j, c \in \{0, 1\}$, there exists $\varepsilon > 0$ such that, for $n$ large enough*

$$P^{(i)}\left(\sup_{k \geq 1, 1 \leq l \leq k} \sup_{u \in T, n^\gamma \leq h \leq |u|} \frac{|A_T^{(j)}(u, k, l, h) - \mu_j(k)h/(2m_j)|}{h^{1/2+\gamma} k^{-M}} \geq 1 \Big| \#T^{(c)} = n\right)$$
$$\leq \exp(-n^\varepsilon).$$

PROOF. Let $\gamma > 0$, $M > 0$ be fixed, and choose $\eta < \gamma^2/M$. By Lemmas 19, 15 and the conditioning argument, we know that $P^{(i)}(\max_{u \in T} c_T(u) \leq n^\eta, \#T \leq Cn | \#T^{(c)} = n) \geq 1 - \exp(-n^\varepsilon)$ for some constants $C, \varepsilon$, and $n$ large enough. Since $A_\mathbf{t}^{(j)}(u, k+1, l, h) = 0$ whenever $\max_{u \in \mathbf{t}} c_\mathbf{t}(u) \leq k$, on the event $\{\max_{u \in T} c_T(u) \leq n^\eta, \#T \leq Cn\}$ we have

$$\sup_{k > n^\eta, 1 \leq j \leq k} \sup_{u \in T, n^\gamma \leq h \leq |u|} k^M h^{-1/2-\gamma} \left| A_T^{(j)}(u, k, l, h) - \frac{\mu_j(k)}{2m_j} h \right|$$

$$\leq \sup_{k\geq n^\eta, u\in T} |u|^{1/2-\gamma} \frac{\mu_j(k)}{2m_j} k^M.$$

Since $|u| \leq \#T \leq Cn \leq Ck^{1/\eta}$ for $k \geq n^\eta$ and since $\mu_j$ has small exponential moments, this is smaller than 1 when $n$ is large enough.

We now estimate $A_T^{(j)}(u, k, l, h)$ for $k \leq n^\eta$. We start with considering forests. Let $k \leq n^\eta$, and $l \in \{1, \ldots, k\}$ be fixed, and $C$ be the same constant as above. By using Lemma 11,

$$P_\infty^{(i)}\left(\sup_{\substack{u \preceq u^{(j)}(Cn), \\ |u|\geq h\geq n^\gamma}} \frac{|A_T^{(j)}(u,k,l,h) - \mu_j(k)h/(2m_j)|}{h^{1/2+\gamma} k^{-M}} \geq 1\right)$$

$$(49) \qquad \leq P_{[n^{1/2+\eta}]}^{(i)}\left(\sup_{\substack{u \preceq u^{(j)}(Cn), \\ n^\gamma \leq h \leq |u|}} \frac{|A_T^{(j)}(u,k,l,h) - \mu_j(k)h/(2m_j)|}{h^{1/2+\gamma} k^{-M}} \geq 1, \mathcal{B}_n\right)$$

$$+ \exp(-n^\varepsilon)$$

where $\mathcal{B}_n = \{\sup_{u \prec u^{(j)}(Cn)} |u| \leq n^{1/2+\eta}, \Upsilon_F(u^{(j)}(Cn)) \leq n^{1/2+\eta}\}$. The probability on the right-hand side can be bounded as follows, using Lemma 12:

$$P_{[n^{1/2+\eta}]}^{(i)}\left(\sup_{\substack{u \preceq u^{(j)}(Cn), \\ n^\gamma \leq h \leq |u| \leq n^{1/2+\eta}}} \frac{|A_T^{(j)}(u,k,l,h) - \mu_j(k)h/(2m_j)|}{h^{1/2+\gamma} k^{-M}} \geq 1\right)$$

$$\leq E_{[n^{1/2+\eta}]}^{(i)}\left[\sum_{u\in F} \mathbb{1}\left\{\sup_{n^\gamma \leq h \leq |u|} \frac{|A_T^{(j)}(u,k,l,h) - \mu_j(k)h/(2m_j)|}{h^{1/2+\gamma} k^{-M}} \geq 1\right\}\right.$$

$$(50) \qquad\qquad\qquad\qquad\qquad\qquad\qquad\qquad\qquad \left.\times \mathbb{1}_{\{|u|\leq n^{1/2+\eta}\}}\right]$$

$$\leq (1 \vee m_i) n^{1/2+\eta}$$

$$\times \sum_{h'=n^\gamma}^{[n^{1/2+\eta}]} \widehat{P}^{(i),h'}\left(\sup_{n^\gamma \leq h \leq h'} \frac{|A_T^{(j)}(V,k,l,h) - \mu_j(k)h/(2m_j)|}{h^{1/2+\gamma} k^{-M}} \geq 1\right).$$

Then, we argue that $A_T^{(j)}(V, k, l, h)$ under $\widehat{P}^{(i),h'}$ is a Binomial random variable $B(m, p)$ with parameters $p = \mu_j(k)/m_j$ and either $m = [h/2 + 1]$ or $m = [(h+1)/2]$ depending on the parity of $i, j, h, h'$. Hoeffding's inequality (33) entails that

$$\widehat{P}^{(i),h'}\left(\sup_{n^\gamma \leq h \leq h'} \frac{|A_T^{(j)}(V,k,l,h) - \mu_j(k)h/(2m_j)|}{h^{1/2+\gamma} k^{-M}} \geq 1\right)$$

$$\leq \sum_{n^\gamma \leq h \leq h'} \widehat{P}^{(i),h'}\left(\left|A_T^{(j)}(V,k,l,h) - \frac{\mu_j(k)h}{2m_j}\right| \geq k^{-M}h^{1/2+\gamma}\right)$$

$$\leq 2h' \max_{n^\gamma \leq h \leq h'} \exp(-k^{-2M}h^{1+2\gamma}/(2m))$$

$$\leq 2h' \exp(-n^{-2M\eta+2\gamma^2}/2).$$

Finally, the expression (50) is bounded by

$$Kn^{1/2+\eta} \sum_{h'=n^\gamma}^{[n^{1/2+\eta}]} h' \exp(-n^{-2M\eta+2\gamma^2}/2) + \exp(-n^\varepsilon),$$

for some $K > 0$ and large $n$, and this is $\leq \exp(-n^{\varepsilon'})$, for large $n$ and some $\varepsilon' > 0$.

This entails that

$$P_\infty^{(i)}\left(\sup_{\substack{k \leq n^\gamma,\, u \prec u^{(j)}([Cn]), \\ 1 \leq l \leq k \quad |u| \geq h \geq n^\gamma}} \frac{|A_T^{(j)}(u,k,l,h) - \mu_j(k)h/(2m_j)|}{h^{1/2+\gamma}k^{-M}} \geq 1\right)$$

$$\leq \exp(-n^{\varepsilon''}),$$

for some $\varepsilon'' > 0$ and $n$ large. To obtain the conditioned statement, we apply our conditioning argument once again. By definition of $C$, $P^{(i)}(\#T > Cn| \#T^{(c)} = n) \leq \exp(-n^\varepsilon)$ for large $n$. So

$$P^{(i)}\left(\sup_{\substack{k \leq n^\eta, \quad u \in T, \\ 1 \leq l \leq k\, n^\gamma \leq h \leq |u|}} \frac{|A_T^{(j)}(u,k,l,h) - \mu_j(k)h/(2m_j)|}{k^{-M}h^{1/2+\gamma}} \geq 1\Big|\#T^{(c)} = n\right)$$

$$\leq P^{(i)}(T^{(c)} = n)^{-1}$$

$$\times \left(P_\infty^{(i)}\left(\sup_{\substack{k \leq n^\eta,\, u \prec u^{(j)}([Cn]), \\ 1 \leq l \leq k \quad n^\gamma \leq h \leq |u|}} \frac{|A_F^{(j)}(u,k,l,h) - \mu_j(k)h/(2m_j)|}{k^{-M}h^{1/2+\gamma}} \geq 1\right)\right.$$

$$\left. + \exp(-n^\varepsilon)\right)$$

$$\leq C'n^{3/2}(\exp(-n^{\varepsilon''}) + \exp(-n^\varepsilon)),$$

for some constant $C' > 0$, which yields the result. □

5.2. *A bound on the Hölder norm of the height process.* The second ingredient which is required to prove Theorem 8 is the following result, showing that the $\alpha$-Hölder norm of the height process under $P^{(i)}(\cdot|\#T^{(j)} = n)$ is tight for any $\alpha < 1/2$.

We start by stating a monotype version of the result we need:

PROPOSITION 21. *Let $\mu$ be a nondegenerate critical offspring distribution which admits some exponential moments. Let $P_r$ be the law of a (monotype) GW forest with offspring distribution $\mu$ and $r$ components. Then for every $r \geq 1$, $\varepsilon > 0$ and $\alpha \in (0, 1/2)$, there exists $C > 0$ such that*

$$(51) \qquad \sup_{n \in \mathbb{N}} P_r \left( \sup_{0 \leq s \neq t \leq 1} \frac{|H^F_{(n-1)s} - H^F_{(n-1)t}|}{\sqrt{n}|s-t|^\alpha} > C \,\Big|\, \#F = n \right) \leq \varepsilon.$$

PROOF. We claim that it is sufficient to prove the statement for $r = 1$, by using an argument similar to that of the proof of Theorem 16 for the more general $r \geq 1$ case. Indeed, under $P_r(\cdot|\#F = n)$, recall that $n^{-1/2} H^F_{(n-1)\cdot}$ has same law as the concatenation of the paths $(n^{-1/2} H^{T_i}_{(n-1)t}, 0 \leq t \leq (\#T_i - 1)/n - 1), 1 \leq i \leq r$ under $P(\cdot|\#T = n + 1, c_T(\varnothing) = r)$, under which the trees $T_i, 1 \leq i \leq r$, are mixtures of independent GW trees conditioned by their sizes, and hence have $\alpha$-Hölder norm $\leq C$ with high probability by the $r = 1$ case. Now, it is elementary to check that if $(\#T_i)^{-1/2} H^{T_i}_{(\#T_i - 1)\cdot}$ has $\alpha$-Hölder norm $\leq C$, then the same is true of $(n^{-1/2} H^{T_i}_{(n-1)t}, 0 \leq t \leq (\#T_i - 1)/(n-1))$, because $n \geq \#T_i$, and the concatenation of the paths still has $\alpha$-Hölder norm $\leq 2^{1-\alpha} C$, hence giving the result for $r \geq 1$.

Hence, we are down to show that for every $\varepsilon > 0$ and $\alpha \in (0, 1/2)$, there exists $C > 0$ such that

$$(52) \qquad \sup_{n \in \mathbb{N}} P \left( \sup_{0 \leq s \neq t \leq 1} \frac{|H^T_{(n-1)s} - H^T_{(n-1)t}|}{\sqrt{n}|s-t|^\alpha} > C \,\Big|\, \#T = n \right) \leq \varepsilon.$$

We define the depth-first traversal, or contour order of a tree $\mathbf{t}$ as a function:

$$F_{\mathbf{t}} : \{0, \ldots, 2\#\mathbf{t} - 2\} \to \{\text{ vertices of } \mathbf{t}\},$$

which we regard as a walk around $\mathbf{t}$, as follows: $F_{\mathbf{t}}(0) = \varnothing$, and given $F_{\mathbf{t}}(i) = z$, choose, if possible, and according to the depth-first order the smallest child $w$ of $z$ which has not already been visited, and set $F_{\mathbf{t}}(i + 1) = w$. If not possible, let $F_{\mathbf{t}}(i + 1)$ be the father of $z$.

For any $0 \leq k \leq 2\#\mathbf{t} - 2$, set $\widehat{H}^{\mathbf{t}}(k) = |F_{\mathbf{t}}(k)|$. The *contour process* $(\widehat{H}^{\mathbf{t}}_s, 0 \leq s \leq 2\#\mathbf{t} - 2)$ is then obtained by interpolating linearly the sequence $(\widehat{H}^{\mathbf{t}}(k))$ between integer abscissa. For any tree $\mathbf{t} \in \mathcal{T}$, $H^{\mathbf{t}}$ is a simple function of $\widehat{H}^{\mathbf{t}}$: let $m_{\mathbf{t}}(0) = 0$, and for any $i \geq 1$, $m_{\mathbf{t}}(i) = \min\{j, j > m_{\mathbf{t}}(i - 1), \widehat{H}^{\mathbf{t}}(j) > \widehat{H}^{\mathbf{t}}(j - 1)\}$, then $H^{\mathbf{t}}(k) = \widehat{H}^{\mathbf{t}}(m_{\mathbf{t}}(k))$. In fact, $m_{\mathbf{t}}(k) = \inf\{j, F_{\mathbf{t}}(j) = u(k)\}$. One may check inductively on $k$ that,

$$(53) \qquad m_{\mathbf{t}}(k) + H^{\mathbf{t}}(k) = 2k \qquad \text{for any } k \geq 0.$$

We will prove (52), using a similar property for the contour process. The first arguments can be found in [13], Lemma 1. Gittenberger [12] proved

(in a stronger form) that for all $s$, $t$, $a > 0$

$$P\left(\left|\frac{\widehat{H}^T_{(2n-2)s} - \widehat{H}^T_{(2n-2)t}}{\sqrt{n}}\right| \geq a \Big| \#T = n\right) \leq C_1|s-t|^{-1}\exp(-C_2 a|s-t|^{-1/2}),$$

which gives, for any $p > 0$, $E(|n^{-1/2}(\widehat{H}^T_{(2n-2)s} - \widehat{H}^T_{(2n-2)t})|^p|\#T = n) \leq C(p)|s-t|^{p/2-1}$. By applying the uniform version of Kolmogorov's criterion given in [24], Theorem 3.4.16, to this estimate for large enough $p$, this ensures that for every $\alpha < 1/2$, the family $(n^{-1/2}\widehat{H}^T_{(2n-2)\cdot})$ is uniformly Hölder continuous under $P(\cdot|\#T = n)$ with exponent $\alpha$ (we write $\alpha$-UHC), that is, for every $\varepsilon > 0$ there exists a finite real number $C_\varepsilon$ such that, for every $n$,

$$P\left(\left|\frac{\widehat{H}^T_{(2n-2)s} - \widehat{H}^T_{(2n-2)t}}{\sqrt{n}}\right| \leq C_\varepsilon |t-s|^\alpha \text{ for all } s,t \in [0,1] \Big| \#T = n\right) \geq 1 - \varepsilon.$$

On the other hand, by a slight adaptation of the second proof of Lemma 1 in [13], we get that the proposition holds with the hypothesis $\alpha < 1/4$ instead of $\alpha < 1/2$. Indeed, the argument given in [13], which deals with the contour process, entirely rests on an exponential inequality linking this process to the so-called depth-first walk, and according to [18], Theorem 2, this inequality is also satisfied for the height process instead of the contour.

We now argue that if $(n^{-1/2}\widehat{H}^T_{(2n-2)\cdot})$ is $\beta$-UHC and $(n^{-1/2}H^T_{(n-1)\cdot})$ is $\alpha$-UHC for any $\beta < 1/2, \alpha < 1/4$, then $(n^{-1/2}H^T_{(n-1)\cdot})$ is $\alpha$-UHC for any $\alpha < 1/2$, which will end the proof. Assume that $n^{-1/2}|\widehat{H}^T_{(2n-2)s} - \widehat{H}^T_{(2n-2)t}| \leq c_1|t-s|^{1/2-a}$ and $n^{-1/2}|H^T_{(n-1)s} - H^T_{(n-1)t}| \leq c_2|t-s|^{\alpha-b}$, for any $s,t \in [0,1]$, and for some $\alpha < 1/2$ (this is true for $\alpha = 1/4$, and any $a,b > 0$ with probability close to 1, for some $c_1$ and $c_2$). Now, let $s$ and $t$ be such that $(n-1)s$ and $(n-1)t$ are two different integers. We have

$$\left|\frac{H^T_{(n-1)s} - H^T_{(n-1)t}}{\sqrt{n}}\right| = \left|\frac{\widehat{H}^T_{m_T((n-1)s)} - \widehat{H}^T_{m_T((n-1)t)}}{\sqrt{n}}\right|$$

$$\leq c_1 \left|\frac{m_T((n-1)s) - m_T((n-1)t)}{2n-2}\right|^{1/2-a}.$$

By (53), this is smaller than

$$c_1\left|s - t + \frac{H^T_{(n-1)t} - H^T_{(n-1)s}}{2n-2}\right|^{1/2-a} \leq c_1\left||s-t| + c_2\frac{|t-s|^{\alpha-b}}{\sqrt{n}}\right|^{1/2-a}$$

$$\leq c_1||t-s| + c_2|t-s|^{\alpha-b+1/2}|^{1/2-a},$$

since $n^{-1/2} \leq |t-s|^{1/2}$. For $s,t \in [0,1]$, $|t-s| \leq |t-s|^{\alpha-b+1/2}$, and then $(n^{-1/2}H^T_{(n-1)\cdot})$ is $(\alpha-b+1/2)(1/2-a)$-UHC. Since this holds for any $a > 0$ and $b > 0$, and since $\phi: \alpha \mapsto (\alpha + 1/2)1/2$ is increasing and contracting, $(n^{-1/2} H^T_{(n-1)}.)$ is $c$-UHC for any $c$ smaller than the fixed point of $\phi$ which is $1/2$. $\square$

PROPOSITION 22. *Let $(\mu_0, \mu_1)$ be a critical nondegenerate offspring distribution that admits some exponential moments. Let $i, j \in \{0, 1\}$. Then for every $\varepsilon > 0, \alpha \in (0, 1/2)$, there exists $C > 0$ with*

$$\sup_{n \in \mathbb{N}} P^{(i)} \left( \sup_{0 \leq s \neq t \leq 1} \frac{|H^T_{(\#T-1)s} - H^T_{(\#T-1)t}|}{\sqrt{n}|s-t|^\alpha} > C \,\Big|\, \#T^{(j)} = n \right) \leq \varepsilon.$$

PROOF. We prove this only for $i = 0, j = 1$, which is the hardest case of both. We fix $\alpha \in (0, 1/2)$. We assume that $(\#T - 1)s$ and $(\#T - 1)t$ are integer. Recall the notation $J^{(1)}_T(k), 0 \leq k \leq \#T - 1$, and extend this into a function $(J^{(1)}_T(t) = J^{(1)}_T([t]), 0 \leq t \leq \#T - 1)$. We bound

$$
\begin{aligned}
\left| \frac{H^T_{(\#T-1)s} - H^T_{(\#T-1)t}}{\sqrt{n}} \right| &\leq \left| \frac{H^T_{(\#T-1)s} - 2H^{\Gamma'(1T)}_{J^{(1)}_T((T-1)s)-1}}{\sqrt{n}} \right| \\
&\quad + 2 \left| \frac{H^{\Gamma'(1T)}_{J^{(1)}_T((\#T-1)s)-1} - H^{\Gamma'(1T)}_{J^{(1)}_T((\#T-1)t)-1}}{\sqrt{n}} \right| \\
&\quad + \left| \frac{H^T_{(\#T-1)t} - 2H^{\Gamma'(1T)}_{J^{(1)}_T((\#T-1)t)-1}}{\sqrt{n}} \right|.
\end{aligned}
\tag{54}
$$

Recall from the proof of Proposition 17 that the law of $\pi(1T)$ under $P^{(0)}(\cdot | \#T^{(1)} = n)$ is a mixture of the form $d\lambda_n(r) P^{(1)}_r(d\mathbf{t} | \#F^{(1)} = n)$, where the laws $\lambda_n, n \geq 1$ are tight, and that $\Gamma'(1T)$ is a monotype GW forest with a $\lambda_n$-distributed number of tree components and conditioned to have $n$ vertices. Therefore, by Proposition 21, with probability $> 1 - \varepsilon$ and for some $C > 0$, the middle term of (54), is bounded by

$$C \left| \frac{J^{(1)}_T((\#T-1)s) - J^{(1)}_T((\#T-1)t)}{n-1} \right|^\alpha \leq C \left| \frac{(\#T-1)s - (\#T-1)t}{n-1} \right|^\alpha,$$

since $J^{(1)}$ is a counting process. This is $\leq C'|t-s|^\alpha$ with probability $1 - 2\varepsilon$ for some $C' > 0$, valid for all $n$ large (by Lemma 15).

Next, using (47), the two other terms of (54) are bounded above by a constant multiple of

$$n^{-\alpha} \sup_{0 \leq k \leq \#T-1} n^{\alpha - 1/2} |H^{\Gamma'(1T)}_k - H^{\Gamma'(1T)}_{k+1}|,$$

with the convention that the second term in the absolute value is 0 if $k = \#T - 1$. By Lemma 18, under $P^{(0)}(\cdot|\#T^{(1)} = n)$, the quantity in the supremum converges to 0 in probability. Thus, for every $n$ large, and $s \neq t$ such that $(\#T-1)s, (\#T-1)t$ are integers, we have under $P^{(0)}(\cdot|\#T^{(1)} = n)$, fixing $\varepsilon > 0$, for all large $n$ and with probability $\geq 1 - \varepsilon$,

$$\left|\frac{H^T_{(\#T-1)s} - H^T_{(\#T-1)t}}{\sqrt{n}}\right| \leq n^{-\alpha} + C'|s-t|^\alpha$$
$$\leq D^{-\alpha}(\#T-1)^{-\alpha} + C'|s-t|^\alpha$$
$$\leq C''|s-t|^\alpha,$$

where we have used that with high probability, $\#T \leq Dn$ for some constant $D > 0$, see the proof of Lemma 15, and the fact that $(\#T-1)^{-1} \leq |s-t|$ for our choice of $s, t$. Finally, this shows the result for all $s, t$ in $\{k(\#T-1)^{-1}, k \in \mathbf{Z}_+\}$ and large $n$, and the result follows from the following elementary lemma, and then taking $C''$ even larger to fit to all $n \geq 1$. The case $i = j$, which we leave to the reader, is similar but easier since it makes use only of trees and the mapping $\Gamma$, rather than forests and the mapping $\Gamma'$. □

LEMMA 23. *Let $\alpha \in (0,1)$. If $f(k/n), 0 \leq k \leq n$ satisfies $|f(k/n) - f(k'/n)| \leq C((k-k')/n)^\alpha$ for every $0 \leq k, k' \leq n$, then the linear interpolation $(f(t), 0 \leq t \leq 1)$ satisfies $|f(t) - f(s)| \leq 3C|s-t|^\alpha$ for every $0 \leq s, t \leq 1$.*

5.3. *Tightness of the label process.* The first step of the proof of Theorem 8 is:

PROPOSITION 24. *Under the hypotheses of Theorem 8, for $i, j \in \{0, 1\}$, the sequence of laws of the processes $(n^{-1/4} S^{T,L}_{(\#T-1)s}, 0 \leq s \leq 1)$ under $\mathbb{P}^{(i)}(\cdot|\#T^{(j)} = n)$ for $n \geq 1$, is tight in $\mathcal{C}([0,1])$.*

PROOF. Our proof follows closely the arguments of [13, 19]. Fix $\varepsilon > 0$. Our goal is to show that there exists $C_1, \beta > 0$ such that for $n$ large enough,

$$(55) \quad \mathbb{P}^{(i)}\left(\sup_{0 \leq s,t \leq 1} \frac{|S^{T,L}_{(\#T-1)s} - S^{T,L}_{(\#T-1)t}|}{|t-s|^\beta} \leq C_1 n^{1/4} \middle| \#T^{(j)} = n\right) \geq 1 - \varepsilon.$$

Since the moment condition (20) is satisfied, there are constants $\eta, C_2, D > 0$ such that $M_0^k \vee M_1^k \leq C_2 k^D$ for every $k \geq 1$, where $M_c^k = \langle \nu_c^k, |x|^{4+\eta}\rangle$. We first choose $\alpha < 1/2$ so that $\alpha(4 + \eta) > 2$ and $M > D + 2$. We know from Proposition 22 that there exists $C_3 > 0$ such that

$$(56) \quad \mathbb{P}^{(i)}\left(\sup_{0 \leq s,t \leq 1} \frac{|H^T_{(\#T-1)s} - H^T_{(\#T-1)s}|}{n^{1/2}|s-t|^\alpha} \leq C_3 \middle| \#T^{(j)} = n\right) \geq 1 - \varepsilon,$$

for all $n$. Let $B_n$ be the intersection of the corresponding event and of the events

$$\left\{\max_{u\in T} c_T(u) \leq n^\gamma\right\}, \qquad \{\#T \leq C_4 n\}$$

and

$$\left\{\max_{u\in T}\max_{k\geq 1, 1\leq l\leq k}\max_{h\geq n^\gamma}\max_{c\in\{0,1\}} \frac{|A_T^{(c)}(u,k,l,h) - \mu_c(k)h/(2m_c)|}{h^{1/2+\gamma}k^{-M}} \leq 1\right\},$$

where $\gamma$ is such that $(D+3)\gamma < 1/2 - \alpha$, and $C_4 > 0$ is chosen so that the probability $P^{(i)}(B_n|\#T^{(j)} = n)$ is $\geq 1 - \varepsilon$ for large $n$, which is possible by Lemmas 13, 15 and Proposition 20.

Notice that by definition of $S^{T,L}$ and Lemma 23, it suffices to show (55) for all $n$ large and $s \neq t$ such that $(\#T - 1)s$ and $(\#T - 1)t$ are integers, which we suppose from now on. We let $m = (\#T - 1)s, m' = (\#T - 1)t$, and $u = u(m+1), u' = u(m'+1)$.

By definition, $S_m^{T,L} - S_{m'}^{T,L} = L(u) - L(u')$. If we let $Y_v = L(v) - L(\neg v)$, and if $\check{u} = \check{u}(m, m')$ denotes the most recent common ancestor to $u$ and $u'$ (i.e., their longest common prefix), we have

$$S_m^{T,L} - S_{m'}^{T,L} = \sum_{v\vdash u, |v|>|\check{u}|} Y_v - \sum_{v\vdash u', |v|>|\check{u}|} Y_v.$$

It is then classical that the number $R_T(m, m') = |u| + |u'| - 2|\check{u}(m, m')|$ of terms involved in these two sums, which informally is the length of the path of edges going from $u$ to $u'$ in $T$, is bounded according to the formula

(57) $$|R_T(m, m') - |H_m^T + H_{m'}^T - 2\check{H}_{m,m'}^T|| \leq 2,$$

where $\check{H}_{m,m'}^T$ is the infimum of $H^T$ between the points $m$ and $m'$. Indeed, if $u$ is an ancestor of $u'$ or conversely, then this expression is exactly $|H_m^T - H_{m'}^T| = ||u| - |u'||$, and otherwise, assuming $m < m'$, $\check{H}_{m,m'}^T$ is equal to the height of the first child $v$ of $\check{u}$ such that $u' \in vT_v$, hence is $|\check{u}| + 1$. In particular, it holds that for large enough $n$ and some $C_5 > 0$, on $B_n$,

(58) $$R_T(m, m') \leq 2 + C_3\sqrt{n}(|s - r|^\alpha + |r - t|^\alpha) \leq C_5\sqrt{n}|s - t|^\alpha,$$

where in the intermediate step $r \in \{k(\#T - 1)^{-1}, k \in \{0, \ldots, \#T - 1\}\}$ lies between $s$ and $t$ and is such that $H_{(\#T-1)r}^T = \check{H}_{m,m'}^T$. We also used the fact that $2 \leq \sqrt{n}|s - t|^\alpha$ for large $n$, since $|s - t| \geq (\#T - 1)^{-1} \geq (C_4 n)^{-1}$ under $B_n$.

Recall that under $\mathbb{P}^{(i)}$, given $T$, the increments of the label process $Y_u = L(v) - L(\neg v), v \in T$ (with the convention $Y_\varnothing = 0$), are such that $(Y_{v1}, \ldots, Y_{vc_T(v)})$, $v \in T$, are independent with respective laws $\nu_{|v|+i}^{c_T(v)}$. By splitting the involved sums according to the shape of $T$, we obtain, whenever $u = \check{u}l(m)w$ and $u' = \check{u}l(m')w'$ for $l(m), l(m') \in \mathbb{N}$ and $w, w' \in \mathcal{U}$,

$$
\begin{aligned}
(59) \quad & S_m^{T,L} - S_{m'}^{T,L} \\
&= (Y_{\check{u}l(m)} - Y_{\check{u}l(m')}) \\
&\quad + \sum_{k \geq 1} \sum_{1 \leq l \leq k} \sum_{c \in \{0,1\}} \sum_{v \vdash u, v \neq u, v \in T^{(c)}} Y_{vl} \mathbb{1}_{\{|v| > |\check{u}|, c_T(v) = k, u \in vlT_{vl}\}} \\
&\quad - \sum_{k \geq 1} \sum_{1 \leq l \leq k} \sum_{c \in \{0,1\}} \sum_{v \vdash u', v \neq u', v \in T^{(c)}} Y_{vl} \mathbb{1}_{\{|v| > |\check{u}|, c_T(v) = k, u' \in vlT_{vl}\}}.
\end{aligned}
$$

Notice that the last sum of the second line has $A_T^{(c)}(u, k, l, h)$ terms [resp. $A_T^{(c)}(u', k, l, h')$ in the third line] where $h = |u| - |\check{u}| - 1$ (resp. $h' = |u'| - |\check{u}| - 1$). Moreover, all the terms involved in the two last lines of (59) are independent and independent of the terms of the first line, with respective laws the $l$th marginal of $\nu_c^k$. The only two terms that bear some dependence are the ones displayed on the first line.

We now use an inequality due to Rosenthal [21], Theorem 2.10, which states that if $X_1, \ldots, X_n$ are independent and centered (but not necessarily identically distributed) under some probability law $\widetilde{P}$, then there exist universal constants $C(p), p \geq 2$, such that

$$
(60) \quad \widetilde{E}[|X_1 + \cdots + X_n|^p] \leq C(p) n^{p/2-1} \sum_{k=1}^{n} \widetilde{E}[|X_k|^p].
$$

This gives, still denoting $h = |u| - |\check{u}| - 1, h' = |u'| - |\check{u}| - 1$, and for $p = 4 + \eta$,

$$
\mathbb{E}^{(i)}[|S_m^{T,L} - S_{m'}^{T,L}|^p | T]
$$

$$
\leq C(p) R_T(m, m')^{p/2-1}
$$

$$
(61) \quad \times \begin{pmatrix} \mathbb{E}^{(i)}[|Y_{\check{u}l(m)} - Y_{\check{u}l(m')}|^p | T] \\ + \displaystyle\sum_{1 \leq k \leq \max_{u \in T} c_T(u)} \sum_{1 \leq l \leq k} \sum_{c \in \{0,1\}} A_T^{(c)}(u, k, l, h) \langle \nu_c^k, |x_l|^p \rangle \\ + \displaystyle\sum_{1 \leq k \leq \max_{u \in T} c_T(u)} \sum_{1 \leq l \leq k} \sum_{c \in \{0,1\}} A_T^{(c)}(u', k, l, h') \langle \nu_c^k, |x_l|^p \rangle \end{pmatrix}
$$

$$
\leq C(p) C_2 R_T(m, m')^{p/2-1}
$$

$$
\times \begin{pmatrix} 2^p c_T(\check{u})^D \\ + \displaystyle\sum_{1 \leq k \leq \max_{u \in T} c_T(u)} k^D \sum_{1 \leq l \leq k} \sum_{c \in \{0,1\}} A_T^{(c)}(u, k, l, h) \\ + \displaystyle\sum_{1 \leq k \leq \max_{u \in T} c_T(u)} k^D \sum_{1 \leq l \leq k} \sum_{c \in \{0,1\}} A_T^{(c)}(u', k, l, h') \end{pmatrix}
$$

which we want to bound on $B_n$. On the latter event, we have $c_T(\check{u})^D \leq n^{D\gamma}$, and by (58) it holds that for $n$ large, and every $s,t$ (satisfying the above constraints) the quantity $2^p C(p) C_2 R_T(m,m')^{p/2-1} c_T(\check{u})^D$ is bounded by

(62) $C_6 n^{p/4-1/2+D\gamma} |s-t|^{\alpha p/2-\alpha} \leq C_6 n^{p/4} |s-t|^{\alpha p/2} n^{D\gamma-1/2} |s-t|^{-\alpha},$

where $C_6 = 2^p C(p) C_2 C_5^{p/2-1}$. On the event $B_n$, and by the assumption on $s,t$, $|s-t| \geq \#T^{-1} \geq (C_4 n)^{-1}$, so if we combine this with the fact that $D\gamma - 1/2 < -\alpha$, we obtain that there exists $\varepsilon' > 0$ with $n^{D\gamma-1/2}|s-t|^{-\alpha} \leq C_4^\alpha n^{-\varepsilon'}$. Since this quantity converges to 0 uniformly on $s,t$, this shows that $2^p C(p) C_2 R_T(m,m')^{p/2-1} c_T(\check{u})^D$ is bounded by $n^{p/4}|s-t|^{\alpha p/2}$ for every $s,t$ and $n$ large enough.

We are now facing several possibilities in handling the rest of (61). On the event that $h \leq n^\gamma$, the term $A^{(c)}(u,k,l,h)$ is bounded by $n^\gamma$, and therefore, on $B_n$,

$$C(p)C_2 R_T(m,m')^{p/2-1} \sum_{1 \leq k \leq n^\gamma} k^D \sum_{1 \leq l \leq k} \sum_{c \in \{0,1\}} A_T^{(c)}(u,l,k,h)$$
$$\leq 2C(p)C_2 C_5^{p/2-1} n^{p/4-1/2}|s-t|^{\alpha p/2-\alpha} n^{(D+3)\gamma}.$$

This quantity is analogous to (62), and by our choice of $\gamma$, it is bounded by $n^{p/4}|s-t|^{\alpha p/2}$ for large $n$, by the same argument as above.

Alternatively, on the event that $h \geq n^\gamma$ and on $B_n$, we can bound $A_T^{(c)}(u,k,l,h)$ above by the quantity $h\mu_c(k)/(2m_c) + h^{1/2+\gamma}k^{-M}$. Since by definition $h \leq R_T(m,m')$, since $\mu_0, \mu_1$ have some exponential moments, and since by our choice of $M$ the sequence $k^{D-M+1}, k \geq 1$, is summable, it follows that there exists some constant $C_7 \in (0, +\infty)$ such that on $B_n \cap \{h \geq n^\gamma\}$,

$$\sum_{k \leq n^\gamma} k^D \sum_{1 \leq l \leq k} \sum_{c \in \{0,1\}} A_T^{(c)}(u,l,k,h) \leq C_7 R_T(m,m') + C_7 R_T(m,m')^{1/2+\gamma}$$

$$\leq 2C_7 R_T(m,m'),$$

because $R_T(m,m')$ is an integer and $\gamma < 1/2$. Still on the event $h \geq n^\gamma$, it follows that the middle term of (61) is bounded by

$$2C(p)C_2 C_7 R_T(m,m')^{p/2} \leq 2C(p)C_2 C_7 C_5 n^{p/4}|s-t|^{\alpha p/2}.$$

Putting things together, we obtain the existence of a constant $C_8 > 0$ such that for every large $n$, the middle term of (61) is bounded by $C_8 n^{p/4}|s-t|^{\alpha p/2}$ on $B_n$. For the same reason, the third term of (61) is bounded by the same quantity, and from the discussion on the term involving $c_T(\check{u})^D$ we finally obtain that for some $C_9 > 0$ and $n$ large, for every $s,t \in \{k(\#T-1)^{-1}, k \in \{0,\ldots,\#T-1\}\}$,

$$\mathbb{E}^{(i)}\left[\left(\frac{|S_m^{T,L} - S_{m'}^{T,L}|}{n^{1/4}}\right)^p \Big| \{\#T^{(j)} = n\} \cap B_n\right] \leq C_9 |s-t|^{\alpha p/2}.$$

By applying Lemma 23, and since $\alpha p/2 > 1$, we have obtained that there exists a finite constant $C_{10} > 0$ and $\eta' > 0$ such that for every $0 \leq s, t \leq 1$ and $n$ large enough,

$$\mathbb{E}^{(i)}\left[\left(\frac{|S^{T,L}_{(\#T-1)s} - S^{T,L}_{(\#T-1)t}|}{n^{1/4}}\right)^{4+\eta}\bigg|\{\#T^{(j)} = n\} \cap B_n\right] \leq C_{10}|s-t|^{1+\eta'}.$$

By Kolmogorov's criterion [24], Theorem 3.4.16, this is enough to conclude that (55) holds with any $\beta \in (0, \eta'/(4+\eta))$, when replacing the conditioning event $\{\#T^{(j)} = n\}$ by $B_n \cap \{\#T^{(j)} = n\}$. Since $P^{(i)}(B_n|\#T^{(j)} = n) \geq 1 - \varepsilon$ for large $n$, we thus obtain (55) with lower bound $(1-\varepsilon)^2$ instead of $1-\varepsilon$. This ends the proof. $\square$

5.4. *Finite-dimensional convergence.* The goal of this section is to prove the following proposition.

PROPOSITION 25. *Let $(\mu_0, \mu_1)$ be nondegenerate critical and admitting some exponential moments. Consider nondegenerate spatial displacement laws $\nu_0^k, \nu_1^k, k \geq 1$, that are centered, and suppose that the hypotheses of Theorem 8 hold. Then the convergence of finite-dimensional marginals holds for the label process, that is, the second component in Theorem 8, jointly with the convergence in distribution of the height process.*

In fact, this statement is true under the slightly weaker hypothesis that the variance $\Sigma_0^k \vee \Sigma_1^k$ of the spatial displacement is a $O(k^D)$ for some $D > 0$, and does not really require the full $4 + \eta$ moment hypothesis of Theorem 8. However, this extra assumption is going to simplify the beginning of the proof, where we use the result of Proposition 24. Notice that the constant $\Sigma$ of Theorem 8 is finite because of the hypothesis on the growth of the moments of order $4 + \eta$ implies that the growth of $\Sigma_i^k$ is at most polynomial, and since $\mu_0, \mu_1$ have some exponential moments, the sums $\sum_{k \geq 1} (\Sigma_i^k)^2 \mu_i(k)$ are finite.

The intuition for the proof of one-dimensional convergence is the following. Given the height process $H^T$, if we take the path from the root to the vertex encoded by a given time $t$, then this path has length of about $H^T_{(\#T-1)t}$ which is of the order of $h = C n^{1/2} B_t$, where $B$ is a Brownian excursion and $C$ is a scaling constant. Among the vertices of this ancestor line, we know from Lemma 20 that a proportion $\mu_j(k)/(2m_j)$ are of type $j$, have $k$ children, and have the property that $u$ is a descendant of the $l$th of these children; these will contribute to a spacial displacement whose distribution is the $l$th marginal of $\nu_j^k$. Since the variance of this is $(\Sigma_j^{k,l})^2$, it has to be expected that the total spatial displacement, once rescaled by $n^{1/4}$ will be asymptotically Gaussian with variance $\Sigma^2 C B_t$, where $\Sigma$ is defined at (19).

PROOF OF PROPOSITION 25. We make a preliminary remark. Because we already know that the laws of $(n^{-1/4} S^{T,L}_{(\#T-1)t}, 0 \leq t \leq 1)$ under $P^{(i)}(\cdot | \#T^{(j)} = n)$ form a tight family, the family of laws of $((n^{-1/2} H^T_{(\#T-1)t})_{0 \leq t \leq 1}, (n^{-1/4} S^{T,L}_{(\#T-1)t})_{0 \leq t \leq 1})$ under $P^{(i)}(\cdot | \#T^{(j)} = n)$ is also tight. Hence, up to extracting a subsequence, we know that these two processes jointly converge to some limit, whose first component is a scaled Brownian excursion thanks to Proposition 17. To prove the proposition, it suffices to show that the only possible limiting distribution is the (properly scaled) head of the Brownian snake $\mathbb{N}^{(1)}$. So, we take such a subsequence in the first place, assume convergence in distribution to some process $(B, S')$, and our goal is to show that given $B$, $S'$ has the law described around (4) up to scaling constants. Let $0 < t_1 < t_2 < \cdots < t_q < 1$ be some fixed real numbers. We will prove that $(n^{-1/2} H^T_{(\#T-1)t})_{0 \leq t \leq 1}, (n^{-1/4} S^{T,L}_{[(\#T-1)t_r]})_{1 \leq r \leq q})$ under $P^{(i)}(\cdot | \#T^{(j)} = n)$ converges in distribution to the corresponding marginal of the head of Brownian snake, which is sufficient to conclude. Throughout the proof, we will assume $q \geq 2$.

Thanks to Skorokhod's representation theorem, we may assume that the convergence of the processes $(n^{-1/2} H^T_{(\#T-1)t}, 0 \leq t \leq 1)$ to a scaled Brownian excursion is almost sure. That is, we can assume that we are given a sequence $(H^n, n \geq 1)$ of processes on some probability space $(\Omega, \mathcal{F}, P)$, with same respective laws as $(n^{-1/2} H^T_{(\#T-1)t}, 0 \leq t \leq 1)$ under $P^{(i)}(\cdot | \#T^{(j)} = n)$, and which converges a.s. for the supremum norm to a process $(B_t, 0 \leq t \leq 1)$, which has same law as $2\sigma^{-1}\sqrt{1+m_j} e$ under $\mathbb{N}^{(1)}$. For every $n \geq 1$, the function $H^n$ determines a unique random tree $T^n$ whose height process is $\sqrt{n} H^n((\#T^n - 1)^{-1}k), 0 \leq k \leq \#T^n - 1$ [notice that the renormalization constant $\#T^n - 1$ can be recovered from $n$ and $H^n$, as $(\#T^n - 1)n^{-1/2}$ is the slope of $H_n$ at $0+$]. By Lemma 20 and the Borel–Cantelli lemma, it holds that, a.s., for any $\gamma > 0$ and $n$ large enough,

$$(63) \quad \max_{c \in \{0,1\}} \sup_{\{k \geq 1, 1 \leq k \leq l\}} \sup_{\{u \in T^n, n^\gamma \leq h \leq |u|\}} \frac{|A^{(c)}_{T^n}(u,k,l,h) - \mu_c(k)h/(2m_c)|}{h^{1/2+\gamma} k^{-M}} \leq 1.$$

Next, the times $t_1, \ldots, t_q$ determine vertices $u^r_n = u([(\#T^n - 1)t_r]), 1 \leq r \leq q$, in $T^n$. We let $\check{u}_n(r, r')$ be the most recent common ancestor to $u^r_n, u^{r'}_n$ in $T^n$. We re-index the set $\mathcal{V}^n = \{u^r_n, 1 \leq r \leq q, \check{u}_n(r, r'), 1 \leq r, r' \leq q\}$ as $\{v^w_n, w \in \widetilde{T}^n\}$, where $\widetilde{T}^n \in \mathcal{T}$, and in such a way that the depth-first order and genealogical structure on $\widetilde{T}^n$ is compatible with the depth-first order and genealogical structure on $\mathcal{V}^n$. Specifically, we let $v^\varnothing_n$ be the least element of $\mathcal{V}^n$, which is the most recent common ancestor to all of $u^1_n, \ldots, u^q_n$, then recursively, $v^{w1}_n, v^{w2}_n, \ldots$ are the descendants of $v^w_n$ in $\mathcal{V}^n$, ranked in depth-first order, and such that no ancestor of $v^{wl}_n$ which is younger than $v^w_n$ belongs to $\mathcal{V}^n$. By convention, we let $v^{\neg \varnothing}_n = \varnothing$. Our aim is now to explain that the sequence $(\widetilde{T}_n, n \geq 1)$ is asymptotically constant a.s., and equal to some random tree, which has $q$ leaves and is binary, that is, vertices have either no child or two children. Informally, this implies that for large enough $n$, the "geometry" of the list of random variables involved in the computations of the labels of the vertices in $\mathcal{V}^n$ is eventually fixed.

We define random times $s_n^w, s^w, w \in \widetilde{T}^n$ recursively as follows (these times will be defined only for $n$ large enough). Let

$$s_n^\varnothing = \inf\left\{ t \in [t_1, t_q] : H_t^n = \min_{t_1 \leq s \leq t_q} H_s^n \right\},$$

and let $\overline{s}_n^\varnothing$ be the corresponding quantity with a sup in place of the inf. Notice that, $H_{s_n^\varnothing}^n$ converges to $\min_{t_1 \leq s \leq t_q} B_s > 0$ a.s. Since local minima of a Brownian excursion are a.s. pairwise distinct, it is then elementary to deduce that $s_n^\varnothing, \overline{s}_n^\varnothing \to s^\varnothing$ where by definition $s^\varnothing$ is the only time in $[t_1, t_q]$ where $B$ attains $\min_{t_1 \leq s \leq t_q} B_s$. Moreover, a.s., $t_1 < \cdots < t_r < s^\varnothing < t_{r+1} < \cdots < t_q$ for some $1 \leq r \leq q-1$, and for $n$ large enough it also holds that $t_1 < \cdots < t_r < s_n^\varnothing \leq \overline{s}_n^\varnothing < t_{r+1} < \cdots < t_q$. We let $r = r_1, q - r = r_2$, $t_1^1 = t_1, \ldots, t_{r_1}^1 = t_r$ and $t_1^2 = t_{r+1}, \ldots, t_{r_2}^2 = t_q$.

Then, given

$$s_n^w, \overline{s}_n^w, s^w, (t_1^{w1}, \ldots, t_{r_{w1}}^{w1}), (t_1^{w2}, \ldots, t_{r_{w2}}^{w2}),$$

have been defined, where $w$ is a word with letters in $\{1, 2\}$, we distinguish two cases. If $r_{w1} > 1$, let

$$s_n^{w1} = \inf\left\{ t \in [t_1^{w1}, t_{r_{w1}}^{w1}] : H_t^n = \min_{t_1^{w1} \leq s \leq t_{r_{w1}}^{w1}} H_s^n \right\},$$

and let $\overline{s}_n^{w1}$ be the corresponding quantity with a sup in place of the inf. For the same reasons as above it holds that $s_n^{w1}, \overline{s}_n^{w1} \to s^{w1}$ where $s^{w1}$ is the only time in $[t_1^{w1}, t_{r_{w1}}^{w1}]$ where $B$ attains $\min_{t_1^{w1} \leq s \leq t_{r_{w1}}^{w1}} B_s$. Moreover, a.s., $t_1^{w1} < \cdots < t_r^{w1} < s^{w1} < t_{r+1}^{w1} < \cdots < t_{r_{w1}}^{w1}$ for some $1 \leq r \leq r_{w1} - 1$, and for $n$ large enough it also holds that $t_1^{w1} < \cdots < t_r^{w1} < s_n^{w1} \leq \overline{s}_n^{w1} < t_{r+1} < \cdots < t_{r_{w1}}^{w1}$. We let

$$r_{w11} = r, \qquad r_{w12} = r_{w1} - r, \qquad t_1^{w11} = t_1^{w1}, \ldots, t_{r_{w11}}^{w11} = t_r^{w1}$$

and

$$t_1^{w12} = t_{r+1}^{w1}, \ldots, t_{r_{w12}}^{w12} = t_{r_{w1}}^{w1}.$$

Definitions are similar for $s_n^{w2}, \overline{s}_n^{w2}, s^{w2}, (t_1^{w21}, \ldots, t_{r_{w21}}^{w21}), (t_1^{w22}, \ldots, t_{r_{w22}}^{w22})$ whenever $r_{w2} > 1$. In the case $r_{w1} = 1$, we simply let $s_n^{w1} = (\#T^n - 1)^{-1}[(\#T^n - 1)t_1^{w1}]$, and $s^{w1} = t_1^{w1}$, and similarly if $r_{w2} = 1$.

By inspection of this recursive construction, notice that the set of $w \in \mathcal{U}$ such that $s_n^w$ is defined is exactly $\widetilde{T}^n$, which is therefore independent of $n$ (provided $n$ is large enough), and equal to some binary tree $\widetilde{T}$. Moreover, it holds that $||v_n^w| - \sqrt{n}H_{s_n^w}^n| \leq 1$ for every $w \in \widetilde{T}$, by the same arguments as those leading to (57).

Now, we reintroduce the labels in $T^n$ by assuming that $(\Omega, \mathcal{F}, P)$ also supports random variables $(Y_{vl}^n, 1 \leq l \leq c_{T^n}(v)), v \in T^n$, used as spatial displacements, that are, conditionally on $T^n$, independent and independent of $B$ with respective laws $\nu_{i+|v|}^{c_{T^n}(v)}$. We let $L_n$ be the associated label on $T^n$ with $L_n(\varnothing) = 0$, and we use a truncation procedure, that is we choose $C$ large and write [remembering (59)]

$$L_n(u) = L_n^C(u) + \widetilde{L}_n^C(u), \qquad u \in T^n,$$

where

$$L_n^C(u) = \sum_{v \vdash u} Y_v^n \mathbb{1}_{\{|v|>0, c_T^n(v) \leq C\}},$$

and $\widetilde{L}_n^C(u)$ is the similar sum with $c_T^n(v) > C$ instead. Then, the random variables $L_n^C(v_n^w) - L_n^C(v_n^{\neg w}), w \in \widetilde{T}$ can be written in the form

$$Y_{v_n^{\neg w}l(w)}^n \mathbb{1}_{\{c_{T^n}(v_n^{\neg w}) \leq C\}} + \sum_{v \vdash v_n^w, |v| > |v_n^{\neg w}|+1} Y_v^n \mathbb{1}_{\{c_{T^n}(v) \leq C\}},$$

whenever $l(w) \in \mathbb{N}$ is such that $v_n^w \in v_n^{\neg w} l(w) T_{v_n^{\neg w}l(w)}^n$. Notice that given $T^n$, all these terms are independent as $w$ ranges in $\widetilde{T}$, except maybe for the first term which is displayed to the left of the sum. Since we rescale by $n^{-1/4}$, this term disappears in the limit (note that its variance is bounded) so that $(n^{-1/4}(L_n^C(v_n^w) - L_n^C(v_n^{\neg w})), w \in \widetilde{T})$, has same limit as the vector of independent components (given $T^n$)

$$\left(n^{-1/4} \sum_{v \vdash v_n^w, |v| > |v_n^{\neg w}|+1} Y_v^n \mathbb{1}_{\{c_{T^n}(v) \leq C\}}, w \in \widetilde{T}\right)$$

(64) $= \left(n^{-1/4} \sum_{k=1}^{C} \sum_{l=1}^{k} \sum_{c \in \{0,1\}} \sum_{\{v \vdash v_n^w, |v| > |v_n^{\neg w}|+1\}} Y_v^n \mathbb{1}_{\{v \in T^{n(c)}, c_{T^n}(v)=k, v_n^w \in vlT_{vl}^n\}},\right.$

$$\left. w \in \widetilde{T}\right).$$

By definition there are $A_{T^n}^{(c)}(v_n^w, k, l, h_n^w)$ terms in the last sum, where $h_n^w = |v_n^w| - |v_n^{\neg w}| - 1$. Since $n^{-1/2}h_n^w$ has same limit as $|H_{s_n^w}^n - H_{s_n^{\neg w}}^n|$ as $n \to \infty$, which is given by $B_{s^w} - B_{s^{\neg w}}$, and which is $> 0$ a.s., we obtain that asymptotically $h_n^w > n^\gamma$ for any fixed $0 < \gamma < 1/2$. Therefore, by (63) it

holds that a.s., for any $\varepsilon > 0$, any $k \leq C, l \leq k$ and for $n$ large

$$(1-\varepsilon)\frac{\mu_c(k)}{2m_c}(B_{s^w} - B_{s^{\neg w}})n^{1/2} \leq A_{T^n}^{(c)}(v_n^w, k, l, h_n^w)$$

$$\leq (1+\varepsilon)\frac{\mu_c(k)}{2m_c}(B_{s^w} - B_{s^{\neg w}})n^{1/2}.$$

It then follows from the central limit theorem applied to (64) that given $B$ (of which $\widetilde{T}$ is a measurable functional), the vector $n^{-1/4}(L_n^C(v_n^w) - L_n^C(v_n^{\neg w}), w \in \widetilde{T})$ converges in distribution to a random vector $(N^C(w), w \in \widetilde{T})$, where the components $N^C(w)$ are independent, centered, Gaussian and have variances

$$\operatorname{Var} N^C(w) = (B_{s^w} - B_{s^{\neg w}})\frac{1}{2}\sum_{k=1}^{C}\sum_{l=1}^{k}\sum_{c\in\{0,1\}}(\Sigma_c^{k,l})^2\frac{\mu_c(k)}{m_c}$$

$$=: (B_{s^w} - B_{s^{\neg w}})\Sigma_C^2.$$

Notice that $\Sigma_C^2 \uparrow \Sigma^2$ as $C \uparrow \infty$, so $(N^C(w), w \in \widetilde{T})$ in turn converges in distribution to a random vector which conditionally on $B$ is constituted of independent Gaussian components $(N(w), w \in \widetilde{T})$, with respective variances $\Sigma^2(B_{s^w} - B_{s^{\neg w}})$ as $C \to \infty$.

Assume for a moment that for every $\varepsilon > 0$, a.s.,

(65)   $$\lim_{C\to\infty}\limsup_{n\to\infty} P\left(\max_{w\in\widetilde{T}}|\widetilde{L}_n^C(v_n^w) - \widetilde{L}_n^C(v_n^{\neg w})| > \varepsilon n^{1/4}\bigg|B\right) = 0.$$

From the fact that $L_n(v_n^w) - L_n(v_n^{\neg w}) = L_n^C(v_n^w) - L_n^C(v_n^{\neg w}) + \widetilde{L}_n^C(v_n^w) - \widetilde{L}_n^C(v_n^{\neg w})$, this implies that conditionally on $B$, $n^{-1/4}(L_n(v_n^w) - L_n(v_n^{\neg w}), w \in \widetilde{T})$ converges to $(N(w), w \in \widetilde{T})$. Indeed, it is an elementary exercise that if $X_n = X_n^C + Y_n^C \in \mathbb{R}^d$, where $X_n^C \to X^C$ as $n \to \infty$, $X^C \to X$ as $C \to \infty$, both in distribution, and $\lim_C \limsup_n P(|Y_n^C| > \varepsilon) = 0$, then $X_n \to X$ as $n \to \infty$ in distribution.

It follows that conditionally on $B$, the vector $n^{-1/4}(L_n(u_n^r), 1 \leq r \leq q)$ is asymptotically a Gaussian vector $(S(t_1), \ldots, S(t_q))$, since

$$L_n(u_n^r) = \sum_{w' \vdash w}(L_n(v_n^{w'}) - L_n(v_n^{\neg w'})),$$

whenever $u_n^r = v_n^w$. Moreover, still given $B$, we have that if $t_r = s^w, t_{r'} = s^{w'}$, then

$$\operatorname{cov}(S_{t_r}, S_{t_{r'}}) = \operatorname{cov}\left(\sum_{w''\vdash w} N(w''), \sum_{w''\vdash w'} N(w'')\right).$$

By independence of the $N(w)$ given $B$, if $w \wedge w'$ is the most recent common ancestor to $w, w'$, we obtain

$$\operatorname{cov}(S_{t_r}, S_{t_{r'}}) = \operatorname{Var}\left(\sum_{w'' \vdash w \wedge w'} N(w'')\right) = \Sigma^2 B_{s^{w \wedge w'}} = \Sigma^2 \check{B}(t_r, t_{r'}),$$

as $s^{w \wedge w'}$ is the unique point of $[\min(s^w, s^{w'}), \max(s^w, s^{w'})]$ such that $B_{s^{w \wedge w'}} = \check{B}(s^w, s^{w'})$. Since $B$ has same law as $2\sigma^{-1}\sqrt{1+m_j}e$ under $\mathbb{N}^{(1)}$, it follows that $B, S$ has the claimed law.

To prove (65), notice that since given $B$, the set $\widetilde{T}$ is finite, it suffices to prove the result for some fixed $w \in \widetilde{T}$. Now, conditionally on $B, T^n$, the sequence $(\widetilde{L}_n^C(v), v \vdash v_n^w)$, where the ancestors of $v_n^w$ are ranked in depth-first order, has independent and centered increments so by Chebyshev's inequality,

$$P(|\widetilde{L}_n^C(v_n^w) - \widetilde{L}_n^C(v_n^{\neg w})| > \varepsilon n^{1/4}|B, T^n)$$
$$\leq n^{-1/2}\varepsilon^{-2} E[|\widetilde{L}_n^C(v_n^w) - \widetilde{L}_n^C(v_n^{\neg w})|^2 | B, T^n].$$

By the independence of increments of the spatial displacement,

$$E[|\widetilde{L}_n^C(v_n^w) - \widetilde{L}_n^C(v_n^{\neg w})|^2 | B, T^n]$$
$$= \sum_{k>C} \sum_{1 \leq l \leq k} \sum_{c \in \{0,1\}} A_{T^n}^{(c)}(v_n^w, k, l, h_n^w)(\Sigma_c^{k,l})^2,$$

where $h_n^w = |v_n^w| - |v_n^{\neg w}|$. Now, we know that $n^{-1/2}h_w^n$ converges to $B_{s^w} - B_{s^{\neg w}} > 0$ as $n \to \infty$, and therefore, by (63) for $\gamma = 1/8$ and $n$ large enough, it holds that the last expression is bounded by (for any $C_1 > B_{s^w}$)

$$C_1\sqrt{n} \sum_{k>C} \sum_{c \in \{0,1\}} \frac{\mu_c(k)}{2m_c}(\Sigma_c^k)^2 + n^{3/8} \sum_{k>C} \sum_{c \in \{0,1\}} k^{-M}(\Sigma_c^k)^2,$$

and the second term is bounded by $C_2 n^{3/8}$ for some constant $C_2 > 0$, because $\mu_c$ has some exponential moments, and $\Sigma_c^k = O(k^D)$ for some $D > 0$. Finally, we obtain that

$$P\left(\max_{w \in \widetilde{T}}|\widetilde{L}_n^C(v_n^w) - \widetilde{L}_n^C(v_n^{\neg w})| > \varepsilon n^{1/4}\Big| B\right)$$
$$\leq \#\widetilde{T} C_1 \varepsilon^{-2} \sum_{k>C} \sum_{c \in \{0,1\}} \frac{\mu_c(k)}{2m_c}(\Sigma_c^k)^2 + C_2 \varepsilon^{-2} n^{-1/8}.$$

Letting $n \to \infty$, this converges to

$$\#\widetilde{T} C_1 \varepsilon^{-2} \sum_{k>C} \sum_{c \in \{0,1\}} \frac{\mu_c(k)}{2m_c}(\Sigma_c^k)^2,$$

which in turn has limit 0 as $C \to \infty$, implying (65). $\square$

**6. Convergence to the Brownian map.** The aim of this section is to discuss the convergence of bipartite maps to the Brownian map, introduced in [20]. We refer to this paper to the construction of the notions of abstract maps, and to the combinatorial considerations leading the authors in a first step to show that any quadrangulation is a tree $\mathcal{D}$, (the doddering tree of [20]) glued with the help of a second tree $\mathcal{G}$ (the gluer tree of [20]), and to show that this construction passes to the continuous limit in a certain sense. The major part of the construction in [20] may be generalized without any problem to bipartite maps; in what follows, we will mainly point out the differences in the construction.

The first and major difference with the study of quadrangulation, is the use of the BDFG bijection, instead of Schaeffer's one.

*Bipartite maps described by a pair of trees.* We first present in a few words the application $\Psi^{-1}$ of Bouttier, Di Francesco and Guitter [6] (see also an illustration Figure 3).

Recall the considerations of Section 2.3, and consider a planar embedding of some $(\mathbf{t}, \ell) \in \overline{\mathbb{T}}$ in the plane, with at least two vertices. We let $v$ and $w$ be the vertices of this embedded graph that correspond to the words $\varnothing, 1$, and root the graph at the edge $(vw)$. Until the end of the paper, we slightly improperly keep the notations $\mathbf{t}, \mathbf{t}^{(0)}, \mathbf{t}^{(1)}, \ell$ for this embedded (rooted) graph, the vertices at even (resp. odd) heights, and a label function defined on $\mathbf{t}^{(0)}$ [the labels of $\mathbf{t}^{(1)}$ are not used in the construction].

The construction of $(\mathbf{m}, e, \mathfrak{r}) = \Psi^{-1}(\mathbf{t}, \ell)$ is done as follows. First add $-\min \ell + 1$ to the labels of the vertices of $\mathbf{t}^{(0)}$. Each vertex $u$ of $\mathbf{t}^{(0)}$ with $k$ neighbors determines $k$ *corners* which are delimited by the $k$ edges emanating from $u$. To each such corner $C$ corresponding to a vertex of $\mathbf{t}^{(0)}$ with label $l \geq 2$, we associate its successor $s(C)$ defined as the first encountered corner of $\mathbf{t}^{(0)}$ with label $l-1$ when going clockwise around the tree (there is always a successor). Then, we draw a blue edge between each corner of $\mathbf{t}^{(0)}$ with label $l \geq 2$ and its successor within the external face of $\mathbf{t}$ and in such a way that no two edges intersect, which is always possible. Finally, we add in complement of the graph in the plane, an extra vertex $\mathfrak{r}$, and add a blue edge between each the corners of $\mathbf{t}^{(0)}$ with label 1 and $\mathfrak{r}$. The map $(\mathbf{m}, e, \mathfrak{r})$ is the map having as set of edges the blue edges, and $e$ is the first edge of $\mathbf{m}$ that starts from $v$ to the left of $(vw)$, oriented so that $(\mathbf{m}, e, \mathfrak{r}) \in \mathcal{M}_+$.

Recall the definition of the depth first traversal $F_\mathbf{t}$ introduced in the proof of Proposition 21, and of the contour process $\widehat{H}^\mathbf{t}$. The durations of $F_\mathbf{t}$ and $\widehat{H}^\mathbf{t}$ is $2(\#\mathbf{t} - 1)$. The vertices of $\mathbf{t}^{(0)}$ are visited at times $0, 2, \ldots, 2(\#\mathbf{t} - 1)$. The labels of the vertices of $\mathbf{t}^{(0)}$ are encoded thanks to

$$\mathbf{R}^{\mathbf{t},\ell}(k) = \ell(F_\mathbf{t}(2k)) \qquad \text{for } 0 \leq k \leq \#\mathbf{t} - 1.$$

We extend $\mathbf{R}^{\mathbf{t},\ell}$ linearly between successive integers. We have COROLLARY 26. *Let $\mathbf{q}$ be a regular critical weight sequence, and let $(\mu_0, \mu_1)$ and $(\nu_0^k, \nu_1^k, k \geq 0)$ be the offspring distributions and spatial displacement laws that are associated with $\mathbf{q}$ as in Proposition 7. Then, the conclusion of Theorem 8 still holds with $n^{-1/4}\mathbf{R}^{T,L}((\#T-1)t)$ instead of $n^{-1/4}S^{T,L}_{(\#T-1)t}$, with the scaling constants given in Section 3.2.*

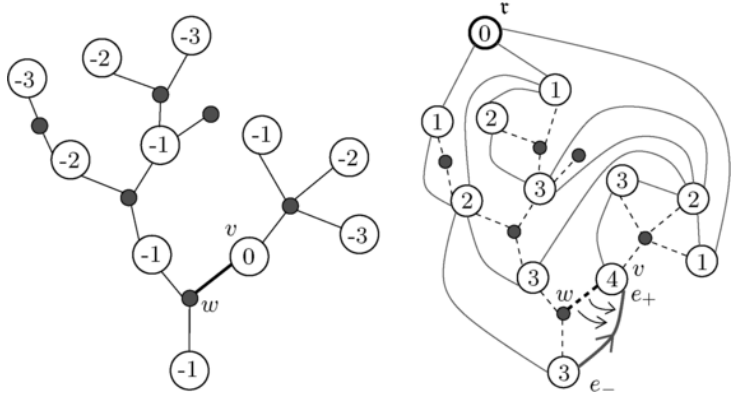

FIG. 3. *Illustration of $\Psi^{-1}$. The black vertices of $\mathbf{t}$ corresponds to the elements of $\mathbf{t}^{(1)}$, the labeled vertices to the elements of $\mathbf{t}^{(0)}$. The two arrows explain how to choose the root of $\mathbf{m}$. It remains to remove the labels, the black vertices and the dotted lines.*

This result is a consequence of two classical steps: firstly, let $R^{\mathbf{t},\ell}$ be the linear interpolation of $R^{\mathbf{t},\ell}(k) = \ell(F_\mathbf{t}(k))$, the label process associated with the depth first traversal [here take again the convention that a vertex of $\mathbf{t}^{(1)}$ has the same label as its father]. The uniform distance between $n^{-1/4}S^{T,L}_{(\#T-1)t}$ and $n^{-1/4}R^{T,L}(2(\#T-1)t)$ goes to 0 in probability. Secondly, for any integer $k$, $\mathbf{R}^{\mathbf{t},\ell}(k) = R^{\mathbf{t},\ell}(2k)$. This shows that the uniform distance between $n^{-1/4}\mathbf{R}^{T,L}((\#T-1)t)$ and $n^{-1/4}R^{T,L}(2(\#T-1)t)$ goes to 0 in probability.

For $\theta \in \{0, 2, \ldots, 2(\#\mathbf{t}-1)-2\}$, we denote by $\mathbf{t}_{(\theta)}$ the element of $\overline{\mathbb{T}}$ obtained from $\mathbf{t}$ by rerooting at the edge $(F_\mathbf{t}(\theta), F_\mathbf{t}(\theta+1))$, and with label function $\ell - \ell(F_\mathbf{t}(\theta))$. The label of the root-vertex of $\mathbf{t}_{(\theta)}$ is 0, and $\mathbf{t}$ and $\mathbf{t}_{(\theta)}$ are equal as unrooted unlabeled trees. Let $\oplus$ denote the addition modulo $2(\#\mathbf{t}-1)$. For any $0 \leq i \leq 2(\#\mathbf{t}-1)$,

$$\widehat{H}^{\mathbf{t}_{(\theta)}}(i) = \widehat{H}^\mathbf{t}(\theta \oplus i) + \widehat{H}^\mathbf{t}(\theta)$$
$$- 2\min\{\widehat{H}^\mathbf{t}(j), (\theta \oplus i) \wedge \theta \leq j \leq (\theta \oplus i) \vee \theta\},$$

and for any $0 \leq i \leq \#\mathbf{t}-1$,

$$\mathbf{R}^{\mathbf{t}_{(\theta)}}(i) = \ell(F_\mathbf{t}(\theta \oplus 2i)) - \ell(F_\mathbf{t}(\theta)) = \mathbf{R}^\mathbf{t}((\theta \oplus 2i)/2) - \mathbf{R}^\mathbf{t}(\theta/2).$$

When exploring $\mathbf{t}_{(\theta)}$, $v$ is visited at time $2(\#\mathbf{t}-1)-\theta$ and $w$ at time $2(\#\mathbf{t}-1)-\theta+1$. Hence, the variable $X(\theta) = 2\#\mathbf{t} - \theta$ suffices to recover $v$ and $w$.

We now exhibit the two trees from which the description of bipartite maps can be done.

Let $\Theta_{\mathbf{t}} = \inf\{\theta, \ell(F_{\mathbf{t}}(\theta)) = \min \ell(F_{\mathbf{t}})\}$, be the first visit time of a vertex with minimum label [this vertex is then a vertex of $\mathbf{t}^{(0)}$]. The labeled tree $\mathbf{t}_{(\Theta_{\mathbf{t}})}$ has nonnegative labels and $\Psi^{-1}(\mathbf{t})$ will be built from $(\Theta_{\mathbf{t}}, \mathbf{t}_{(\Theta_{\mathbf{t}})})$ (see Figure 4).

Add in the plane the point $u = N_{\#\mathbf{t}} = (\#\mathbf{t}, 0)$, and for $1 \leq i \leq \#\mathbf{t} - 1$ draw the vertex $N_i = (i, \mathbf{R}^{\mathbf{t}(\Theta)}(i) + 1)$. Then, for $1 \leq j \leq \#\mathbf{t} - 1$, an edge is added between the vertices $N_j$ and $N_{j'}$, where $j'$ is the successor of $j$, in $\mathbf{t}_{(\Theta_{\mathbf{t}})}$. The edges are drawn so that they do not cross, and in such a way that the edge $(N_j, N_{j'})$ surrounds from above the edges that start from abscissas lying between $j$ and $j'$. The set of vertices and edges thus drawn is a tree which we call $\mathcal{D}$, see Figure 5 and [20] for a proof. Up to a time inversion, the process $\mathbf{R}^{\mathbf{t}(\Theta_{\mathbf{t}})} + 1$ is the height process of $\mathcal{D}$.

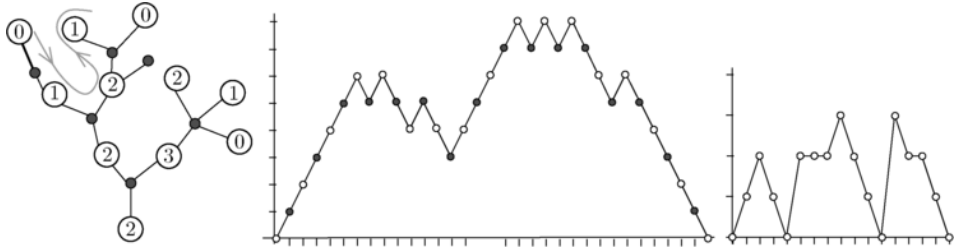

FIG. 4. *Rerooting on the first minimum of the tree $\mathbf{t}$ of Figure 3. To the right of the tree are represented the corresponding $\widehat{H}^{\mathbf{t}(\theta)}$ and $\mathbf{R}^{\mathbf{t}(\theta)}$.*

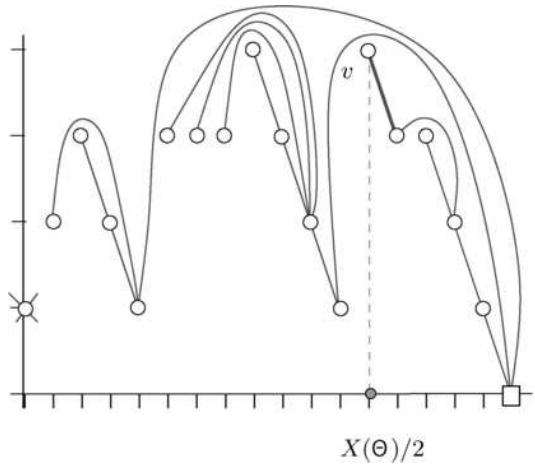

FIG. 5. *The tree $\mathcal{D}$.*

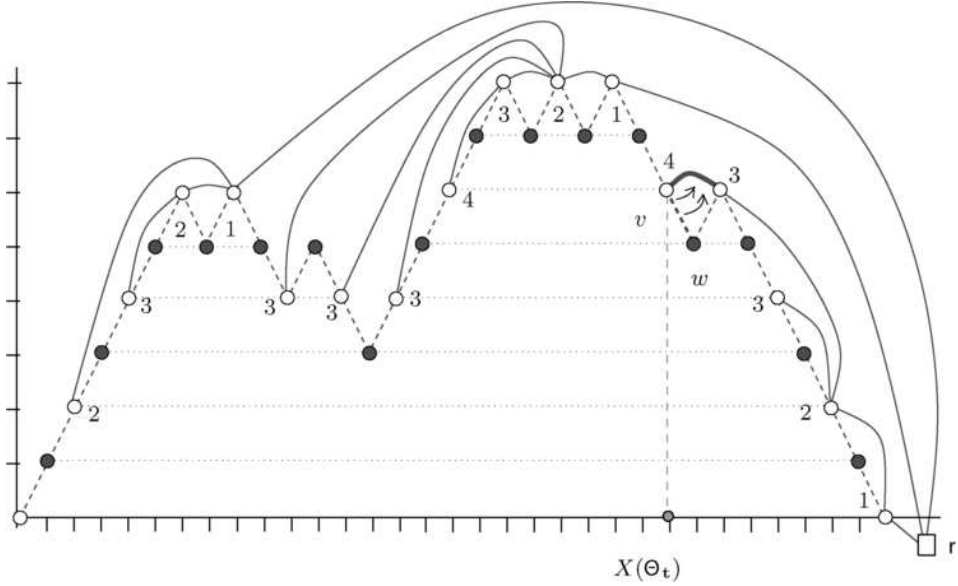

Fig. 6. *Drawing $\mathcal{D}$ on $\mathcal{G}$.*

We denote by $\mathcal{G}$ the tree $\mathbf{t}_{(\Theta_\mathbf{t})}$, whose contour process is $\widehat{H}^{\mathbf{t}(\Theta_\mathbf{t})}$. Each vertex of $\mathcal{D}$ (but the root) corresponds to a corner of a vertex in $\mathcal{G}^{(0)}$: for $j \geq 1$, the vertex $N_j$ of $\mathcal{D}$ corresponds to the vertex visited at time $2j$ for the contour order in $\mathcal{G}$. To get $\mathbf{m}$, some vertices of $\mathcal{D}$ have to be identified: glue the vertices of $\mathcal{D}$ that correspond to corners of the same vertex of $\mathcal{G}$ in such a way that the edges do not intersect. A nice way to do this is to draw $\mathcal{D}$ on the contour process of $\mathcal{G}$ as on Figure 6: place the root of $\mathcal{D}$ in the plane (not on the graph of $\widehat{H}^{\mathcal{G}}$). Then, for $0 \leq i \leq \#\mathbf{t} - 2$, place the $i+1$th vertex of $\mathcal{D}$ on the $2\#\mathbf{t} - 2i - 2$ th corner of $\mathcal{G}^{(0)}$. Then, use a deformation of the plane in order to glue together the corners of $\mathcal{G}$, corresponding to the same vertices. They are specified by horizontal dotted lines on Figure 6. The variable $X(\Theta_\mathbf{t})$ finally allows to find the root of $\mathbf{m}$.

Some changes appear when compared to [20]. Here the maps are both rooted and pointed instead of being only rooted [the variable $X(\Theta_\mathbf{t})$ allows to handle this], here, the natural traversal for both trees is the clockwise traversal, and here $\#\mathcal{D} = \#\mathcal{G}$ instead of $\#\mathcal{D} = 2\#\mathcal{G}$ for quadrangulations. Also, when conditioning by the number of faces (or vertices), the size of $\mathcal{D}$ and $\mathcal{G}$ are random, whereas they are deterministic in the case of quadrangulations. Last, only half of the vertices of $\mathcal{G}$ are used here instead of all of them.

It is then possible to adapt the notion of abstract maps introduced in [20] in order to handle these differences: a component encoding the distinguished point is added, the contour processes can be taken in the space of continuous function with compact support instead of $C[0,1]$, and the convention on the traversal order on trees can be adapted to the present setting. Apart from these technicalities, it remains to get the asymptotics of the two trees $\mathcal{D}$ and $\mathcal{G}$ under the considered distributions $P_{\mathbf{q}}^{F=n}$ or $P_{\mathbf{q}}^{S=n}$.

Let $\mathbb{H}$ be the states space of the tour of the Brownian snake (it is the states space of $(e, r)$ as defined around formula (4), see [19, 20] for more details). We recall the operation of rerooting of a normalized labeled tree (see [16, 20]) defined for any $\theta \in [0, 1]$ by

$$J^{(\theta)} : \mathbb{H} \longmapsto \mathbb{H},$$
$$(\zeta, f) \longmapsto J^{(\theta)}(\zeta, f) = (\zeta^{(\theta)}, f^{(\theta)}),$$

where for any $x \in [0, 1]$,

(66)
$$f^{(\theta)}(x) = f(\theta + x) - f(\theta),$$
$$\zeta^{(\theta)}(x) = \zeta(\theta + x) + \zeta(\theta) - 2\check{\zeta}(\theta \oplus x, \theta),$$

where the additions in the arguments are modulo 1. This may be understood as follows. Suppose $(\zeta, f)$ is the encoding of a labeled tree $(\mathbf{t}, \ell)$: $\zeta$ is the (renormalized) contour process of $\mathbf{t}$, and $f$ is the (renormalized) label process associated with $(\mathbf{t}, \ell)$. Then $(\zeta^{(\theta)}, f^{(\theta)})$ is the encoding of the labeled tree $(\mathbf{t}', \ell')$ which is obtained from $(\mathbf{t}, \ell)$ by rerooting $\mathbf{t}$ on the corner that is visited at time $\theta$, and adding $-f(\theta)$ to all labels [this fixes $\ell'(\text{root}(\mathbf{t}')) = 0$]. We are particularly interested by the rerooting on $I(f) = \inf \operatorname{Arg\,min} f$, the first minimum of the label process:

$$\Phi : \mathbb{H} \longmapsto [0, 1] \times \mathbb{H},$$
$$(\zeta, f) \longmapsto (I(f), (\zeta^+, f^+)) := (I(f), (\zeta^{(I(f))}, f^{I(f)})).$$

The application $\Phi$ is invertible. Note that it would not be without the first coordinate $I(f)$. The pair $(e^+, r^+)$ corresponds to the head of the Brownian snake $(e, r)$ under $\mathbb{N}^{(1)}$. We refer to [17] and [16] for properties of $(e^+, r^+)$ and its occurrence as a limit of conditioned spatial trees.

LEMMA 27. *Under* $\mathbb{N}^{(1)}$, $I(r)$ *is uniform on* $[0, 1]$ *and independent of* $(e^+, r^+)$.

PROOF. First, according to Lemma 16 in [20] (see also [17], Proposition 2.5), $\# \operatorname{Arg\,min} r = 1$ a.s. The law of $(e, r)$ is preserved by rerooting (see [20]) and $I(r^{(\theta)}) = I(r) - \theta \mod 1$. Then $I(r)$ is uniform in [0,1]. Now, let us check the independence. Suppose that r reaches its minimum once. For any $x \in [0, 1)$, $\Phi(e^{(x)}, r^{(x)}) = (\theta - x \mod 1, (e^+, r^+))$. Hence, in each class stable by rerooting, the positive representative $(e^+, r^+)$ is independent of $I(r)$. □

The asymptotics of the trees $\mathcal{D}$ and $\mathcal{G}$, that are sufficient to get a generalization of the convergence of rescaled bipartite maps to the Brownian maps, are given by the following proposition.

PROPOSITION 28. *Let* $\mathbf{q}$ *be a regular critical weight sequence. Under* $P_{\mathbf{q}}^{F=n}$ *(resp. $P_{\mathbf{q}}^{S=n}$), the process* $(\frac{\Theta_T}{2(\#T-1)}, \frac{\widehat{H}^{T(\Theta_T)}(2(\#T-1)\cdot)}{n^{1/2}}, \frac{\mathbf{R}^{T(\Theta_T)}((\#T-1)\cdot)}{n^{1/4}})$ *converges in distribution to*

$$\left(U, \frac{4}{\sqrt{(Z_{\mathbf{q}}-1)\rho_{\mathbf{q}}}}\mathrm{e}^+, \left(\frac{4\rho_{\mathbf{q}}}{9(Z_{\mathbf{q}}-1)}\right)^{1/4}\mathrm{r}^+\right) \qquad under\ \mathbb{N}^{(1)}$$

$$\left[resp.\ \left(U, \frac{4}{\sqrt{\rho_{\mathbf{q}}}}\mathrm{e}^+, \left(\frac{4\rho_{\mathbf{q}}}{9}\right)^{1/4}\mathrm{r}^+\right)\ under\ \mathbb{N}^{(1)}\right],$$

*where* $U$ *is an uniform random variable independent of* $(\mathrm{e}^+, \mathrm{r}^+)$, *with the constants given in Section* 3.2.

PROOF. First, under $P_{\mathbf{q}}^{F=n}$ (resp. $P_{\mathbf{q}}^{S=n}$), the process $\Phi(\frac{\widehat{H}^T(2(\#T-1)\cdot)}{n^{1/2}}, \frac{\widehat{R}^T(2(\#T-1)\cdot)}{n^{1/4}})$ converges in distribution to

$$\Phi\left(\frac{4}{\sqrt{(Z_{\mathbf{q}}-1)\rho_{\mathbf{q}}}}\mathrm{e}, \left(\frac{4\rho_{\mathbf{q}}}{9(Z_{\mathbf{q}}-1)}\right)^{1/4}\mathrm{r}\right) \qquad \text{under } \mathbb{N}^{(1)}$$

$$\left[\text{resp. } \Phi\left(\frac{4}{\sqrt{\rho_{\mathbf{q}}}}\mathrm{e}, \left(\frac{4\rho_{\mathbf{q}}}{9}\right)^{1/4}\mathrm{r}\right) \text{ under } \mathbb{N}^{(1)}\right].$$

Indeed, the applications Arg min and then $\Phi$ are continuous on the space of continuous functions that reach their minimum once, and r reaches a.s. its minimum once (see [17, 20]). The conclusion follows from Theorem 8 and Lemma 27. □

**Acknowledgments.** Grateful thanks are due to the two anonymous referees for extremely patient and thorough readings of earlier versions of this work. Their numerous comments were most appreciated and allowed to greatly improve the exposition.

## REFERENCES

[1] ALDOUS, D. J. (1991). The continuum random tree. II. An overview. In *Stochastic Analysis (Durham, 1990). London Math. Soc. Lecture Note Ser.* **167** 23–70. Cambridge Univ. Press. MR1166406
[2] ALDOUS, D. J. (1993). The continuum random tree. III. *Ann. Probab.* **21** 248–289. MR1207226

[3] AMBJØRN, J., DURHUUS, B. and JONSSON, T. (1997). *Quantum Geometry. A Statistical Field Theory Approach.* Cambridge Univ. Press. MR1465433

[4] ANGEL, O. (2003). Growth and percolation on the uniform infinite planar triangulation. *Geom. Funct. Anal.* **13** 935–974. MR2024412

[5] ANGEL, O. and SCHRAMM, O. (2003). Uniform infinite planar triangulations. *Comm. Math. Phys.* **241** 191–213. MR2013797

[6] BOUTTIER, J., DI FRANCESCO, P. and GUITTER, E. (2004). Planar maps as labeled mobiles. *Electron. J. Combin.* **11** 1–27. MR2097335

[7] CHASSAING, P. and DURHUUS, B. (2006). Local limit of labelled trees and expected volume growth in a random quadrangulation. *Ann. Probab.* **34** 879–917. MR2243873

[8] CHASSAING, P. and SCHAEFFER, G. (2004). Random planar lattices and integrated superBrownian excursion. *Probab. Theory Related Fields* **128** 161–212. MR2031225

[9] DUQUESNE, T. (2003). A limit theorem for the contour process of conditioned Galton–Watson trees. *Ann. Probab.* **31** 996–1027. MR1964956

[10] DUQUESNE, T. and LE GALL, J.-F. (2002). Random trees, Lévy processes and spatial branching processes. *Astérisque* **281**. MR1954248

[11] FELLER, W. (1971). *An Introduction to Probability Theory and Its Applications II*, 2nd ed. Wiley, New York. MR0270403

[12] GITTENBERGER, B. (2004). A note on "State spaces of the snake and its tour—convergence of the discrete snake," by J.-F. Marckert and A. Mokkadem. *J. Theoret. Probab.* **16** 1063–1067. MR2033198

[13] JANSON, S. and MARCKERT, J.-F. (2005). Convergence of discrete snakes. *J. Theoret. Probab.* **18** 615–645. MR2167644

[14] KURTZ, T., LYONS, R., PEMANTLE, R. and PERES, Y. (1997). A conceptual proof of the Kesten–Stigum theorem for multi-type branching processes. In *Classical and Modern Branching Processes* (*Minneapolis, MN, 1994*). *IMA Vol. Math. Appl.* **84** 181–185. Springer, New York. MR1601737

[15] LE GALL, J.-F. (1999). *Spatial Branching Processes, Random Snakes and Partial Differential Equations.* Birkhäuser, Basel. MR1714707

[16] LE GALL, J.-F. (2006). Conditional limit theorem for tree-indexed random walks. *Stochastic Process. Appl.* **116** 539–567. MR2205115

[17] LE GALL, J.-F. and WEILL, M. (2006). Conditioned Brownian trees. *Ann. Inst. H. Poincaré Probab. Statist.* **42** 455–489. MR2242956

[18] MARCKERT, J.-F. and MOKKADEM, A. (2003). The depth first processes of Galton–Watson trees converge to the same Brownian excursion. *Ann. Probab.* **31** 1655–1678. MR1989446

[19] MARCKERT, J.-F. and MOKKADEM, A. (2004). States spaces of the snake and its tour—convergence of the discrete snake. *J. Theoret. Probab.* **16** 1015–1046. MR2033196

[20] MARCKERT, J.-F. and MOKKADEM, A. (2006). Limit of normalized random quadrangulations: The Brownian map. *Ann. Probab.* **34** 2144–2202.

[21] PETROV, V. V. (1995). *Limit Theorems of Probability Theory.* Oxford Univ. Press, New York. MR1353441

[22] PITMAN, J. (2006). *Combinatorial Stochastic Processes. Lecture Notes in Math.* **1875**. Springer, Berlin. MR2245368

[23] SCHAEFFER, G. (1998). Conjugaison d'arbres et cartes combinatoires aléatoires. Ph.D. thesis, Univ. Bordeaux I. Available at http://www.lix.polytechnique.fr/~schaeffe/Biblio/PhD-Schaeffer.ps.

[24] Stroock, D. W. (1993). *Probability Theory, an Analytic View*. Cambridge Univ. Press. MR1267569

CNRS & LaBRI  
Université Bordeaux 1  
351 cours de la Libération  
33405 Talence Cedex  
France  
E-mail: marckert@labri.fr  
URL: www.labri.fr/perso/marckert/

CNRS & Laboratoire de Mathématique  
Université de Paris-Sud  
Bât. 425  
91405 Orsay Cedex  
France  
E-mail: Gregory.Miermont@math.u-psud.fr  
URL: www.math.u-psud.fr/~miermont